\documentclass[12pt]{amsart}
\usepackage{placeins}
\usepackage{amssymb}
\usepackage[usenames]{color}

\newtheorem{theorem}{Theorem}
\newtheorem{example}[theorem]{Example}
\newtheorem{proposition}[theorem]{Proposition}

\newtheorem{corollary}[theorem]{Corollary}
\newtheorem{remark}[theorem]{Remark}
\newtheorem{definition}{Definition}
\newtheorem{lemma}[theorem]{Lemma}
\newtheorem{problem}[theorem]{Problem}
\newcommand{\F}{{\mathcal F}}
\newcommand{\A}{{\mathcal A}}

\newcommand{\N}{{\mathcal N}}
\newcommand{\M}{{\mathfrak M}}

\begin{document}

{

\title{Equilevel algebras}

\author{V.A.~Vassiliev}
\address{Faculty of Mathematics and Computer Science, Weizmann Institute of Science,
7610001, Herzl St 234, Rehovot, Israel}
\email{vavassiliev@gmail.com}

\begin{abstract}
{\em Singular knots} are smooth maps $S^1 \to {\mathbb R}^3$ that have self-inter\-sec\-tions or points at which the derivative vanishes. 
We describe an algebraic classification of their singularities, and study the topological properties of the corresponding stratification of the {\em discriminant subset} $\Sigma \subset C^\infty (S^1, {\mathbb R}^3)$ consisting of all singular knots. This classification is defined by a system of special subalgebras of the function space $C^\infty(S^1, {\mathbb R})$. These are either defined by the {\em chord diagrams}, that is, by finite sets of conditions of the form $f(x_i) = f(\tilde x_i)$, $\{x_i, \tilde x_i\} \subset S^1$, or are the limit positions of such subalgebras in the space of all subspaces of a fixed codimension in $C^\infty(S^1, {\mathbb R})$. These limits arise at various collisions of the points $x_i, \tilde x_i$ that define the chord diagrams. 

For each natural number $k$, the set of such codimension-$k$ subalgebras is a $2k$-dimensional compact semialgebraic variety with a canonical $k$-dimensional vector bundle on it. 
We describe the natural stratification of these varieties for $k \leq 3$ and compute their cohomology rings and characteristic classes of canonical vector bundles. We also find many cohomology classes of sets of these subalgebras of arbitrary codimensions, and prove geometric corollaries concerning the topology of corresponding discriminant strata, in particular on their intersections with the finite-dimensional approximating subspaces of the knot space. 
\end{abstract}


\thanks{This work was supported by the Absorption Center in Science of the Ministry of Immigration and Absorption of the State of Israel} 

\maketitle

\section{Introduction}

\subsection{Singular knots, function algebras, and discriminant stratification}  The space $C^\infty(S^1, {\mathbb R}^3)$ consists of the {\em knots}, i.e., smooth embeddings $f: S^1 \to {\mathbb R}^3$, and {\em singular knots}, which have self-intersections or points of vanishing derivative. The set of all the singular knots is called the
{\em discriminant} and is denoted by $\Sigma$. It is a hypersurface in $C^\infty(S^1, {\mathbb R}^3)$; the isotopy classes of knots are the connected components of the complementary set $C^\infty(S^1, {\mathbb R}^3) \setminus \Sigma$. The nonsingular points of the discriminant are maps with only one transverse self-intersection point. 
More complicated discriminant points are the maps with many simple self-intersections, cusps, self-intersections of higher multiplicity, self-tangencies, etc. 

All numerical knot invariants (modulo the constants) can be realized by the linking numbers in $C^\infty(S^1, {\mathbb R}^3)$ with appropriate closed hypersurfaces contained in the discriminant. Selecting and describing such hypersurfaces that define invariants requires studying the discriminant strata of arbitrary codimensions. Moreover, this study is useful in constructing higher cohomology classes of the spaces of knots in arbitrary spaces ${\mathbb R}^n$, $n \geq 3$. 

It is convenient to speak of the discriminant strata in the terms of subalgebras of the function space $C^\infty(S^1, {\mathbb R})$. For instance, typical discriminant points are related with the subalgebras of codimension one, which are defined by single {\em chords}, i.e. unordered pairs of points in $ S^1$. The algebra corresponding to a chord
 $\{x, \tilde x\}$ consists of functions that satisfy the equation $f(x)=f(\tilde x)$. A map $S^1 \to {\mathbb R}^3$ has a self-intersection if all its three coordinate components belong to such an algebra. The set of maps satisfying this condition for a given chord $\{x, \tilde x\}$ is a subspace of codimension three in $C^\infty(S^1, {\mathbb R}^3)$. Almost all of the discriminant is swept out by the family of these subspaces parameterized by the two-dimensional configuration space of all chords.
The simplest $k$-fold self-intersections of the discriminant are described by {\em $k$-chord diagrams} (see e.g. \cite{ks}, \cite{bl}, \cite{CDM}), i.e., collections of $k$ non-intersecting pairs of points $\{x_i, \tilde x_i\}$, $i=1, \dots, k$, and by the corresponding subalgebras of codimension $k$, which are defined by $k$ conditions $f(x_i)=f(\tilde x_i)$. 

More complicated discriminant strata and the related subalgebras
occur at the degenerations of the chord diagrams. As the points of $k$-chord diagrams approach one another in $S^1$, the corresponding subalgebras in $C^\infty(S^1, {\mathbb R})$ move continuously in the space of all subspaces of codimension $k$ and converge to the algebras associated with more degenerate maps. 
For instance, when two points $x_t, \tilde x_t$ of the same chord depending on parameter $t$ converge to a point $\check x $, the algebras defined by conditions $f(x_t)=f(\tilde x_t)$ tend to the algebra defined by the condition $f'(\check x)=0$. When a point of one chord meets a point of a different chord, the codimension-two algebras defined by the pairs of chords tend to the algebra defined by a condition of the form $f(x_1)=f(x_2) = f(x_3)$. These two simplest degenerations of chord diagrams provide the so-called 1- and 3-term relations in the theory of finite-type invariants. 

More complicated collisions of chord points yield more complicated algebras. For example, the discriminant stratum consisting of singular knots with self-tangencies is defined by a family of codimension-two subalgebras that depend on two points $A \neq B \in S^1$ and a number $\alpha \neq 0$: each algebra consists of functions with $f(A)=f(B)$ and $f'(A) = \alpha f'(B)$. Also, 
there is a family of codimension-three subalgebras that can be defined by a triple of points, $A, B, $ and $C$, in $S^1$ and a point $(\alpha:\beta:\gamma) \in {\mathbb R}P^2$. Each algebra in this class consists of all functions that satisfy three independent conditions
\begin{equation}
\label{eone}
f(A)=f(B)=f(C), \qquad \alpha f'(A) + \beta f'(b) + \gamma f'(C) =0. 
\end{equation} 
These algebras are the limits of those defined by three chords when the endpoints of each chord collide with some endpoints of two other chords. Some other subalgebras of codimension three can be defined by two points $A, B \in {\mathbb R}^1$ and one point $(\alpha:\beta:\gamma) \in {\mathbb R}P^2$; they consist of all functions $f: {\mathbb R}^1 \to {\mathbb R}$ such that 
\begin{equation} 
\label{twone}
f(A) = f(B), \qquad f'(A) = 0, \qquad \alpha f'''(A) + \beta f''(A) + \gamma f'(B) =0.
\end{equation} 

{\em Equilevel algebras} of codimension $k$ in $C^\infty(S^1, {\mathbb R})$ are subalgebras that can be obtained as limits of subalgebras defined by $k$-chord diagrams
when the endpoints of these chords collide in some way; see \S~\ref{maindef} for the exact definition. The algebraic classification of these algebras provides a classification of singular knots. For example, singular knots whose singularities are defined by the algebras of type (\ref{eone}) have triple self-intersection points with linearly dependent derivatives of three local branches. 

The discriminant strata are classified according to the deepest (of greatest codimension) subalgebras containing all three coordinate functions defining their singular knots. For example, the only two strata related with algebras of codimension one but not two contain the singular knots with simple (=transverse) self-intersections or {\em simple cusps}  (i.e., points where the first derivative vanishes, but the second and the third derivatives are non-proportional). The nine strata related with algebras of codimension two  but not three  contain, respectively, the singular knots with 1) two simple self-intersections, 2) a simple self-intersection and a simple cusp, 3) a triple point with linearly independent  derivatives, 4) and 5) self-tangencies with equal or opposite directions, 6) two simple cusps, 7) a simple cusp generically intersected by a different branch of the curve, 8) a point where the first derivative vanishes, the second and the third derivatives are proportional, and the fourth derivative is linearly independent of them, 9) a point where both first and second derivatives vanish but the third does not. All these strata are parameterized by the involved points of the circle. The strata 4), 5), and 8) have an additional module: a proportionality coefficient of derivatives. For another example of the moduli of codimension-3 algebras, see the coefficients $(\alpha:\beta:\gamma) \in {\mathbb R}P^2$ in (\ref{eone}) and (\ref{twone}).

A coarser classification of knot singularities was used in \cite{ks}, but singular knots of the same class of that classification can look very different geometrically.

For any natural number $k$, we denote by $CD_k(S^1)$ the space of all codimension-$k$ subalgebras in $C^\infty(S^1, {\mathbb R})$ defined by collections of chords, and denote by $\overline{CD}_k(S^1)$ the space of all codimension-$k$ equilevel subalgebras, i.e. the closure of $CD_k(S^1)$ in the space of codimension-$k$ subalgebras. This space $\overline{CD}_k(S^1)$ is a compact $2k$-dimensional semialgebraic variety. The {\em tautological bundle} on this variety is the vector bundle whose fibers are these subalgebras. 
It is also convenient to consider the $k$-dimensional {\em normal bundle}, that is, the bundle of quotient spaces of $C^\infty(S^1, {\mathbb R})$ by the fibers of the tautological bundle.

We classify all equilevel algebras of codimension $\leq 3$ in $C^\infty(S^1, {\mathbb R})$ and compute the mod 2 and rational cohomology rings of the spaces $\overline{CD}_k(S^1)$, $k \leq 3$. We also describe many higher-dimensional homology and cohomology classes, including the characteristic classes of the normal bundles, of the spaces $\overline{CD}_k(S^1)$ for arbitrary $k$.

\begin{example} \rm
The space $CD_1(S^1)$ of single chords (i.e., ``one-chord diagrams'') on $S^1$ is homeomorphic to the open M\"obius band. 
Its closure $\overline{CD}_1(S^1)$ is homeomorphic to the closed M\"obius band. 
The added points are all subalgebras parameterized by the points $A \in S^1$ and consisting of functions that satisfy the condition $f'(A)=0$.
A cell decomposition of the space $\overline{CD}_1(S^1)$ consists of four cells, which are denoted by
\unitlength 1mm
\begin{picture}(8,6)
\put(4,2){\circle{6}}
\put(4,-0.9){\line(0,1){5.8}}
\put(7,2){\circle*{1}}
\end{picture}, 
\begin{picture}(8,6)
\put(4,2){\circle{6}}
\put(1.2,2){\line(1,0){6}}
\put(7,2){\circle*{1}}
\end{picture}, 
\begin{picture}(8,6)
\put(4,2){\circle{6}}
\put(0.3,0.8){\small$\ast$}
\put(7,2){\circle*{1}}
\end{picture}, and 
\begin{picture}(8,6)
\put(4,2){\circle{6}}
\put(7,1){\small $\ast$}
\put(7,2){\circle*{1}}
\end{picture} ,
\ of dimensions 2, 1, 1, and 0. They consist of algebras defined, respectively, by chords $(A, B)$ that do not contain the distinguished point $\bullet \in S^1$, by chords that contain this point, by the conditions $f'(A)=0$ for any $A \neq \bullet$, and by the unique such condition for $A=\bullet$. The mod 2 cell complex with these generators has the boundary maps $\partial_2 \left( \begin{picture}(8,4)
\put(4,1){\circle{6}}
\put(4,-1.9){\line(0,1){5.8}}
\put(7,1){\circle*{1}}
\end{picture} \right) = \begin{picture}(8,4)
\put(4,1){\circle{6}}
\put(0.3,-0.2){\small$\ast$}
\put(7,1){\circle*{1}}
\end{picture}$, $\partial_1 \left( \begin{picture}(8,4)
\put(4,1){\circle{6}}
\put(1.2,1){\line(1,0){6}}
\put(7,1){\circle*{1}}
\end{picture} \right) = \partial_1 \left( \begin{picture}(8,4)
\put(4,1){\circle{6}}
\put(0.3,-0.2){\small$\ast$}
\put(7,1){\circle*{1}}
\end{picture}\right) = 0$. 
Thus, $H_2\left(\overline{CD}_1(S^1), {\mathbb Z}_2\right) \simeq 0$, and the group $H_1\left(\overline{CD}_1(S^1), {\mathbb Z}_2\right) \simeq {\mathbb Z}_2$ is generated by the class of the cell \begin{picture}(8,4)
\put(4,1.5){\circle{6}}
\put(1.2,1.5){\line(1,0){6}}
\put(7,1.5){\circle*{1}}
\end{picture}. Alternatively, this group is generated by the cycle 
\begin{picture}(8,4)
\put(4,1.5){\circle{6.3}}
\put(4,-1.5){\vector(0,1){6}}
\put(4,4.5){\vector(0,-1){6}}
\put(7,1.5){\circle*{1}}
\end{picture}
parameterized by the space $S^1/{\mathbb Z}_2 \equiv {\mathbb R}P^1$ of diameters of $S^1$. Each element of this cycle is the algebra of functions that have equal values at the endpoints of a fixed diameter. The normal bundle on $\overline{CD}_1(S^1)$ is non-trivial and even not orientable.
\end{example}

The natural cell decompositions of the spaces $\overline{CD}_2(S^1)$ and $\overline{CD}_3(S^1)$ consist of 28 and 260 cells, respectively. Their cohomology rings are described in \S~\ref{toco} below.

\subsection{Finite-dimensional approximations of the knot space, characteristic classes of canonical bundles,  and non-stability of discriminant strata}

The characteristic classes of the normal bundles on $\overline{CD}_k(S^1)$ participate in the study of finite-dimensional approximations of the space of knots, and specifically of the stabilization of topological properties of intersections of discriminant strata with the approximating spaces. Let ${\mathcal L} $ be a generic vector subspace of dimension $N$ in $ C^\infty(S^1, {\mathbb R}^3)$. If $N$ is sufficiently large with respect to $k$, $N \geq 5k$, then each codimension-$3k$ space of singular knots defined by a codimension-$k$ equilevel algebra (e.g. by $k$ independent chords) intersects the space ${\mathcal L}$ along a subspace of codimension exactly $3k$. However, for slightly smaller $N$ there necessarily exist such spaces with smaller codimensions of intersection with ${\mathcal L}$. 

For a basic example, the non-orientability of the normal bundle on $\overline{CD}_1(S^1)$ implies that there are no three-dimensional subspaces in $C^\infty(S^1, {\mathbb R}^3)$ all of whose non-zero points are the knots. In other words, for each three-dimensional vector subspace ${\mathcal L} \subset C^\infty(S^1, {\mathbb R}^3)$  there exists a pair of opposite points $(A, B)$ in $S^1$ such that $f(A)=f(B)$ for some non-zero map $f \in {\mathcal L}$. Analogously, for any seven-dimensional subspace ${\mathcal L} \subset C^\infty(S^1, {\mathbb R}^3)$ there necessarily exists a pair of distinct chords such that the subspace in ${\mathcal L}$ consisting of singular knots $S^1 \to {\mathbb R}^3$ that glue the endpoints of each chord is at least two-dimensional. This follows from the non-triviality of the second Stiefel--Whitney class of the normal bundle on $CD_2(S^1)$. See Theorems \ref{corm} and \ref{corm3} below, and also \cite{pacific} for other examples.

\subsection{Order complex and discriminant resolution}
The set of all equilevel algebras of arbitrary codimensions is a partially ordered set, where the order is defined by the incidence of these algebras. Therefore, the topological order complex of the spaces $\overline{CD}_k(S^1)$ can be defined in the following way. Consider the join 
$$ \overline{CD}_1(S^1) * \overline{CD}_2(S^1) * \dots * \overline{CD}_k(S^1)$$
of these spaces (that is, the naturally topologized union of $(k-1)$-dimensional simplices, whose vertices are the points of the $k$ different spaces $\overline{CD}_i(S^1)$) and define the order complex $\lozenge_k$ as its subset consisting only of those simplices whose vertices are all incident to each other. Each complex $\lozenge_k$ is canonically embedded into $\lozenge_{k+1}$, and the total order complex $\lozenge$ is the union of all of them.
Similar to \cite{ks}, this order complex is the basis for a simplicial resolution of the discriminant set, and, consequently, for constructing the cohomology classes of the knot space. This resolution is a subspace of the product $\lozenge \times C^\infty(S^1, {\mathbb R}^3)$. For each codimension-$k$ equilevel algebra $\A$, we take
the {\em subordinate subcomplex} $\lozenge(\A)\subset \lozenge_k$ consisting of simplices all of whose vertices are the algebras containing $\A$. 
 The simplicial resolution $\sigma$ is defined as the union of all products $\lozenge(\A) \times \A^3 \subset \lozenge \times (C^\infty(S^1, {\mathbb R}))^3$. It is naturally filtered: the space $\sigma_p$ is the union of all products $\lozenge(\A) \times \A^3$ over all algebras of codimension $\leq p$. The order $p$ cohomology classes of the knot space can be defined as linking numbers with projections of closed cycles of finite codimension contained in $\sigma_p$. To correctly determine the linking index with infinite-dimensional cycles, one must use finite-dimensional approximations of the knot space and examine the stabilization of the topological properties of their intersections with the discriminant and its resolution. This raises the issues discussed in the previous subsection.

Each subspace $\sigma_p \setminus \sigma_{p-1}$ of this filtration is a vector bundle with fibers $\A^3$ and a finite-dimensional base. Thus, the characteristic classes of these bundles (which reduce to the tautological bundles on $\overline{CD}_p(S^1)$) are an essential part of their topological structure.

The same construction with $\A^3$ replaced by $\A^n$ can be applied to resolving singular knot spaces in arbitrary spaces ${\mathbb R}^n$, $n>3$. 

Almost all of this theory can be extended to the function spaces on arbitrary manifolds instead of $S^1$. This approach works in studying spaces of singular maps when the singularities are defined in terms of the simultaneous behavior of these maps at pairs of distinct points: for example, as self-intersections of the image (as in the case of knots), or as {\em non-strictly Morse functions} (i.e., functions $M \to {\mathbb R}$ having multiple critical values).

\subsection{An analog for ternary singularities} A similar theory considers subalgebras defined by collections of conditions of the form $f(x_i)=f(\tilde x_i) = f(\hat x_i),$ $i=1, \dots, k$, and by degenerations of such conditions; see arXiv: 2506.12865.
These algebras arise in the study of discriminant sets and their complements when the condition of being singular is expressed in terms of the simultaneous behavior of maps at some three points of the source manifold. A sample problem of this type is the classification of plane curves without triple self-intersections, see e.g. \cite{A}, \cite{FT}, \cite{Merx}, \cite{MD}. 

Conversely, in the more elementary theory, when the notion of a singular map is defined in the terms of its behavior at single points, the analog of the space of chord diagrams is the configuration space $B(M, k)$ of subsets of cardinality $k$ in $M$, and the analog of the space $\overline{CD}_k(M)$ is its closure in the space of all ideals of codimension $k$ in the space $C^\infty(M, {\mathbb R})$.

\medskip

Throughout the paper, the symbol $\Box_N$ indicates the end of the proof of proposition $N$, $\Box_{N(i)}$ marks the end of the proof of statement $i$ of proposition $N$.

\subsection{Topological computations}
\label{toco}

\subsubsection{Topology of the spaces $\overline{CD}_2(S^1)$ of codimension-2 equilevel algebras}

The natural classification of these algebras is given in \S~\ref{cells2}.

\begin{theorem}
\label{mthmadd}
The mod 2 cohomology groups of the space $\overline{CD}_2(S^1)$
are isomorphic to ${\mathbb Z}_2$ in dimensions 0, 1, 2, and 3, and are trivial in all other dimensions. 
\end{theorem}

The proof of this theorem and a list of basic cycles of corresponding homology groups are given in \S~\ref{cells2}. All other propositions \ref{mthmmult}--\ref{corm} of the present subsubsection are proved in \S\S~\ref{twoser}, \ref{other2} (or immediately). 
\medskip

Let $W$ be the generator of the group $H^1\left(\overline{CD}_2(S^1), {\mathbb Z}_2\right)$.

\begin{theorem}
\label{mthmmult}
In the ring $H^*\left(\overline{CD}_2(S^1), {\mathbb Z}_2\right)$,

1. $W^2 \neq 0$,

2. $W^3 = 0$. 
\end{theorem}

\begin{corollary}
\label{cor0}
The group $H^i\left(\overline{CD}_2(S^1), {\mathbb Q}\right)$ is isomorphic to ${\mathbb Q}$ for $i=0$ and $3$, and is trivial for all other $i$.
\end{corollary}

\noindent
{\it Proof.} $\overline{CD}_2(S^1)$ is a connected variety, so $H^0\left(\overline{CD}_2(S^1), {\mathbb Q}\right) \simeq {\mathbb Q}.$
The dimensions of the rational cohomology groups are no greater than those with ${\mathbb Z}_2$ coefficients. Theorem \ref{mthmadd} and statement 1 of Theorem \ref{mthmmult} implies that the group $H^1(\overline{CD}_2(S^1), {\mathbb Z})$ is a 
 torsion group. By Theorem \ref{mthmadd}, the Euler characteristic of $\overline{CD}_2(S^1)$ is 0, so its rational Betti numbers satisfy the conditions $b_1=0, b_2 \leq 1, b_3 \leq 1,$ and $ b_3-b_2=1$.
\hfill $\Box_{\ref{cor0}}$

\begin{definition} \rm
\label{defcan}
The {\it normal vector bundle} $\N_k$ on $\overline{CD}_k(S^1)$ is the $k$-dimensional vector bundle whose fiber at an equilevel subalgebra is the quotient space of $C^\infty(S^1, {\mathbb R}) $ by that subalgebra.
\end{definition}

\begin{proposition} \rm
\label{proban}
The first Stiefel--Whitney class of the normal vector bundle $\N_1$ on $\overline{CD}_1(S^1)$ is nontrivial. 
\end{proposition}

\noindent
{\it Proof.} The fiber of the tautological bundle at a point $\{A, B\} \subset S^1$ is the space of functions $f:S^1 \to {\mathbb R}$ satisfying the condition $f(A) = f(B)$.
Ordering the endpoints $(A, A+\pi)$ of an arbitrary diameter of $S^1$ defines an orientation of its normal bundle $\N_1$ at the corresponding point of $\overline{CD}_1(S^1)$: we say that the positive part of the fiber consists of the cosets of functions with $f(A+\pi) > f(A)$. Moving the point $A$ continuously by the angle $\pi$ changes this orientation. \hfill $\Box_{\ref{proban}}$

\begin{theorem}
\label{SWth}
1. The first Stiefel--Whitney class of the bundle $\N_2$ on $\overline{CD}_2(S^1)$ is nontrivial. 

2. The second Stiefel--Whitney class of the bundle $\N_2$ is also nontrivial. 
\end{theorem}

\begin{corollary}
\label{corSW}
The total Stiefel--Whitney class of the bundle $\N_2$ on $\overline{CD}_2(S^1)$ is equal to $1+W+W^2$. \hfill $\Box_{\ref{corSW}}$
\end{corollary}

These calculations imply the following Borsuk--Ulam-type statement.

\begin{proposition}
\label{BU}
For any pair of linearly independent smooth functions $f, g: S^1 \to {\mathbb R}$, there exists a nontrivial linear combination $\lambda f+ \mu g$ with real coefficients $\lambda$ and $ \mu$ whose derivative vanishes at a pair of opposite points of the circle. For a {\em generic} pair of functions \ $f$ \ and \ $g$, the number of linear combinations satisfying this condition $($and considered up to multiplication by nonzero constants$)$ is odd.
\end{proposition}

For any pair of distinct chords, the set of all maps $S^1 \to {\mathbb R}^3$ that take equal values at the endpoints of each chord is a vector subspace of codimension six in $C^\infty(S^1, {\mathbb R}^3)$. Correspondingly, the intersection of this space with a generic seven-dimensional vector subspace ${\mathcal L}^7 \subset C^\infty(S^1, {\mathbb R}^3)$ has dimension 1. 

\begin{theorem}
\label{corm}
For any seven-dimensional subspace ${\mathcal L}^7 \subset C^\infty(S^1, {\mathbb R}^3)$, there exist pairs of distinct chords in $S^1$ such that the space of maps $f \in {\mathcal L}^7$ that take equal values at the endpoints of each chord is at least two-dimensional. Moreover, the set of pairs of chords that satisfy this condition is at least two-dimensional.
\end{theorem} 

\begin{remark} \rm  Many estimates of this type are proved in \cite{pacific}. Theorem \ref{corm} improves the estimate from \cite{pacific} that states only that  for any six-dimensional subspace ${\mathcal L}^6 \subset C^\infty(S^1, {\mathbb R}^3)$ there exist pairs of distinct chords in $S^1$ such that the space of maps $f \in {\mathcal L}^6$ having equal values at the endpoints of each chord is at least one-dimensional.
\end{remark}

\subsubsection{Topology of the space $\overline{CD}_3(S^1)$}
\label{res3}

\begin{theorem}
\label{thmhom}
$H^0\left(\overline{CD}_3(S^1), {\mathbb Z}_2\right) \simeq 
H^1\left(\overline{CD}_3(S^1), {\mathbb Z}_2\right) \simeq {\mathbb Z}_2$.

$H^2\left(\overline{CD}_3(S^1), {\mathbb Z}_2\right) \simeq
H^3\left(\overline{CD}_3(S^1), {\mathbb Z}_2\right) \simeq {\mathbb Z}_2^2.$

$H_i\left(\overline{CD}_3(S^1), {\mathbb Z}_2\right) \simeq 0 $ \ for all \ $i \geq 4$.
\end{theorem}

The cycles generating the corresponding homology groups will be described in Propositions \ref{hpro3}, \ref{pro18}, and \ref{probet}.

\begin{theorem}
\label{rat3}
$H^0\left(\overline{CD}_3(S^1), {\mathbb Q} \right) \simeq 
H^3\left(\overline{CD}_3(S^1), {\mathbb Q} \right) \simeq {\mathbb Q}$; all other rational cohomology groups of the space $\overline{CD}_3(S^1)$ are trivial.
\end{theorem}

\begin{theorem}
\label{thmsw}
1. The group $H^1\left(\overline{CD}_3(S^1), {\mathbb Z}_2\right)$ is generated by the first Stiefel--Whitney class $w_1(\N_3)$ of the normal bundle. 

\noindent
2. The group $H^2\left(\overline{CD}_3(S^1), {\mathbb Z}_2\right)$ is generated by the classes $w_1^2(\N_3)$ and $w_2(\N_3)$. 

\noindent
3. The class $w_3(\N_3)$ is a nontrivial element of the group $H^3\left(\overline{CD}_3(S^1), {\mathbb Z}_2\right)$.

\noindent
4.This class $w_1^3(\N_3)$ is equal to $w_3(\N_3)$.

\noindent
5.
The class $w_1(\N_3) \smile w_2(\N_3)$ is trivial.
\end{theorem}

\begin{remark} \rm
The multiplication in the ring $H^*\left(\overline{CD}_3(S^1), {\mathbb Z}_2\right)$ is completely determined by Theorem \ref{thmsw}.
\end{remark}

\begin{corollary}
\label{cor71}
$H^1(\overline{CD}_3(S^1), {\mathbb Z}_4) \simeq {\mathbb Z}_2$.
$H^2(\overline{CD}_3(S^1), {\mathbb Z}_4) \simeq {\mathbb Z}_2 \oplus {\mathbb Z}_2$.

\noindent
$H^3(\overline{CD}_3(S^1), {\mathbb Z}_4) \simeq {\mathbb Z}_2 \oplus {\mathbb Z}_4.$ 
\end{corollary}

\noindent
{\it Proof.} Consider the exact sequence of the cohomology groups of the space $\overline{CD}_3(S^1)$ defined by the short
 exact sequence of coefficients $0 \to {\mathbb Z}_2 \to {\mathbb Z}_4 \to {\mathbb Z}_2 \to 0$. Its differentials $H^i\left(\overline{CD}_3(S^1), {\mathbb Z}_2\right) \to H^{i+1}\left(\overline{CD}_3(S^1), {\mathbb Z}_2\right)$ are the Bockstein operators $Sq^1$.
For any vector bundle,
$Sq^1(w_1^2)=0$ and $Sq^1(w_2) = w_3 + w_1 \smile w_2$ (see e.g. \cite{borel}). By statements 3 and 5 of Theorem \ref{thmsw}, $w_3(\N_3)+w_1(\N_3)\smile w_2(\N_3) \neq 0$, so the map $Sq^1: H^2\left(\overline{CD}_3(S^1), {\mathbb Z}_2\right) \to H^{3}\left(\overline{CD}_3(S^1), {\mathbb Z}_2\right)$ is of rank 1. Statements 1 and 2 of the same theorem imply that the map $Sq^1: H^1(\overline{CD}_3(S^1), {\mathbb Z}_2) \to H^{2}(\overline{CD}_3(S^1), {\mathbb Z}_2)$ is also nontrivial. 
Thus, $H^1\left(\overline{CD}_3(S^1), {\mathbb Z}_4\right) \simeq {\mathbb Z}_2$, the group $H^2\left(\overline{CD}_3(S^1), {\mathbb Z}_4\right)$ consists of four elements, and the group $H^3\left(\overline{CD}_3(S^1), {\mathbb Z}_4\right)$ consists of eight elements. Also, the Euler characteristic of $\overline{CD}_3(S^1)$ is equal to 0 and $H^0\left(\overline{CD}_3(S^1), {\mathbb Z}_4\right) \simeq {\mathbb Z}_4$. Therefore, the Betti numbers of $\overline{CD}_3(S^1)$ with ${\mathbb Z}_4$ coefficients in dimensions $0, 1, 2, 3$ can only be $1, 0, 0, 1$.
\hfill $\Box_{\ref{cor71}}$
\medskip

Theorem \ref{rat3} follows immediately from this corollary and Euler characteristic considerations. \hfill $\Box_{\ref{rat3}}$
\medskip

Theorems \ref{thmhom} and \ref{thmsw} are proved in \S~\ref{homapp} using the CW-structure of the space $\overline{CD}_3(S^1)$ described in \S\S~\ref{cells}--\ref{ic}.

\subsubsection{Lower estimates of the groups $H^*\left(\overline{CD}_k(S^1), {\mathbb Z}_2\right)$ for arbitrary $k$}

 For each topological space $X$, the {\em unordered $k$-configuration space} $B(X,k)$ is the space of all subsets of cardinality $k$ in $X$, supplied with the natural topology.
The rings $H^*(B({\mathbb R}^2, k), {\mathbb Z}_2)$ for all $k$ have been computed in \cite{fuks}. In particular, it was shown there that the groups $H^i(B({\mathbb R}^2, k), {\mathbb Z}_2)$
are non-trivial for all $i \leq k - I_2(k)$, where $I_2(k)$ is the number of ones in the binary decomposition of $k$. The cohomology ring $H^*\left(\overline{CD}_k(S^1), {\mathbb Z}_2\right)$ is ``not smaller'' than $H^*(B({\mathbb R}^2, k), {\mathbb Z}_2)$
in the following exact sense.

\begin{theorem}
\label{estim}
The ring $H^*(B({\mathbb R}^2, k), {\mathbb Z}_2)$ is isomorphic to a quotient ring of the subring in $H^*\left(\overline{CD}_k(S^1), {\mathbb Z}_2\right)$ generated by the Stiefel-Whitney classes of the normal bundle.
\end{theorem}

This theorem implies that the spaces $\overline{CD}_k(S^1)$ have nontrivial cohomology classes in all dimensions up to $k-I_2(k)$.
The following Theorem \ref{estim2} improves this estimate to the number $k$. However, Theorem \ref{mthmadd} shows that even this estimate is not sharp for $k=2$.

\begin{theorem}
\label{estim2}
For any natural $k$, 

1$)$ the class $w_k(\N_k) \in H^k(\overline{CD}_k(S^1), {\mathbb Z}_2)$ is nontrivial;

2$)$ all classes $w_1^i(\N_k) \in H^i(\overline{CD}_k(S^1), {\mathbb Z}_2),$ $i \leq k-1$, are nontrivial

3$)$ the class $ w_{k-1} (\N_k)$ is nontrivial, moreover it is different from the class
$w_1^{k-1}(\N_k)$ for all $k>2$.
\end{theorem}

\begin{remark} \rm
The second statement of Theorem \ref{SWth} and the third statement of Theorem \ref{thmsw} are special cases of the first statement of Theorem \ref{estim2}.
\end{remark}

\begin{theorem}
\label{estim3}
If $k$ is an even number greater than 2, then the group $H^k(\overline{CD}_k, {\mathbb Z}_2)$ is at least two-dimensional.
\end{theorem}

\begin{theorem}
\label{corm3}
For any vector subspace ${\mathcal L}^N \subset C^\infty(S^1, {\mathbb R})$ of dimension $N
\in [k, 2k-2]$, there exists an equilevel subalgebra of codimension $k$ in $C^\infty(S^1, {\mathbb R})$such that the dimension of its intersection with ${\mathcal L}^N$ is greater than $N-k$.
\end{theorem}

 Theorem \ref{estim} is proved in \S~\ref{prestim}. Theorems \ref{estim2}, \ref{estim3}, and \ref{corm3} are proved in \S~\ref{twoser}.

\section{Formal definition and basic properties of equilevel algebras}
\subsection{Main definition}
\label{maindef}

Let $M$ be a non-singular real algebraic manifold. 
Consider the {\em ordered configuration space} $F(M, 2k) \subset M^{2k},$ i.e., the set of all sequences 
\begin{equation}
\label{ocd}
(x_1, \tilde x_1; x_2, \tilde x_2; \dots; x_k, \tilde x_k)
\end{equation}
 of $2k$ points of $M$ whose elements are all distinct. Each point of this space defines a subalgebra of codimension $k$ in $C^\infty(M, {\mathbb R})$ consisting of functions that satisfy all $k$ equations $f(x_i)=f(\tilde x_i)$.

Consider an arbitrary germ of an algebraic parametric curve $\nu: [0,\varepsilon) \to M^{2k}$ such that all points $\nu(t),$ $t \in (0, \varepsilon),$ belong to the subspace $F(M, 2k)$. For any $t \in (0,\varepsilon)$, denote by $\F_t$ the subalgebra defined as in the previous paragraph by the ordered configuration $\nu(t)$. 

Define $\F_0$ as the space of all functions $f \in C^\infty(M, {\mathbb R})$ such that there exists a family of smooth functions $f_t \in C^\infty( M , {\mathbb R})$, $ t \in [0,\varepsilon)$, depending continuously in the $C^{2k}(M, {\mathbb R})$-topology on the parameter $t$, and such that $f_0 =f$ and $f_t \in \F_t$ for each $t \in (0,\varepsilon). $

 Clearly, the space $\F_0$ is also a subalgebra in $C^\infty(M, {\mathbb R}).$ 

\begin{definition} \rm
Codimension-$k$ {\em Equilevel algebras}  are arbitrary subalgebras $\F_0 \subset C^\infty(M, {\mathbb R})$ that can be defined by algebraic curves $\nu$ in this way. The set of all such algebras is denoted by $\overline{CD}_k(M)$.
\end{definition}

\begin{remark} \rm
In general, such a limit subalgebra $\F_0$ does not uniquely define the point $\nu(0) \in M^{2k}$ even up to the reorderings of its elements $x_i$ and $ \tilde x_j$. For example, the subalgebra of class $\overline{CD}_2(S^1)$ consisting of all functions that satisfy the conditions $f(A)=f(B)=f(C)$ for three distinct points $A, B, C \in S^1$ can be defined by this construction with the point $\nu(0) \in (S^1)^4$ equal to any one of the points $(A, B; A, C)$, $(A, B; B, C)$, $(A, C; B, C)$, etc. 
\end{remark}

 However, the set of all points of $M$ that are involved in the sequence $\nu(0)$ for any such definition of a given equilevel subalgebra is the same for all curves $\nu$ defining this subalgebra. This set is called the {\em support} of that subalgebra. Clearly, the support of an equilevel algebra of codimension $k$ consists of at most $2k$ points.

The following is an important intermediate class of algebras between $\overline{CD}_k(M)$ and the class of algebras defined by chord diagrams (\ref{ocd}).

\begin{definition}
\label{defcd}
The subalgebras of class ${CD}_k(M)$ in $C^\infty(M, {\mathbb R})$ are the subalgebras of codimension $k$ which can be defined by several conditions $f(x_i) = f(\tilde x_i)$, $x_i \neq \tilde x_i \in M$.
\end{definition}

\begin{remark} \rm In this definition, it is not required that all $2k$ points $x_i$, $\tilde x_j$ of these conditions are different. For example, the algebra defined by two conditions $f(A)=f(B)=f(C)$ also belongs to $CD_2(M)$.)
\end{remark}

For each point $A \in M$, denote by $\M_A$ the ideal in $C^\infty(M, {\mathbb R})$ consisting of all functions vanishing at $A$.

\begin{proposition}
\label{trui}
Suppose that an equilevel subalgebra ${\mathcal F} \subset C^\infty(M, {\mathbb R})$ of class $\overline{CD}_k(M)$ has the above definition  via a curve $\nu:[0, \varepsilon) \to M^{2k}$, and its support consists of the points $A_1, \dots, A_r$, where each point $A_i$ participates in the sequence $\nu(0)$ with some multiplicity $m_i$, \ $m_1 + \dots + m_r =2k$. Then, the subalgebra ${\mathcal F}$ contains 

a$)$ all constant functions,

b$)$ the ideal 
\begin{equation} \M_{A_1}^{m_1} \cap \dots \cap \M_{A_r}^{m_r}
\label{ide}
\end{equation} in $C^\infty(M, {\mathbb R})$, that is, the set of all functions vanishing at each point $A_i$ together with all partial derivatives up to degree $m_i-1$.
\end{proposition}

\noindent
{\it Proof.} a) All algebras ${\mathcal F}_t$ defined by configurations $\nu(t),$ $t>0$, contain the space of constant functions.

b) Let $f \in C^\infty(M, {\mathbb R})$ be an arbitrary function from the ideal (\ref{ide}). We must show that it is approximated by a continuous family of functions from the algebras defined by the configurations $\nu(t)$, $t >0$, as $t$ tends to $+0$. Using additivity and the partition of unity, it is sufficient to consider only the case when $f$ is identically zero everywhere outside a small neighborhood of one point
$A_i$ of the support of $\F$, and in the local coordinates $z_1, \dots, z_{\dim M}$ in this neighborhood this function is divisible by a monomial of degree $m_i$. Thus, $f$ has the form
\begin{equation}
f_0 \cdot z_{p_1} \cdot z_{p_2}\cdot \ldots \cdot z_{p_{m_i}}, \label{onne} \end{equation}
where $f_0$ is a smooth function that is identically zero outside the coordinate neighborhood of $A_i$, and $p_1, \dots, p_{m_i}$ are arbitrary indices ranging from 1 to $\dim M$. For any $t \in (0, \varepsilon)$, let $\tilde A_{i, 1}(t), \dots , \tilde A_{i, {m_i}}(t)$ be $m_i$ points of the configuration $\nu(t)$ tending continuously to the point $A_i$.
The desired family of functions $f_t: M \to {\mathbb R}$ that tends to the function (\ref{onne}) is then given by
$$f_0 \cdot (z_{p_1} -z_{p_1}(\tilde A_{i, 1}(t))) \cdot (z_{p_2} -z_{p_2}(\tilde A_{i, 2}(t))) \cdot \ldots \cdot (z_{p_{m_i}} - z_{p_{m_i}} (\tilde A_{i, {m_i}}(t))). 
\eqno{\mbox{$\Box_{\ref{trui}}$}}
$$

\subsection{Fighting the infinite-dimensional issues}
\label{findim}

In this subsection, we show that every space $\overline{CD}_k(M)$ can be realized as a semialgebraic subvariety of a finite-dimensional nonsingular algebraic manifold.

\begin{definition} \rm
A vector subspace $ \F \subset C^\infty(M, {\mathbb R})$ is {\em of codimension $k$} if, for any vector subspace ${\mathcal L}^n \subset C^\infty(M, {\mathbb R})$ of finite dimension $n \geq k$, the intersection $\F \cap {\mathcal L}^n$ is at least $(n-k)$-dimensional, and there exist $k$-dimensional vector subspaces ${\mathcal L}^k$ intersecting $\F$ at the origin only. The subspace ${\mathcal L}^n,$ $ n \geq k$, is {\it transversal} to a subspace of codimension $k$ if the dimension of their intersection is exactly $n-k$. 
\end{definition}

It is easy to see that if the finite-dimensional subspace ${\mathcal L}$ is transversal to a subspace ${\mathcal F}$ of finite codimension, then ${\mathcal L}$ is also transversal to any larger subspace ${\mathcal F}' \supset {\mathcal F}$. 

\begin{definition} \rm 
A finite-dimensional vector subspace ${\mathcal L} \subset C^\infty(M, {\mathbb R})$ is called {\it $k$-representative} if
\begin{equation} 
\label{condition}
\mbox{
\fbox{
\parbox{11cm}{
for any $2k$ different points of $M$ and any collection of $(2k-1)$-jets of functions $M \to {\mathbb R}$ at these points, there exists a function $f \in {\mathcal L}$ whose $(2k-1)$-jets at these $2k$ points coincide with the given ones.}
}
}
\end{equation}
\end{definition}

\begin{lemma} 
\label{glemma}
For any algebraic manifold $M$ and any natural number $k$, the $k$-representative subspaces exist. Moreover, for any such $M$ 
and $k$, there exists a number $d(k, M)$ such that for any $d \geq d(k,M)$ the collections of $d$ smooth functions $M \to {\mathbb R}$, whose linear combinations form $k$-representative subspaces, are dense in the $C^\infty$-topology in the space of all collections of $d$ smooth functions.
\end{lemma}

\noindent
{\it Proof.} This follows immediately from the Thom's multijet transversality theorem, see \cite{GG}. \hfill $\Box_{\ref{glemma}}$

\begin{lemma}
\label{le24}
Let ${\mathcal L}$ be a finite-dimensional $k$-representative subspace in $C^\infty(M, {\mathbb R})$, then

1. Every equilevel algebra of class $\overline{CD}_k(M)$ defined as in \S~\ref{maindef} by an algebraic curve $\nu: [0, \varepsilon) \to M^{2k}$ is a subspace of codimension $k$ in $C^\infty(M, {\mathbb R})$.

2. Each such equilevel algebra intersects the subspace ${\mathcal L}$ along a subspace of codimension $k$ in ${\mathcal L}$.

3. The intersections of different codimension $k$ equilevel algebras with the subspace ${\mathcal L}$ are different subspaces of codimension $k$ in ${\mathcal L}$.
\end{lemma}

\noindent
{\it Proof.} Every configuration $\nu(t)$, $t>0$, arising in the above definition of the equilevel algebra, defines a point of the Grassmann manifold $G_{\dim {\mathcal L} -k}({\mathcal L})$, namely the intersection of ${\mathcal L}$ and the corresponding subalgebra $\F_t$. These points corresponding to different positive values of $t$ form a semialgebraic curve in this compact algebraic manifold and thus have a limit in it as $t \to +0$. This limit point is a subspace of codimension $k$ in ${\mathcal L}$. 

All points of this limit subspace belong to the equilevel subalgebra $\F_0$ defined by the curve $\nu$. Indeed, for any non-zero point $f$ of this subspace consider the $k$-dimensional transversal slice of this subspace in ${\mathcal L}$. For each $t \in [0, \varepsilon)$ define the function $f_t $ as the intersection point of this slice with the subspace $\F_t$. On the other hand, the points of ${\mathcal L}$ not in the limit subspace do not belong to $\F_0$, so this subspace is the intersection $\F_0 \cap {\mathcal L}$. Due to property (\ref{condition}), the space ${\mathcal L}$ is transversal to the ideal (\ref{ide}) defined by the point $\nu(0)$, and thus also to the subalgebra $\F_0$ containing this ideal. Therefore, the codimension of the subalgebra $\F_0$ in $C^\infty(M, {\mathbb R})$ is also equal to $k$. 

Equilevel algebras with different supports obviously have different intersections with ${\mathcal L}$. All equilevel algebras with the same support $(A_1, \dots, A_r)$ are completely defined by the quotient spaces of their factorization by the ideal $$\M_{A_1}^{2k} \cap \dots \cap \M_{A_r}^{2k}; $$
these quotient spaces are determined by the intersections of these algebras with the space ${\mathcal L}$. \hfill $\Box_{\ref{le24}}$
\medskip

Thus, the entire space $\overline{CD}_k(M)$ can be identified with a semialgebraic subvariety of a finite-dimensional Grassmann manifold: each subalgebra $\F$ is associated with its intersection $ \F \cap {\mathcal L}$ with a $k$-representative space ${\mathcal L}$. According to the {\em curve selection lemma} (see \cite{M}, \S~3), if $M$ is compact, then this subvariety is also compact. 
The topological structure of the space $\overline{CD}_k(M)$ is induced from the topology of the Grassmann manifold and is independent on the choice of the space ${\mathcal L}$. Indeed, each larger space ${\mathcal L}' \supset {\mathcal L}$ defines the same topology of $\overline{CD}_k(M)$, and any two finite-dimensional subspaces of $C^\infty(M,{\mathbb R})$ are subspaces of a common larger finite-dimensional subspace.

\subsection{Two problems} Is it true that the intersection of two equilevel algebras is also an equilevel algebra (at least for one-dimensional basic manifold $M$)?

More generally, is it true that any codimension-$k$ subalgebra in $C^\infty(S^1, {\mathbb R})$ defined in the terms of the behavior of its functions at some $r$ points and satisfying the conclusions of Proposition \ref{trui} is an equilevel algebra?  If not, then possibly it would be useful to consider also the theory where the spaces $\overline{CD}_k(S^1)$ are replaced with the spaces of all such algebras.

\subsection{Proof of Theorem \ref{estim}}
\label{prestim}
Let us consider the circle $S^1$ as the line ${\mathbb R}^1$ with one added point.
Denote by
${\mathcal B}_k$ the set of all unordered collections of $k$ points $(x, y) \in {\mathbb R}^2 \subset (S^1)^2,$ $x< y$, such that the system of corresponding $k$ conditions of the form $f(x)=f(y)$ defines a subalgebra of codimension $k$ in $C^\infty({S}^1, {\mathbb R})$. Each such collection is a $k$-configuration in the half-plane $\{(x, y): x < y \} \subset {\mathbb R}^2$, so it is a point of $B({\mathbb R}^2, k)$. Let $\psi$ be the map ${\mathcal B}_k \to CD_k({S}^1) \subset \overline{CD}_k(S^1)$ that maps any such configuration to the subalgebra defined by it. 

Recall that the {\em regular vector bundle} on $B({\mathbb R}^2,k)$ (in particular on ${\mathcal B}_k$) is the $k$-dimensional vector bundle, whose fiber over a $k$-point configuration is the space of ${\mathbb R}$-valued functions on that configuration. 

Denote by ${\mathfrak R}$ the subring of the ring $H^*({\mathcal B}_k, {\mathbb Z}_2)$ multiplicatively generated by the Stiefel--Whitney classes of the regular bundle. Consider the maps \begin{equation}
H^*(B({\mathbb R}^2,k), {\mathbb Z}_2) \to H^*({\mathcal B}_k, {\mathbb Z}_2) \stackrel{\psi^*}\gets H^*(\overline{CD}_k(S^1), {\mathbb Z}_2)
\label{incl}
\end{equation} defined respectively by the inclusion ${\mathcal B}_k \to B({\mathbb R}^2, k)$ and by the map $\psi$. Theorem \ref{estim} follows immediately from the following lemma.

\begin{lemma}
\label{estilem}
1. The left map $($\ref{incl}$)$ is monomorphic.

2. The ring ${\mathfrak R}$ coincides with image of the ring 
$H^*(B({\mathbb R}^2,k), {\mathbb Z}_2) $ under this map.

3. The ring ${\mathfrak R}$ belongs to the image of the subring of $H^*(\overline{CD}_k(S^1), {\mathbb Z}_2)$ generated by Stiefel--Whitney classes of the normal bundle under the right map $($\ref{incl}$)$.
\end{lemma}

\noindent
{\it Proof.} 1. Statement 1 is Lemma 16 of \cite{pacific}.

2. In \cite{fuks}, it is proven that the entire ring $H^*(B({\mathbb R}^2,k), {\mathbb Z}_2)$ is multiplicatively generated by Stiefel--Whitney classes of the regular bundle. Hence, the subring of $H^*({\mathcal B}_k, {\mathbb Z}_2)$
generated by these classes coincides with the image of the map (\ref{incl}). 

3. By Lemma 13 of \cite{pacific}, the regular bundle on ${\mathcal B}_k$ is isomorphic to the bundle induced by the map $\psi$ from the normal bundle $\N_k$ on the space $CD_k(S^1) \subset \overline{CD}_k(S^1)$. Therefore, statement 3 follows from the naturality of Stiefel--Whitney classes. \hfill $\Box_{\ref{estilem}}$ $\Box_{\ref{estim}}$ 

\section{Classification of equilevel algebras of codimension two in $C^\infty(S^1, {\mathbb R})$. Cell structure on $\overline{CD}_2(S^1)$.}
\label{cells2}

\subsection{} Consider the circle $S^1$ as the quotient space of the segment $[0, 2\pi]$ bu the relation $\{0\}=\{2\pi\}$. 

We divide the equilevel subalgebras of $C^\infty(S^1, {\mathbb R})$ into two types: those algebras
whose supports do not contain or do contain the distinguished point $\bullet = \{0, 2\pi\}\in S^1$. We will describe 14 families of such algebras of codimension two of either of these two types. These families are all diffeomorphic to open balls of certain dimensions and define the structure of a $CW$-complex on the space $\overline{CD}_2(S^1)$. Algebras that can be obtained from each other by orientation-preserving diffeomorphisms $(S^1, \bullet) \to (S^1, \bullet)$ belong to the same cell. Thus, the points of their supports that differ from $\bullet$ are among the parameters of these cells. Moreover, uniform (in a sense) algebras with the same support also lie in the same cell, and their moduli make the remaining part of the parameters of the cells.

Each cell of the first type (consisting of algebras, whose supports do not contain the point $\bullet$) is related to a cell of the second type: the elements of the latter cell are the algebras obtained from the former algebras by shifting the argument space $S^1$ along its orientation until the support point with the largest value of the coordinate comes to the distinguished point $\bullet$. The dimension of the resulting cell is one less than that of the initial cell. 

\subsection{About the notation}
\label{onnot}

In the figures illustrating the subalgebras of codimension 2, the chords connecting two points, $A$ and $B$, of the circle $S^1$ denote the condition $f(A)=f(B)$. A tripod with feet at the points $A$, $B$, and $C$ denotes the condition $f(A)=f(B)=f(C)$. An asterisk $\ast$ at point $A$ denotes the condition $f'(A)=0$, and a double asterisk denotes the condition $f'(A)=f''(A)=0$. A circled asterisk $\circledast$ at point $A$ denotes the family of algebras depending on a parameter $\alpha \in {\mathbb R}$ and defined by the pair of conditions $f'(A)=0, f'''(A) = \alpha f''(A)$. 

The double chord connecting points $A$ and $B$ denotes a family of pairs of conditions $f(A)=f(B)$, $f'(A)=\alpha f'(B)$, depending on a real parameter $\alpha \neq 0$. This notation is always accompanied by a plus or minus sign indicating the sign of $\alpha$ for all algebras forming the cell. Each algebra defined in this way by the data $A,$ $ B,$ and $\alpha$ is called $\gimel(\alpha; A, B)$. Clearly, $\gimel (\alpha; A, B) \equiv \gimel(\alpha^{-1}; B, A)$. These algebras with $\alpha $ equal to $0$ or $\infty$ are also defined: the algebra $\gimel(0; A, B)$ or $\gimel(\infty; A, B)$ consists of functions that satisfy the conditions $f(A)=f(B)$ and $f'(A)=0$ or $f'(B)=0$. However, these exceptional algebras belong to different cells, which are denoted by a single chord, one of whose endpoints is marked by an asterisk.

Each of the pictures of Theorem~\ref{mainprop} denotes a family of algebras (or of families of algebras with parameter $\alpha$) parameterized by points of their supports (that is, ends of chords and points marked by asterisks) not coinciding with the distinguished point $\bullet$.

\unitlength 0.7mm
\begin{theorem}
\label{mainprop}
The variety $\overline{CD}_2(S^1)$ has the structure of a CW-complex with:

\begin{itemize}
\item
three 4-dimensional cells \ 
\begin{picture}(12,8)
\put(6,3){\circle{10}} 
\put(10.9,3){\circle*{1}}
\put(2.4,-0.5){\line(6,1){8.2}}
\put(2.1,6.2){\line(5,-1){8.7}}
\end{picture} \ 
\begin{picture}(12,8)
\put(6,3){\circle{10}} 
\put(10.9,3){\circle*{1}}
\put(3,-1.1){\line(-1,6){1.1}}
\put(8,-1.7){\line(1,4){1.9}}
\end{picture} \ 
\begin{picture}(12,8)
\put(6,3){\circle{10}} 
\put(10.9,3){\circle*{1}}
\bezier{80}(3,-1.1)(4.5,1)(5.4,2.2)
\bezier{80}(9.0,7.1)(7.5,5)(6.3,3.4)
\put(8.85,-1.1){\line(-3,4){6}}
\end{picture} \ parameterized by the endpoints of the chords,
\item nine 3-dimensional cells: six cells \
\begin{picture}(12,9)
\put(6,3){\circle{10}} 
\put(10.9,3){\circle*{1}}
\put(5,6.6){\small $\ast$}
\put(2,0){\line(6,1){8.7}}
\end{picture} 
\ \begin{picture}(12,8)
\put(6,3){\circle{10}} 
\put(10.9,3){\circle*{1}}
\put(5,-3){\small $\ast$}
\put(1.8,6){\line(6,-1){9}}
\end{picture} 
\ \begin{picture}(12,8)
\put(6,3){\circle{10}} 
\put(10.9,3){\circle*{1}}
\put(-0.2,2){\small $\ast$}
\put(3,-1){\line(3,4){6}}
\end{picture} 
\ \begin{picture}(12,8)
\put(6,3){\circle{10}} 
\put(10.9,3){\circle*{1}}
\bezier{120}(1,3)(4,3)(6,3)
\bezier{120}(6,3)(7.2,1)(8.6,-1.3)
\bezier{120}(6,3)(7.2,5)(8.6,7.3)
\end{picture} 
\ \begin{picture}(12,8)
\put(6,3){\circle{10}} 
\put(10.9,3){\circle*{1}}
\put(2.85,-0.95){\line(3,4){6.1}}
\put(3.3,-1.3){\line(3,4){6.1}}
\put(2.8,3.2){\tiny $+$}
\end{picture} 
\ \begin{picture}(12,8)
\put(6,3){\circle{10}} 
\put(10.9,3){\circle*{1}}
\put(2.85,-0.95){\line(3,4){6.1}}
\put(3.3,-1.3){\line(3,4){6.1}}
\put(2.8,3.2){\tiny $-$}
\end{picture} 
of the first type and three cells 
\ \begin{picture}(12,12)
\put(6,3){\circle{10}} 
\put(10.9,3){\circle*{1}}
\put(11,3){\line(-4,-3){6.4}}
\put(1,3){\line(4,3){6.5}}
\end{picture} 
\ \begin{picture}(12,8)
\put(6,3){\circle{10}} 
\put(10.9,3){\circle*{1}}
\put(11,3){\line(-4,3){6.4}}
\put(1,3){\line(4,-3){6.5}}
\end{picture} 
\ \begin{picture}(12,8)
\put(6,3){\circle{10}} 
\put(10.9,3){\circle*{1}}
\bezier{80}(11,3)(8.7,2.5)(6.4,2)
\bezier{80}(5.6,1.9)(3.65,1.5)(1.3,1)
\put(6,-2.1){\line(0,1){10.3}}
\end{picture} \ of the second type,
\item ten 2-dimensional cells: four cells \ 
\begin{picture}(12,11)
\put(6,3){\circle{10}} 
\put(10.9,3){\circle*{1}}
\put(2,-2.2){\small $\ast$}
\put(5,6.6){\small $\ast$}
\end{picture} \ 
\begin{picture}(12,8)
\put(6,3){\circle{10}} 
\put(10.9,3){\circle*{1}}
\put(5,6.6){\small $\ast$}
\put(6.2,7.8){\line(1,-6){1.55}}
\end{picture} \ 
\begin{picture}(12,8)
\put(6,3){\circle{10}} 
\put(10.9,3){\circle*{1}}
\put(5,-3){\small $\ast$}
\put(6.2,-1.8){\line(1,4){2.3}}
\end{picture} \ 
\begin{picture}(12,8)
\put(6,3){\circle{10}} 
\put(10.9,3){\circle*{1}}
\put(1.5,5.7){\small $\circledast$}
\end{picture} \ of the first type and six cells
\ \begin{picture}(11,11)
\put(5,3){\circle{10}} 
\put(10,3){\circle*{1}}
\put(2.9,6.1){\small $\ast$}
\put(10,3){\line(-4,-3){6.4}}
\end{picture} \ 
\begin{picture}(11,8)
\put(5,3){\circle{10}} 
\put(10,3){\circle*{1}}
\put(5,-2.8){\small $\ast$}
\put(10,3){\line(-4,3){6.4}}
\end{picture} \ 
\begin{picture}(11,8)
\put(5,3){\circle{10}} 
\put(10,3){\circle*{1}}
\put(9.7,1.9){\small $\ast$}
\put(5,-2.1){\line(1,5){1.9}}
\end{picture} \ 
\begin{picture}(12,8)
\put(6,3){\circle{10}} 
\put(11,3){\circle*{1}}
\bezier{100}(11,3)(9,3)(6,3)
\bezier{100}(6,3)(4.5,5)(3,7)
\bezier{100}(6,3)(4.5,1)(3,-1)
\end{picture} \ 
\begin{picture}(12,8)
\put(6,3){\circle{10}} 
\put(11,3){\circle*{1}}
\put(11,3.3){\line(-1,0){10.1}}
\put(11,2.7){\line(-1,0){10.1}}
\put(4.6,4){\tiny $+$}
\end{picture} \ 
\begin{picture}(12,8)
\put(6,3){\circle{10}} \put(11,3){\circle*{1}}
\put(11,3.3){\line(-1,0){10.1}}
\put(11,2.7){\line(-1,0){10.1}}
\put(4.6,4){\tiny $-$}
\end{picture} \ of the second type,
\item five 1-dimensional cells: \ \begin{picture}(12,12)
\put(6,3){\circle{10}} 
\put(10.9,3){\circle*{1}}
\put(-0.5,2){\small $\ast$}
\put(1,2){\small $\ast$}
\end{picture} \ and \ 
\begin{picture}(12,8)
\put(6,3){\circle{10}} 
\put(10.9,3){\circle*{1}}
\put(10.9,2){\small $\ast$}
\put(0.1,3){\small $\ast$}
\end{picture} \ 
\begin{picture}(12,8)
\put(6,3){\circle{10}} 
\put(10.9,3){\circle*{1}}
\put(0.2,0.2){\small $\ast$}
\put(11,3){\line(-6,-1){10}}
\end{picture} \ 
\begin{picture}(12,8)
\put(6,3){\circle{10}} 
\put(10.9,3){\circle*{1}}
\put(11,2){\small $\ast$}
\put(11,3){\line(-5,1){9.7}}
\end{picture} \ 
\begin{picture}(12,8)
\put(6,3){\circle{10}} 
\put(10.9,3){\circle*{1}}
\put(9.7,2.2){\small $\circledast$}
\end{picture} \ ,
\item and one 0-dimensional cell \begin{picture}(12,11)
\put(6,3){\circle{10}} \put(10.8,3){\circle*{1}}
\put(9,2){\small $\ast$}
\put(11,2){\small $\ast$}
\end{picture}~. 
\end{itemize}
\medskip

The boundary operators of this cell complex mod 2 are as follows:
\unitlength 0.7 mm

$\partial \Bigl(\begin{picture}(12,8)
\put(6,1){\circle{9.5}} 
\put(10.9,1){\circle*{1}}
\put(2.6,-2.5){\line(6,1){7.7}}
\put(2.4,4.2){\line(5,-1){8.1}}
\end{picture} \Bigr) \ = \ 
\begin{picture}(12,7)
\put(6,1){\circle{10}} 
\put(10.9,1){\circle*{1}}
\put(5,4.6){\small $\ast$}
\put(2.2,-2){\line(6,1){8.4}}
\end{picture} \ + \ 
\begin{picture}(12,7)
\put(6,1){\circle{10}} 
\put(10.9,1){\circle*{1}}
\put(5,-4.9){\small $\ast$}
\put(2.2,4){\line(6,-1){8.4}}
\end{picture} \ + \ 
\begin{picture}(12,7)
\put(6,1){\circle{10}} 
\put(10.9,1){\circle*{1}}
\bezier{120}(1.2,1)(4,1)(6,1)
\bezier{120}(6,1)(7.2,-1)(8.4,-3.1)
\bezier{120}(6,1)(7.2,3)(8.5,5.2)
\end{picture} \ + \ 
\begin{picture}(12,7)
\put(6,1){\circle{10}} 
\put(10.9,1){\circle*{1}}
\put(11,1){\line(-4,-3){6.3}}
\put(1.1,1){\line(4,3){6.3}}
\end{picture} \ + \ \begin{picture}(12,7)
\put(6,1){\circle{10}} 
\put(10.9,1){\circle*{1}}
\put(11,1){\line(-4,3){6.3}}
\put(1.2,1){\line(4,-3){6.2}}
\end{picture}$ 
\bigskip

$\partial \Bigl(\begin{picture}(12,7)
\put(6,1){\circle{10}} 
\put(10.9,1){\circle*{1}}
\put(3,-2.8){\line(-1,6){1.05}}
\put(8,-3.3){\line(1,4){1.8}}
\end{picture}\Bigr) \ = \ 
\begin{picture}(12,7)
\put(6,1){\circle{10}} 
\put(10.9,1){\circle*{1}}
\put(-0.1,0){\small $\ast$}
\put(3.2,-2.9){\line(3,4){5.9}}
\end{picture} \ + \ 
\begin{picture}(12,7)
\put(6,1){\circle{10}} 
\put(10.9,1){\circle*{1}}
\put(11,1){\line(-4,3){6.2}}
\put(1.2,1){\line(4,-3){6.2}}
\end{picture} \ + \ 
\begin{picture}(12,7)
\put(6,1){\circle{10}} 
\put(10.9,1){\circle*{1}}
\put(11,1){\line(-4,-3){6.2}}
\put(1.2,1){\line(4,3){6.2}}
\end{picture} \ + \ 
\begin{picture}(12,7)
\put(6,1){\circle{10}} 
\put(10.9,1){\circle*{1}}
\put(2.9,-2.7){\line(3,4){5.8}}
\put(3.3,-3.1){\line(3,4){5.8}}
\put(2.8,1.2){\tiny $-$}
\end{picture}$ 
\bigskip

 $\partial \Bigl(\begin{picture}(12,7)
\put(6,1){\circle{10}} 
\put(10.9,1){\circle*{1}}
\bezier{80}(3.1,-2.9)(4.5,-1)(5.7,0.6)
\bezier{80}(8.9,4.9)(7.5,3)(6.3,1.4)
\put(8.85,-2.75){\line(-3,4){5.7}}
\end{picture} \Bigr) \ = \ 
\begin{picture}(12,7)
\put(6,1){\circle{10}} 
\put(10.9,1){\circle*{1}}
\bezier{120}(1.2,1)(4,1)(6,1)
\bezier{120}(6,1)(7.2,-1)(8.5,-3)
\bezier{120}(6,1)(7.2,3)(8.5,5.1)
\end{picture} \ + \ 
\begin{picture}(12,7)
\put(6,1){\circle{10}} 
\put(10.9,1){\circle*{1}}
\put(2.9,-2.7){\line(3,4){5.8}}
\put(3.3,-3.1){\line(3,4){5.8}}
\put(2.8,1.2){\tiny $+$}
\end{picture} $
\bigskip

$\partial \Bigl( \begin{picture}(12,7)
\put(6,1){\circle{10}} 
\put(10.9,1){\circle*{1}}
\put(5,4.6){\small $\ast$}
\put(2.2,-2){\line(6,1){8.3}}
\end{picture} \Bigr) \ = \ 
\begin{picture}(12,7)
\put(6,1){\circle{10}} 
\put(10.9,1){\circle*{1}}
\put(2,-4.2){\small $\ast$}
\put(5,4.6){\small $\ast$}
\end{picture} \ + \ 
\begin{picture}(12,7)
\put(6,1){\circle{10}} 
\put(10.9,1){\circle*{1}}
\put(5,4.6){\small $\ast$}
\put(6.2,5.8){\line(1,-6){1.55}}
\end{picture} \ + \ 
\begin{picture}(11,7)
\put(5,1){\circle{10}} 
\put(10,1){\circle*{1}}
\put(2.9,4.2){\small $\ast$}
\put(10,1){\line(-4,-3){6.2}}
\end{picture} \ + \ 
\begin{picture}(11,7)
\put(5,1){\circle{10}} 
\put(10,1){\circle*{1}}
\put(9.6,0){\small $\ast$}
\put(5,-3.7){\line(1,5){1.85}}
\end{picture}$ 
\bigskip

$\partial \Bigl( \begin{picture}(12,7)
\put(6,1){\circle{10}} 
\put(10.9,1){\circle*{1}}
\put(5,-4.8){\small $\ast$}
\put(2.3,4){\line(6,-1){8.2}}
\end{picture} \Bigr) \ = \ 
\begin{picture}(11,7)
\put(5,1){\circle{10}} 
\put(10,1){\circle*{1}}
\put(9.6,0){\small $\ast$}
\put(5,-3.8){\line(1,5){1.85}}
\end{picture} \ + \ 
\begin{picture}(12,7)
\put(6,1){\circle{10}} 
\put(10.9,1){\circle*{1}}
\put(5,-5){\small $\ast$}
\put(6.2,-3.8){\line(1,4){2.25}}
\end{picture} \ + \ 
\begin{picture}(12,7)
\put(6,1){\circle{10}} 
\put(10.9,1){\circle*{1}}
\put(2,-4.1){\small $\ast$}
\put(5,4.6){\small $\ast$}
\end{picture} \ + \ 
\begin{picture}(11,7)
\put(5,1){\circle{10}} 
\put(10,1){\circle*{1}}
\put(5,-4.8){\small $\ast$}
\put(10,1){\line(-4,3){6.2}}
\end{picture}$
\bigskip

$\partial \Bigl( \begin{picture}(12,7)
\put(6,1){\circle{10}} 
\put(10.9,1){\circle*{1}}
\put(0,0){\small $\ast$}
\put(3.1,-2.9){\line(3,4){5.8}}
\end{picture} \Bigr) \ = \ 
\begin{picture}(11,7)
\put(5,1){\circle{10}} 
\put(10,1){\circle*{1}}
\put(5,-4.8){\small $\ast$}
\put(10,1){\line(-4,3){6.2}}
\end{picture} \ + \ 
\begin{picture}(12,7)
\put(6,1){\circle{10}} 
\put(10.9,1){\circle*{1}}
\put(5,-5){\small $\ast$}
\put(6.2,-3.8){\line(1,4){2.25}}
\end{picture} \ + \ 
\begin{picture}(12,7)
\put(6,1){\circle{10}} 
\put(10.9,1){\circle*{1}}
\put(5,4.6){\small $\ast$}
\put(6.2,5.8){\line(1,-6){1.55}}
\end{picture} \ + \ 
\begin{picture}(11,7)
\put(5,1){\circle{10}} 
\put(10,1){\circle*{1}}
\put(2.9,4.1){\small $\ast$}
\put(10,1){\line(-4,-3){6.2}}
\end{picture} \ + \ \begin{picture}(12,7)
\put(6,1){\circle{10}} \put(10.9,1){\circle*{1}}
\put(1.5,3.8){\small $\circledast$}
\end{picture}$ 
\bigskip

$\partial \Bigl( \begin{picture}(12,7)
\put(6,1){\circle{10}} 
\put(10.9,1){\circle*{1}}
\bezier{120}(1.2,1)(4,1)(6,1)
\bezier{120}(6,1)(7.2,-1)(8.4,-3.1)
\bezier{120}(6,1)(7.2,3)(8.4,5.1)
\end{picture} \Bigr) \ = 
\partial \Bigl( \begin{picture}(12,7)
\put(6,1){\circle{10}} 
\put(10.9,1){\circle*{1}}
\put(2.9,-2.7){\line(3,4){5.8}}
\put(3.3,-3.1){\line(3,4){5.8}}
\put(2.8,1.2){\tiny $+$}
\end{picture} \Bigr) \ = \ 
\begin{picture}(12,7)
\put(6,1){\circle{10}} 
\put(10.9,1){\circle*{1}}
\put(5,4.6){\small $\ast$}
\put(6.2,5.8){\line(1,-6){1.55}}
\end{picture} \ + \ 
\begin{picture}(12,7)
\put(6,1){\circle{10}} 
\put(10.9,1){\circle*{1}}
\put(5,-5){\small $\ast$}
\put(6.2,-3.8){\line(1,4){2.25}}
\end{picture} $ 
\bigskip

$\partial \Bigl( \begin{picture}(12,7)
\put(6,1){\circle{10}} 
\put(10.9,1){\circle*{1}}
\put(2.9,-2.7){\line(3,4){5.8}}
\put(3.3,-3.1){\line(3,4){5.8}}
\put(2.8,1.2){\tiny $-$}
\end{picture} \Bigr) \ = \ 
\begin{picture}(12,7)
\put(6,1){\circle{10}} 
\put(10.9,1){\circle*{1}}
\put(5,4.6){\small $\ast$}
\put(6.2,5.8){\line(1,-6){1.55}}
\end{picture} \ + \ 
\begin{picture}(12,7)
\put(6,1){\circle{10}} 
\put(10.9,1){\circle*{1}}
\put(5,-5){\small $\ast$}
\put(6.2,-3.8){\line(1,4){2.25}}
\end{picture} \ + \ 
\begin{picture}(12,7)
\put(6,1){\circle{10}} 
\put(10.9,1){\circle*{1}}
\put(1.6,3.8){\small $\circledast$}
\end{picture}$
\bigskip

$\partial \Bigl( \begin{picture}(12,7)
\put(6,1){\circle{10}} 
\put(10.9,1){\circle*{1}}
\put(2,-4.1){\small $\ast$}
\put(5,4.6){\small $\ast$}
\end{picture}\Bigr) \ = \ 
\begin{picture}(12,7)
\put(6,1){\circle{10}} 
\put(10.9,1){\circle*{1}}
\put(-0.5,0){\small $\ast$}
\put(1,0){\small $\ast$}
\end{picture} $ 
\bigskip

$\partial \Bigl( \begin{picture}(12,7)
\put(6,1){\circle{10}} 
\put(10.9,1){\circle*{1}}
\put(5,4.6){\small $\ast$}
\put(6.2,5.8){\line(1,-6){1.55}}
\end{picture} \Bigr) \ = \ 
\partial \Bigl( \begin{picture}(12,7)
\put(6,1){\circle{10}} 
\put(10.9,1){\circle*{1}}
\put(5,-5){\small $\ast$}
\put(6.2,-3.8){\line(1,4){2.25}}
\end{picture} \Bigr) \ = \ 
\begin{picture}(12,7)
\put(6,1){\circle{10}} 
\put(10.9,1){\circle*{1}}
\put(11,0){\small $\ast$}
\put(11,1){\line(-5,1){9.4}}
\end{picture} \ + \ 
\begin{picture}(12,7)
\put(6,1){\circle{10}} 
\put(10.9,1){\circle*{1}}
\put(0.2,-1.8){\small $\ast$}
\put(11,1){\line(-6,-1){10}}
\end{picture} \ + \ 
\begin{picture}(12,7)
\put(6,1){\circle{10}} 
\put(10.9,1){\circle*{1}}
\put(-0.5,0){\small $\ast$}
\put(1,0){\small $\ast$}
\end{picture}$ 
\bigskip

$\partial \Bigl( \begin{picture}(12,7)
\put(6,1){\circle{10}} \put(10.9,1){\circle*{1}}
\put(1.6,3.8){\small $\circledast$}
\end{picture}\Bigr) \ = \ 
\partial \Bigl( \begin{picture}(12,7)
\put(6,1){\circle{10}} 
\put(10.9,1){\circle*{1}}
\put(-0.5,0){\small $\ast$}
\put(1,0){\small $\ast$}
\end{picture}\Bigr) \ = \ 0$
\bigskip

$\partial \Bigl( \begin{picture}(12,7)
\put(6,1){\circle{10}} 
\put(10.9,1){\circle*{1}}
\put(11,1){\line(-4,-3){6.2}}
\put(1.2,1){\line(4,3){6.1}}
\end{picture} \Bigr) \ = \ 
\begin{picture}(11,7)
\put(5,1){\circle{10}} 
\put(10,1){\circle*{1}}
\put(9.6,0){\small $\ast$}
\put(5,-3.8){\line(1,5){1.85}}
\end{picture} \ + \ 
\begin{picture}(11,7)
\put(5,1){\circle{10}} 
\put(10,1){\circle*{1}}
\put(2.9,4.1){\small $\ast$}
\put(10,1){\line(-4,-3){6.2}}
\end{picture} \ + \ 
\begin{picture}(12,7)
\put(6,1){\circle{10}} 
\put(11,1){\circle*{1}}
\put(11,1.3){\line(-1,0){10}}
\put(11,0.7){\line(-1,0){10}}
\put(5,2){\tiny $-$}
\end{picture}$
\bigskip

$\partial \Bigl( \begin{picture}(12,7)
\put(6,1){\circle{10}} 
\put(10.9,1){\circle*{1}}
\put(11,1){\line(-4,3){6.2}}
\put(1.2,1){\line(4,-3){6.1}}
\end{picture} \Bigr) \ = \ 
\begin{picture}(11,7)
\put(5,1){\circle{10}} \put(10,1){\circle*{1}}
\put(5,-4.8){\small $\ast$}
\put(10,1){\line(-4,3){6.2}}
\end{picture} \ + \ 
\begin{picture}(11,7)
\put(5,1){\circle{10}} 
\put(10,1){\circle*{1}}
\put(9.6,0){\small $\ast$}
\put(5,-3.8){\line(1,5){1.85}}
\end{picture} \ + \ 
\begin{picture}(12,7)
\put(6,1){\circle{10}} 
\put(11,1){\circle*{1}}
\put(11,1.3){\line(-1,0){10}}
\put(11,0.7){\line(-1,0){10}}
\put(5,2){\tiny $-$}
\end{picture}$ 
\bigskip

$\partial \Bigl( \begin{picture}(12,7)
\put(6,1){\circle{10}} 
\put(10.9,1){\circle*{1}}
\bezier{80}(11,1)(8.7,0.5)(6.4,0)
\bezier{80}(5.6,-0.1)(3.65,-0.5)(1.6,-0.9)
\put(6,-3.9){\line(0,1){9.9}}
\end{picture} \Bigr) \ = \ 0$
\bigskip

$\partial \Bigl( \begin{picture}(11,7)
\put(5,1){\circle{10}} 
\put(10,1){\circle*{1}}
\put(2.9,4.1){\small $\ast$}
\put(10,1){\line(-4,-3){6.2}}
\end{picture} \Bigr) \ = \ 
\partial \Bigl( \begin{picture}(11,7)
\put(5,1){\circle{10}} 
\put(10,1){\circle*{1}}
\put(5,-4.8){\small $\ast$}
\put(10,1){\line(-4,3){6.2}}
\end{picture} \Bigr) \ = \ 
\begin{picture}(12,7)
\put(6,1){\circle{10}} 
\put(10.9,1){\circle*{1}}
\put(11,0){\small $\ast$}
\put(11,1){\line(-5,1){9.5}}
\end{picture} \ + \ 
\begin{picture}(12,7)
\put(6,1){\circle{10}} 
\put(10.9,1){\circle*{1}}
\put(0.2,-1.8){\small $\ast$}
\put(11,1){\line(-6,-1){10}}
\end{picture} \ + \ 
\begin{picture}(12,7)
\put(6,1){\circle{10}} 
\put(10.9,1){\circle*{1}}
\put(10.8,-0.1){\small $\ast$}
\put(0.3,1){\small $\ast$}
\end{picture} \ + \ 
\begin{picture}(12,7)
\put(6,1){\circle{10}} 
\put(10.9,1){\circle*{1}}
\put(9.6,0.2){\small $\circledast$}
\end{picture}$
\bigskip

$\partial \Bigl( \begin{picture}(11,7)
\put(5,1){\circle{10}} 
\put(10,1){\circle*{1}}
\put(9.6,0){\small $\ast$}
\put(5,-3.8){\line(1,5){1.87}}
\end{picture} \Bigr) \ = \ 
\begin{picture}(12,7)
\put(6,1){\circle{10}} 
\put(10.9,1){\circle*{1}}
\put(10.8,-0.1){\small $\ast$}
\put(0.1,1){\small $\ast$}
\end{picture} \ + \ 
\begin{picture}(12,7)
\put(6,1){\circle{10}} 
\put(10.9,1){\circle*{1}}
\put(9.6,0.2){\small $\circledast$}
\end{picture}$
\bigskip

$\partial \Bigl( \begin{picture}(12,7)
\put(6,1){\circle{10}} 
\put(11,1){\circle*{1}}
\bezier{100}(11,1)(9,1)(6,1)
\bezier{100}(6,1)(4.5,3)(3.1,4.9)
\bezier{100}(6,1)(4.5,-1)(3.1,-2.9)
\end{picture} \Bigr) \ = \ 
\begin{picture}(12,7)
\put(6,1){\circle{10}} 
\put(10.9,1){\circle*{1}}
\put(0.2,-1.8){\small $\ast$}
\put(11,1){\line(-6,-1){10}}
\end{picture}$
\bigskip

$\partial \Bigl( \begin{picture}(12,7)
\put(6,1){\circle{10}} 
\put(11,1){\circle*{1}}
\put(11,1.3){\line(-1,0){10}}
\put(11,0.7){\line(-1,0){10}}
\put(4.5,2){\tiny $+$}
\end{picture} \Bigr) \ = \ 
\partial \Bigl( \begin{picture}(12,7)
\put(6,1){\circle{10}} 
\put(11,1){\circle*{1}}
\put(11,1.3){\line(-1,0){10}}
\put(11,0.7){\line(-1,0){10}}
\put(4.5,2){\tiny $-$}
\end{picture} \Bigr) \ = \ 
\begin{picture}(12,7)
\put(6,1){\circle{10}} 
\put(10.9,1){\circle*{1}}
\put(0.2,-1.8){\small $\ast$}
\put(11,1){\line(-6,-1){10}}
\end{picture} \ + \ 
\begin{picture}(12,7)
\put(6,1){\circle{10}} 
\put(10.9,1){\circle*{1}}
\put(11,0){\small $\ast$}
\put(11,1){\line(-5,1){9.4}}
\end{picture}$
\bigskip

$\partial \Bigl( \begin{picture}(13,7)
\put(6,1){\circle{10}} \put(10.9,1){\circle*{1}}
\put(10.8,-0.1){\small $\ast$}
\put(0.3,1){\small $\ast$}
\end{picture} \Bigr) \ = \ \partial \Bigl( \begin{picture}(12,7)
\put(6,1){\circle{10}} \put(10.9,1){\circle*{1}}
\put(0.2,-1.8){\small $\ast$}
\put(11,1){\line(-6,-1){10}}
\end{picture} \Bigr) \ = \ \partial \Bigl( \begin{picture}(13,7)
\put(6,1){\circle{10}} \put(10.9,1){\circle*{1}}
\put(11,0){\small $\ast$}
\put(11,1){\line(-5,1){9.4}}
\end{picture} \Bigr) \ = \ \partial \Bigl( \begin{picture}(12,7)
\put(6,1){\circle{10}} \put(10.9,1){\circle*{1}}
\put(9,0.4){$\circledast$}
\end{picture} \Bigr) \ = \ 0$ .
\end{theorem}
\medskip

These formulas immediately imply the following detailing of Theorem \ref{mthmadd}.

\unitlength 0.7 mm
\begin{corollary}
\label{cor43}
$H_4\left(\overline{CD}_2(S^1), {\mathbb Z}_2\right) \simeq 0$. The group $H_3\left(\overline{CD}_2(S^1), {\mathbb Z}_2\right)$ is isomorphic to $ {\mathbb Z}_2$ and is generated by the cycle 
\begin{picture}(12,8)
\put(6,3){\circle{10}} 
\put(10.9,3){\circle*{1}}
\bezier{80}(11,3)(8.7,2.5)(6.4,2)
\bezier{80}(5.6,1.9)(3.65,1.5)(1.6,1.1)
\put(6,-1.9){\line(0,1){9.9}}
\end{picture}. The group $H_2\left(\overline{CD}_2(S^1), {\mathbb Z}_2\right)$ is isomorphic to ${\mathbb Z}_2$ and is generated by the chain 
\ $\begin{picture}(12,10)
\put(6,3){\circle{10}} \put(11,3){\circle*{1}}
\put(11,3.3){\line(-1,0){10}}
\put(11,2.7){\line(-1,0){10}}
\put(4.5,4){\tiny $+$}
\end{picture} \ + \ 
\begin{picture}(12,9)
\put(6,3){\circle{10}} 
\put(11,3){\circle*{1}}
\put(11,3.3){\line(-1,0){10}}
\put(11,2.7){\line(-1,0){10}}
\put(4.5,4){\tiny $-$}
\end{picture}$~. The group $H_1\left(\overline{CD}_2(S^1), {\mathbb Z}_2\right)$ is isomorphic to ${\mathbb Z}_2$ and is generated by either of the cells \ \begin{picture}(12,10)
\put(6,3){\circle{10}} \put(10.8,3){\circle*{1}}
\put(10.8,1.9){\small $\ast$}
\put(0.3,3){\small $\ast$}
\end{picture} \ or \ \begin{picture}(13,9)
\put(6,3){\circle{10}} \put(10.8,3){\circle*{1}}
\put(10,2.4){\small $\circledast$}
\end{picture}~. $H_0(\overline{CD}_2(S^1), {\mathbb Z}_2) \simeq {\mathbb Z}_2$. \hfill $\Box_{\ref{cor43}}$ \label{prmthmadd}$\Box_{\ref{mthmadd}}$
\end{corollary}

\begin{remark} \rm
Alternative cycles that also generate the groups $H_1\left(\overline{CD}_2(S^1), {\mathbb Z}_2\right)$ and $H_2\left(\overline{CD}_2(S^1), {\mathbb Z}_2\right)$ are described in \S\S~\ref{oreal} and \ref{gener}, see Corollary \ref{alter}.
\end{remark} 

\subsection{Proof of Theorem \ref{mainprop}}
\label{prompro}

Let ${\mathcal F}$ be an equilevel algebra of class $\overline{CD}_2(S^1)$.
If its support is a single point $A$, then, by Proposition \ref{trui}, ${\mathcal F}$ contains the ideal $\M_A^4$, the image of the algebra ${\mathcal F} \cap \M_A$ in the quotient space $ \M_A / \M_A^4 \simeq {\mathbb R}^3$ is a subspace of codimension 2 (and hence a line), and the algebra ${\mathcal F}$ is characterized by this subspace. All functions $f \in {\mathcal F}$ necessarily satisfy the condition $f'(A) = 0$, so this line lies in the subspace $ \M^2_A / \M_A^4 \simeq {\mathbb R}^2$. Thus the second linear condition on the functions $f \in {\mathcal F}$ can be formulated in the terms of this subspace and has the form $$ \xi f'''(A) + \eta f''(A)=0$$ for some point $(\xi:\eta) \in {\mathbb R}P^1$. If $\xi \neq 0$ then ${\mathcal F}$ is an algebra of type \begin{picture}(12,9)
\put(6,3){\circle{10}} 
\put(1.6,5.8){\small $\circledast$}
\end{picture} , with $\alpha = - \frac{\eta}{\omega}$. If $\xi=0$ then it is of type \ \begin{picture}(12,9)
\put(6,3){\circle{10}} 
\put(-0.5,2){\small $\ast$}
\put(1,2){\small $\ast$}
\end{picture}~. 
\smallskip

If the support of ${\mathcal F}$ consists of two points $A$ and $B$, then the family of matched configurations $\nu(t) = (x_1(t), \tilde x_1(t); x_2(t), \tilde x_2(t))$ tending to these two points can be only of three types:

1) both points of each chord $(x_1(t), \tilde x_1(t))$ or $(x_2(t), \tilde x_2(t))$ tend to one point,

2) the points of one chord tend to one point (say, $A$), and the points of the other chord tend to both $A$ and $B$,

3) both endpoints of each chord tend to both points of the support: for certainty, let $x_1$ and $x_2$ tend to $A$, while $\tilde x_1$ and $\tilde x_2$ tend to $B$.

In the first case, we have two conditions $f'(A)=0$, $f'(B)=0,$ and hence an algebra of type \begin{picture}(12,9)
\put(6,3){\circle{10}} 
\put(2,-2){\small $\ast$}
\put(5,6.6){\small $\ast$}
\end{picture}. \ In the second case, we have the conditions $f'(A)=0$ and $f(A)=f(B)$, and an algebra of type \begin{picture}(12,9)
\put(6,3){\circle{10}} 
\put(5,6.6){\small $\ast$}
\put(6.2,7.8){\line(1,-6){1.55}}
\end{picture}~. 
In the third case, the limit algebra is of class $\gimel (\alpha; A, B)$ 
where $$\alpha = \lim_{t \to 0}\frac{\tilde x_2(t)- \tilde x_1(t)}{x_2(t)- x_1(t)} \ . $$
So, it is a point of a cell of type \begin{picture}(12,8)
\put(6,3){\circle{10}} 
\put(2.85,-0.75){\line(3,4){5.8}}
\put(3.3,-1.1){\line(3,4){5.8}}
\put(2.8,3.2){\tiny $+$}
\end{picture} \ or
\ \begin{picture}(12,8)
\put(6,3){\circle{10}} 
\put(2.85,-0.75){\line(3,4){5.8}}
\put(3.3,-1.1){\line(3,4){5.8}}
\put(2.8,3.2){\tiny $-$}
\end{picture}
if this limit is finite and is not equal to 0. If this limit is equal to 0 or $\infty$, then we again obtain an algebra of class \begin{picture}(12,9)
\put(6,3){\circle{10}}
\put(5,6.6){\small $\ast$}
\put(6.2,7.8){\line(1,-6){1.5}}
\end{picture}~.

The cases of algebras with supports consisting of three or four points are trivial.

All of these algebras indeed appear as limits of algebras defined by 2-chord diagrams. The only questionable algebras of class \ 
\begin{picture}(12,9)
\put(6,3){\circle{10}} 
\put(1.6,5.8){\small $\circledast$}
\end{picture} \
with any given value of the parameter $\alpha$ appear as the limits of algebras of class \
\begin{picture}(12,10)
\put(6,2){\circle{10}} 
\put(0,2){\small $\ast$}
\put(3.3,-2){\line(3,4){5.9}}
\end{picture} \ defined by the conditions $$f'(A)=0, \ \ f(A-t)=f\left(A+t -\frac{\alpha}{3} t^2\right), \qquad t \to +0.$$ The latter algebras themselves are the limits of some algebras defined by chord diagrams, which hence approximate our algebra of class \begin{picture}(12,9)
\put(6,3){\circle{10}} 
\put(1.6,5.8){\small $\circledast$}
\end{picture} . For example, this algebra is approached by the family of algebras defined by the conditions 
$$f(A+t^3)=f(A-t^3), \ \ f(A-t)=f\left(A+t -\frac{\alpha}{3} t^2\right).$$

This algebra is also the limit of algebras of class \
\begin{picture}(12,10)
\put(6,3){\circle{10}} 
\put(2.85,-0.75){\line(3,4){5.8}}
\put(3.3,-1.1){\line(3,4){5.8}}
\put(2.8,3.2){\tiny $-$}
\end{picture} \
given by the conditions $$f(A-t)=f(A+t), \ \ f'(A-t)= \left(-1+\frac{2 \alpha}{3} t\right) f'(A+t). $$

All other incidences of cells are elementary. 

The definitions of all these cells immediately imply that, for each $i$-dimen\-si\-onal cell $\beta^i_\alpha$ of the list of Theorem \ref{mainprop}, there is a parameterization \begin{equation}
\varphi^i_\alpha: D^i \to \overline{CD}_2(S^1) \subset G_{\dim {\mathcal L}-2}(\mathcal L),
\label{paramcell}
\end{equation}
 where $D^i$ is the standard open unit ball in ${\mathbb R}^i$, ${\mathcal L}$ is the 2-representative subspace in $C^\infty(S^1, {\mathbb R})$ containing $\overline{CD}_2(S^1)$, and $\varphi_\alpha^i $ is a regular algebraic map that diffeomorphically maps $D^i$ onto $\beta^i_\alpha$. It remains to prove that the variety $\overline{CD}_2(S^1)$ has the structure of a CW-complex with these cells. 

\begin{lemma}
\label{le21}
For each $i$-dimensional cell described in Theorem \ref{mainprop}
there exists a continuous map $\chi: \bar {\mathfrak D}^i \to \overline{CD}_2(S^1)$ from a closed ball of dimension $i$, which maps the inner part of this ball homeomorphically onto our cell. \end{lemma}

\noindent
{\it Proof.} Let $K^i \subset D^i \times \overline{CD}_2(S^1)$ be the graph of the map (\ref{paramcell}). By construction and the Tarski--Seidenberg theorem,
it is a  non-singular semialgebraic subvariety of the algebraic manifold \ ${\mathbb R}^i \times G_{\dim {\mathcal L}-2}(\mathcal L)$. Let $\tilde K^i$ be the algebraic closure of $K^i$, and let $Y \subset \tilde K^i$ be the subvariety $\tilde K^i \cap (\partial D^i \times G_{\dim {\mathcal L}-2}(\mathcal L))$ of $\tilde K^i$.
Let \begin{equation}
\label{resol}
\xi: (\Xi, \tilde Y) \to (\tilde K^i, Y)
\end{equation} 
be a resolution of singularities (in the sense of \cite{hir}, \cite{spi}), such that $\Xi$ is a smooth manifold and $\tilde Y$ is a divisor with normal crossings in it. Let $\bar K^i \subset \Xi $ be the (topological) closure of the set $\xi^{-1}(K^i)$. Finally, let $\check K^i$ be the space of pairs of the form (a point of $\bar K^i$, a local connected component of the space $\xi^{-1}(K^i)$ near this point). (I believe that in all our cases $\check K^i = \bar K^i$, but it is tedious and unnecessary to verify this.)
 By construction, $\check K^i$ is a smooth manifold with corners, and its interior part is diffeomorphic to the open ball $D^i$. Moreover, the lifting of the squared norm function from the space ${\mathbb R}^i \supset D^i$ has in it the standard form: for any point of the boundary,
near which the manifold $\check K^i$ is distinguished by the conditions $z_1\geq 0, \dots, z_p \geq 0$ in some local coordinates in $\Xi$, this lifted function has the form $1 - \theta(z_1, \dots, z_i) \cdot z_1^{c_1} \cdot \ldots \cdot z_p^{c_p}$ for some smooth function $\theta >0$ and some natural numbers $c_1, \dots, c_p$. Then the boundary set of $\check K^i$ is homeomorphic to the standard $(i-1)$-dimensional sphere. Indeed, every level set of the lifted norm function of ${\mathbb R}^i$ that is sufficiently close to the boundary contracts homeomorphically onto this boundary set via the corresponding real Clemens structure; but such level sets are the diffeomorphic preimages of concentric spheres in ${\mathbb R}^i$. Thus, the entire variety $\check K^i$ is homeomorphic to the standard closed $i$-dimensional ball
$\bar {\mathfrak D}^i$. The desired map of this variety to $\overline{CD}_2(S^1)$ is defined by the composition of the resolution map (\ref{resol}) and the obvious projection ${\mathbb R}^i \times \overline{CD}_2(S^1) \to \overline{CD}_2(S^1)$. \hfill $\Box_{\ref{le21}}$
\medskip 

To verify that we really obtain a CW-complex, it
remains to prove that the boundary set $\bar \beta_\alpha^i \setminus \beta_\alpha^i$
of each of our cells $\beta_\alpha^i$ belongs to the union of cells of lower dimensions. This follows easily from the construction of these cells. 

All formulas of Theorem~\ref{mainprop} on the boundary operators also follow immediately from the construction of these cells. \hfill $\Box_{\ref{mainprop}}$ 

\subsection{Important filtration in $\overline{CD}_k(S^1)$}
\label{filth}

Let $M$ be an one-dimensional manifold and $\F \in \overline{CD}_k(M)$ be an equilevel subalgebra.

\begin{definition} \rm
\label{muldef}
The {\it multiplicity} of a support point $A$ of the algebra $\F$ is the minimal natural number $m$ such that the algebra $\F$ contains all functions  $f \in \M_A^m$ that are identically zero in neighborhoods of all other support points of $\F$. The {\it multiplicity of the algebra} $\F$ is the sum of the multiplicities of all its support points.
\end{definition}

\begin{example} \rm
All equilevel algebras of codimension two in $C^{\infty}(S^1, {\mathbb R})$ are of multiplicity four, except only for the algebras of classes 
\begin{picture}(12,10)
\put(6,4){\circle{10}} 
\bezier{120}(1.1,4)(4,4)(6,4)
\bezier{120}(6,4)(7.2,2)(8.4,0)
\bezier{120}(6,4)(7.2,6)(8.45,8.1)
\end{picture} ,
\begin{picture}(12,10)
\put(6,4){\circle{10}} 
\put(5,7.6){\small $\ast$}
\put(6.2,8.8){\line(1,-6){1.55}}
\end{picture} , and \ 
\begin{picture}(12,10)
\put(6,4){\circle{10}} 
\put(-0.5,3){\small $\ast$}
\put(1,3){\small $\ast$}
\end{picture}, which are of multiplicity three. In our cell decomposition of $\overline{CD}_2(S^1)$, these algebras of multiplicity three form eight cells
depending on the disposition of their supports with respect to the point $\bullet$.
The union of these cells is a $CW$-subcomplex in $\overline{CD}_k(S^1)$. 
The corresponding cell subcomplex mod 2 is acyclic in all dimensions except for 0 and 1. In these two dimensions, its homology groups are one-dimensional, and the group $H_1$ is generated by the class of the cell 
\begin{picture}(12,10)
\put(6,4){\circle{10}} 
\put(-0.5,3){\small $\ast$}
\put(1,3){\small $\ast$}
\put(11,4){\circle*{2}}
\end{picture}~.
 The quotient complex generated by cells of multiplicity four has mod 2 homology groups of ranks 0, 1, 2, 1 in dimensions 0, 1, 2, 3.
\end{example}

\begin{proposition}
\label{poro29}
For any natural numbers $p$ and $k$, the set of all algebras $\F$ of multiplicity $\leq p$ is closed in $\overline{CD}_k(M)$.
\end{proposition}

\noindent
{\it Proof.} Let $\varkappa: [0,\varepsilon) \to \overline{CD}_k(M)$ be a germ of an algebraic parametric curve such that $\varkappa(0) = \F$. Let $A$ be a support point of $\F$ and $\tilde A_1(t), \dots, \tilde A_r(t)$ be all support points of the algebras $\F_t \equiv \varkappa(t)$ tending to $A$ as $t$ tends to 0. Let $q_1, \dots, q_r$ be the multiplicities of these points for sufficiently small $t>0$. We need to show that the multiplicity of the support point $A$ of $\F$ is at most $q_1 + \dots + q_r$.
Let $f: M \to {\mathbb R}$ be a function that vanishes at $A$ with all derivatives up to degree $q-1$, where $q=q_1+ \dots + q_r,$ and is identically equal to zero in some neighborhoods of all other support points of $\F$. Then, $f$ can be represented in the form $f_0 \cdot x^{q} + \tilde f$ where $x$ is the local coordinate with the origin $A$, $f_0$ is a smooth function identically equal to 0 outside the coordinate neighborhood of $A$, and $\tilde f$ is a smooth function identically equal to 0 near all support points of $\F$ including $A$. For any $t \in (0,\varepsilon)$, take the function $$f_0 \cdot (x - x(\tilde A_1(t)))^{q_1} \cdot \ldots \cdot (x - x(\tilde A_r(t)))^{q_r} + \tilde f. $$ This function belongs to the algebra $\F_t$. These functions tend to $f$ when $t $ tends to $+0$, therefore $f \in \F$.
\hfill $\Box_{\ref{poro29}}$ \medskip

Thus, the sets of algebras of multiplicity $\leq p$ define an increasing filtration on the spaces $\overline{CD}_k(S^1)$. 
An important property of this filtration is that the set of all equilevel algebras of a fixed multiplicity in $\overline{CD}_k(S^1)$ is fibered over the source circle. 
Indeed, we can consider this circle as the set of complex numbers of unit norm. Let us associate to each equilevel subalgebra $\F \in \overline{CD}_k(S^1)$ the product of its support points taken with their multiplicities.
The resulting map $\overline{CD}_k(S^1) \to S^1$ is continuous on the set of equilevel algebras of any fixed multiplicity $p$. The restriction of this map to this set commutes with the simultaneous rotation of the source and the target circles by the angles $\alpha$ and $p \alpha$, respectively.

\begin{remark} \rm
There are other natural filtrations on the spaces $\overline{CD}_k(S^1)$, first of all the two-term filtration with the first term formed by the algebras of second type, and also the filtration by the number of support points.
\end{remark}

\section{Cycles in $\overline{CD}_k(S^1)$ for arbitrary $k$. Proof of Theorems \ref{estim2}, \ref{estim3}, and \ref{corm3}} 
\label{twoser}

\subsection{Cycles $\Sigma$ and $\aleph_k$.}
\label{gener}

Let $A$ and $B$ be two points in $S^1$, and $\Delta$ be an arbitrary local diffeomorphism $ (S^1, A) \to (S^1,B)$. For any natural number $k$, denote by $\gimel_k (\Delta; A, B)$ the subalgebra of codimension $k$ in $C^\infty(S^1, {\mathbb R})$ consisting of all functions $f$ such that the $(k-1)$-jets of the functions $f$ 
and $f \circ \Delta$ at the point $A$ are the same. 

\begin{lemma}
\label{le78}
For any local diffeomorphism $\Delta$, the algebra $\gimel_k (\Delta; A, B)$ belongs to $\overline{CD}_k(S^1)$.
\end{lemma}

\noindent
{\it Proof.} For any $t \in (0, \varepsilon)$, take for $\nu(t)$ the configuration $(A, B; A+t, \Delta (A+t); A+2t, \Delta(A+2t); \dots, A+(k-1)t, \Delta(A+(k-1)t))$. All elements of the algebra $\gimel_k(\Delta; A, B)$ obviously belong to the equilevel algebra defined by this curve $\nu$ as in \S~\ref{maindef}: for any function $f$ from this algebra, it is sufficient to take all functions $f_t$ equal to $f$. Having the same codimensions, these algebras are identical.
\hfill $\Box_{\ref{le78}}$ \medskip

It is easy to see that the algebra $\gimel_k(\Delta; A, B)$
depends only on the $(k-1)$-jet of the local diffeomorphism $\Delta$, so we will also use this notation assuming that $\Delta$ is such a jet, and not a local diffeomorphism itself. 
For $k=2$, the $(k-1)$-jets of local diffeomorphisms are just their differentials, which are described by a single number, $\Delta'(A)$. Thus, the algebras $\gimel_2(\Delta; A,B)$ are the same as the algebras $\gimel(\Delta'(A); A, B)$ defined in \S \ref{onnot}. 

For any two points $A \neq B$, the set of all algebras $\gimel_k (\Delta; A, B)$ is diffeomorphic to two copies of the space ${\mathbb R}^{k-1}$ that are distinguished by the sign of the differential of $\Delta$ at $A$. This set can be identified with a $(k-1)$-dimensional semialgebraic subset of the Grassmann manifold 
\begin{equation}G_{k-1}((\M_A \cap \M_B)/(\M^k_A \cap \M_B^k)) \simeq G_{k-1}({\mathbb R}^{2k-2}).
\label{gra}
\end{equation} 
Namely, each algebra $\F$ from this set is associated with its image in the quotient space $$(\F \cap \M_A \cap \M_B)/(\M_A^k \cap \M_B^k).$$

Denote by $\Sigma_k(A, B)$ the closure of this set in this Grassmannian.

\begin{proposition}
\label{prokuku}
The variety $\Sigma_k(A,B)$ is a cycle in the manifold $($\ref{gra}$)$.
\end{proposition}

\begin{lemma}
\label{lekuku}
Any codimension-$k$ equilevel subalgebra in $C^\infty(S^1, {\mathbb R})$, which
\begin{itemize}
\item[a)]
contains the ideal $\M_A^k \times \M_B^k$ and all the constant functions, and 
\item[b)] contains a function $f$ with $f(A)=f(B)=0$ and $f'(A) \neq 0 \neq f'(B)$, 
\end{itemize}
is an equilevel algebra of class $\gimel_k(\Delta; A, B)$. 
\end{lemma}

\noindent
{\it Proof.} Let $x \equiv \varphi-\varphi(A)$ and $y\equiv \varphi-\varphi(B)$ be local coordinates in $S^1$ with centers at $A$ and $B$. Norming the function $f$ and subtracting its appropriate degrees, we can assume that $f \equiv y$ in a neighborhood of the point $B$. Let $a_1 x + \dots + a_{k-1}x^{k-1}$ be the $(k-1)$-jet of $f$ at the point $A$. Then the functions $1, f, f^2, \dots, f^{k-1}$ all belong to the algebra $\gimel_k(\Delta; A, B)$ where $\Delta(x) \equiv a_1 x + \dots + a_{k-1}x^{k-1}$. These $k$  functions are linearly independent in the space $C^\infty(M, {\mathbb R})/(\M^k_A \cap \M_B^k)$, hence our algebra contains the  algebra  $\gimel_k(\Delta; A, B)$. 
Having the same codimensions, these algebras are identical. \hfill $\Box_{\ref{lekuku}}.$ 
\medskip

\noindent
{\it Proof of Proposition \ref{prokuku}.} By definition, the union of algebras $\gimel_k(\Delta; A, B)$ is a semialgebraic subset  in the manifold (\ref{gra}). Therefore, the sum of maximal strata of its algebraic closure is a cycle of the same dimension $k-1$. 
According to Lemma \ref{lekuku}, all points of this algebraic closure that are not of the class $\gimel_k(\Delta; A,B)$ are the algebras such that the derivatives of 
all their functions vanish at $A$ or $B$. Thus, they form a proper algebraic subset of dimension $\leq k-1$, in particular also define a (probably, empty) cycle of that dimension. Therefore, the closure of the difference of these cycles also defines a cycle. \hfill $\Box_{\ref{prokuku}}$
\medskip

Denote by $\aleph_k$ the union of the cycles $\Sigma_k(A, A+\pi) \subset \overline{CD}_k(S^1)$ over all pairs of {\em opposite} points of $S^1$. 

\begin{lemma}
\label{le79}
The $k$-th Stiefel--Whitney class $w_k(\N_k) \in H^k(\overline{CD}_k(S^1), {\mathbb Z}_2)$ of the normal bundle takes the non-zero value on the cycle $\aleph_k$.
\end{lemma}
 
\noindent
{\it Proof}. 
The restriction of the bundle $\N_k$ to the variety $\aleph_k$ has a cross-section
defined by the cosets of the function \ $\sin \varphi$ modulo the subalgebras defined by the points of the base of this bundle. This cross-section has a single intersection point with the zero section of this bundle. Indeed, such intersection points are the algebras $\gimel_k(\Delta; A, A+\pi) \in \aleph_k$ containing this function. If $\sin \varphi \in \gimel_k(\Delta; A, A+\pi)$, then 
$\sin(A+\pi)=\sin(A)$, i.e., $A = k \pi$, and then necessarily the $(k-1)$-jet of $\Delta$ at $A$ coincides with that of the diffeomorphism mapping each point $ \varphi$ to $\pi -\varphi.$ This intersection of cross-sections is transversal at a regular point of the variety $\aleph_k$. The boundary points of the varieties $\Sigma_k(A, A+\pi)$ (i.e., their points not of class $\gimel_k(\Delta; A, A+\pi))$ are algebras all whose elements $f$ satisfy pairs of conditions $f(A)=f(A+ \pi)$ and either $f'(A)=0$ or $f'(A+\pi)=0$, so our cross-section $\{\sin \varphi\}$ of the bundle $\N_k$ does not vanish at such points. Thus, the mod 2 Euler class of this bundle takes the non-zero value on the cycle $\aleph_k$. \hfill $\Box_{\ref{le79}}$ $\Box_{\ref{estim2}(1)}$ \label{prSWth(3)} $\Box_{\ref{SWth}(2)}$ $\Box_{\ref{thmsw}(3)}$

\begin{lemma}
\label{lem86}
The class $w_{k-1}(\N_k)$ takes the non-zero value on the cycle $\Sigma_k(A, A+\pi)$ for any $A \in S^1$. 
\end{lemma}

\noindent
{\it Proof.} 
Consider the two-dimensional constant bundle on $\overline{CD}_k(S^1)$, whose fibers
consist of the functions \begin{equation} 
\label{fou1a}
\lambda_1 \cos \varphi + \lambda_2 \sin \varphi \ .\end{equation}
Consider the tautological homomorphism of this bundle 
to $\N_k$ that for each point of the base 
 maps each function (\ref{fou1a}) to its coset modulo the corresponding subalgebra. 
For every $A \in [0,\pi]$ there is exactly one algebra $\gimel_k(\Delta; A, A+\pi) \in \Sigma_k(A, A+\pi)$ that has a nontrivial intersection with this subspace: namely, $\Delta$ should be the $(k-1)$-jet of the local diffeomorphism $(S^1, A) \to (S^1, A+\pi)$ that maps each point $A + \tau $ to $A + \pi - \tau$. In the other words, the tautological map of $\Sigma_k(A, A+\pi)$ to the Grassmann manifold of all codimension-$k$ subspaces in $C^\infty(S^1, {\mathbb R})$ has a single intersection point with the Schubert cell of all subspaces not transversal to the two-dimensional plane of functions (\ref{fou1a}). This intersection is transversal, therefore the first obstruction to the existence of two independent sections of the bundle $\N_k$ takes the non-zero value on the fundamental class of the variety $\Sigma_k(A, A+\pi)$. \hfill $\Box_{\ref{lem86}}$ \ $\Box_{\ref{SWth}(1)}$ 

\subsection{Cycles $\divideontimes$ and $\diamondsuit_k$}
\label{dvd}

Let $(A_1, A_2, \dots, A_k)$ be a set of $k$ different points in $S^1,$ and 
$\sigma = (\lambda_1: \lambda_2 : \ldots : \lambda_k)$ be an arbitrary point of ${\mathbb R}P^{k-1}$. Define the algebra $\xi (\sigma; A_1, \dots, A_k)$ by the conditions 
\begin{equation}
f(A_1) = f(A_2) =\dots = f(A_k), \ \ \ \lambda_1 f'(A_1) + \lambda_2 f'(A_2) + \ldots + \lambda_k f'(A_k) = 0.
\label{equ117}
\end{equation}

\begin{lemma}
\label{lekaka}
Each algebra $\xi (\sigma; A_1, \dots, A_k)$ belongs to the variety $\overline{CD}_k(S^1)$.
\end{lemma}

\noindent
{\it Proof}. Take $\nu(t) \equiv (A_1, A_2 + t \lambda_2; A_2, A_3 + t \lambda_3; \dots ; A_k, A_1+ t \lambda_1)$, cf. \S~\ref{sss811}. \hfill $\Box_{\ref{lekaka}}$
\medskip

For each set of points $A_1, \dots, A_k$, denote by $\divideontimes(A_1, \dots, A_k)$ the set of all algebras $\xi (\sigma; A_1, \dots, A_k)$ with arbitrary $\sigma$. This set is obviously diffeomorphic to ${\mathbb R} P^{k-1}$.

\begin{lemma}
\label{lem87}
In the restriction to any variety $\divideontimes(A_1, \dots, A_k)$, 

1$)$ all classes $w_i(\N_k)$, $i \geq 2,$ are trivial, 

2$)$ all classes $w_1^{r}(\N_k)$, $r = 1, 2, \dots, k-1$, are nontrivial.
\end{lemma} 

\noindent
{\it Proof.} On this variety, the bundle $\N_k$ splits into the direct sum of 

1) the constant $(k-1)$-dimensional bundle normal to the subspace in $C^\infty(S^1, {\mathbb R})$ defined by $k-1$ conditions $f(A_1) = \dots = f(A_k)$, and 

2) the one-dimensional bundle which is normal to the bundle of hyperplanes defined by the right equations (\ref{equ117}), in particular is equivalent to the tautological bundle over ${\mathbb R}P^{k-1}$. \hfill $\Box_{\ref{lem87}} \ \Box_{\ref{SWth}(1)} \ \Box_{\ref{estim2}(2)}$
\medskip

Statement 3 of Theorem \ref{estim2} follows immediately from Lemmas \ref{lem86} and \ref{lem87}. \hfill $\Box_{\ref{estim2}(3)}$
\medskip

Define the manifold $\diamondsuit_k \subset \overline{CD}_k(S^1)$ as the union of the manifolds 
\begin{equation}
\label{dz}
\divideontimes\left(A, A+ 2\pi/k, A+4\pi/k, \dots, A+2(k-1)\pi/k\right)
\end{equation}
 over all $A \in [0, 2\pi/k]$.

\begin{theorem}
\label{th41}
The group $H_3(\overline{CD}_3(S^1), {\mathbb Z}_2)$ is generated by the classes of cycles $\aleph_3$ and $\diamondsuit_3$.
\end{theorem}

This theorem is proved in \S~\ref{homapp}, see Propositions \ref{hpro3}, \ref{pro18}, and \ref{probet}.

 The manifolds (\ref{dz}) can be considered as the fibers of the fiber bundle 
\begin{equation}
\label{fb9}
\diamondsuit_k \to S^1/{\mathbb Z}_k 
\end{equation}
taking each such fiber to the class of the point $A$.

\begin{lemma}
\label{prt23}
The restrictions of the 
classes $w_i(\N_k)$ to the manifold $\diamondsuit_k$ are trivial for all $i>2$ if $k$ is even and for all $i>1$ if $k$ is odd.
\end{lemma}

\noindent
{\it Proof.} In the restriction to $\diamondsuit_k$, the bundle $N_k$ splits into the sum of 

1) the $(k-1)$-dimensional vector bundle that is constant along the fibers of the bundle (\ref{fb9}) and, for all points of any such fiber (\ref{dz}),
is normal to the subspace in $C^\infty(S^1, {\mathbb R})$ defined by the conditions $$f(A) = f(A+2\pi/k) = \dots = f(A+2(k-1)\pi/k);$$ 

2) the one-dimensional vector bundle whose restriction to such a fiber is isomorphic to the tautological bundle on the space ${\mathbb R} P^{k-1}$.

The first vector bundle is lifted from a bundle over the one-dimensional base $S^1/{\mathbb Z}_k$, so its classes $w_i$ with $i>1$ are trivial. For odd $k$ this bundle is orientable; thus, its first Stiefel--Whitney class is also trivial. The second bundle is one-dimensional; thus its classes $w_j$ with $j>1$ are also trivial. The statement of Lemma \ref{prt23} now follows from the Whitney multiplication formula, see \S~4 in \cite{MS}.
\hfill $\Box_{\ref{prt23}}$

\subsection{Proof of Theorem \ref{estim3}}
If $k$ is even, then the first bundle from the previous proof is non-orientable, and its first Stiefel--Whitney class in $H^1(\diamondsuit_k, {\mathbb Z}_2)$ is Poincar\'e dual to an arbitrary fiber of the fiber bundle (\ref{fb9}). Therefore, the cup product of the class $w_1$ of the first vector bundle and the class $w_1^{k-1}$ of the second one takes the non-zero value on the fundamental cycle of the manifold $\diamondsuit_k$. In particular, this cycle defines a nontrivial element of the group $H_k(\overline{CD}_k, {\mathbb Z}_2)$ for even $k$. If $k$ is even and greater than $2$, then this element is independent of the fundamental class of the variety $\aleph_k$, because, by Lemmas \ref{le79} and \ref{prt23}, the class $w_k(\N_k)$ takes different values on these two cycles. \hfill $\Box_{\ref{estim3}}$

\subsection{Proof of Theorem \ref{corm3}}

If all codimension $k$ equilevel subalgebras intersect the space ${\mathcal L}^N$ transversally, then these intersections form a $(N-k)$-dimensional vector bundle on $\overline{CD}_k(S^1)$. The sum of this bundle and the bundle $\N_k$ is isomorphic to the constant bundle with fiber ${\mathcal L}^N$. Therefore, the total Stiefel--Whitney class of this bundle of intersection spaces is the inverse of that of $\N_k$. By the proof of Lemma \ref{lem87}, the total Stiefel--Whitney class of the restriction of the bundle $\N_k$ 
to an arbitrary variety $\divideontimes(A_1, \dots, A_k) \simeq {\mathbb R}P^{k-1}$ is equal to $1+ \alpha$, where $\alpha$ is the generator of the group $H^1({\mathbb R}P^{k-1}, {\mathbb Z}_2)$. Thus, the restriction of the total Stiefel--Whitney class of the bundle of intersections to this variety is equal to $(1+\alpha)^{-1}= 1 + \alpha + \dots + \alpha^{k-1}$. The dimension of this bundle cannot be smaller than $k-1$, and the dimension of the space ${\mathcal L}^N$ cannot be smaller than $k-1+ \dim \N_k = 2k-1$. \hfill $\Box_{\ref{corm3}}$

\begin{problem} \rm
Investigate the set of all equilevel algebras of codimension six with four-point support that are defined by the curves $\nu: [0,\varepsilon) \to (S^1)^{12}$ such that the limit positions of six chords look as {\huge $\boxtimes$}.
\end{problem}

\section{Proof of other results on the cohomology of $\overline{CD}_2(S^1)$}
\label{other2}

\subsection{Cycle $\Gamma$}
\label{oreal}

\begin{definition} The one-dimensional cycle $\Gamma \subset \overline{CD}_2(S^1)$
is the set of all algebras consisting of functions whose derivative vanishes at some two opposite points of $S^1$.
\end{definition}

This cycle is parameterized by the space $S^1/{\mathbb Z_2}$ of such pairs of opposite points.

\begin{proposition}
\label{nz}
The first Stiefel--Whitney class of the normal bundle $\N_2$ on $\overline{CD}_2(S^1)$ takes the non-zero value on the cycle $\Gamma$.
\end{proposition}

\noindent
{\it Proof.}
Ordering the endpoints of an arbitrary diameter $(A, A+ \pi)$ of $S^1$ specifies a canonical frame and thus an orientation of the normal bundle $\N_2$ over the corresponding point of the manifold $\Gamma \subset \overline{CD}_2(S^1)$. The first (respectively, second) vector of this frame is the class of any function with $f'(A+\pi)=1, f'(A)=0$ (respectively, $f'(A+\pi)=0, f'(A)=1$). Moving the point $A$ by the angle $\pi$ permutes these two basic vectors and thus breaks the orientation of the bundle.
\hfill $\Box_{\ref{nz}}$ $\Box_{\ref{SWth}(1)}$

\begin{corollary}
\label{alter}
The cycle $ \Gamma$ is homologous to the cycle \ \begin{picture}(12,9)
\put(6,3){\circle{10}} \put(10.8,3){\circle*{1}}
\put(10.8,1.9){\small $\ast$}
\put(0.3,3){\small $\ast$}
\end{picture}~. The cycle $\aleph_2$ is homologous to the cycle $ \begin{picture}(12,9)
\put(6,3){\circle{10}} \put(11,3){\circle*{1}}
\put(11,3.3){\line(-1,0){9.8}}
\put(11,2.7){\line(-1,0){9.8}}
\put(5,4){\tiny $+$}
\end{picture} \ + \ \begin{picture}(12,10)
\put(6,3){\circle{10}} \put(11,3){\circle*{1}}
\put(11,3.3){\line(-1,0){9.8}}
\put(11,2.7){\line(-1,0){9.8}}
\put(5,4){\tiny $-$}
\end{picture} $~. 
\end{corollary}

\noindent {\it Proof.} In both cases, the two classes being compared are nontrivial elements of a group that is isomorphic to ${\mathbb Z}_2$. \hfill $\Box_{\ref{alter}}$
\medskip

\noindent
{\it Proof of Proposition \ref{BU}.} If a two-dimensional subspace of $C^\infty(S^1, {\mathbb R})$ does not contain nontrivial functions with derivatives vanishing at opposite points of the circle, then this subspace defines a trivialization of the restriction of the normal bundle to the manifold 
$\Gamma \subset \overline{CD}_2(S^1)$; this contradicts Proposition \ref{nz}. \hfill $\Box_{\ref{BU}}$

\subsection{Proof of Theorem \ref{mthmmult}}
\label{prmthmmult}

The 2-dimensional cycle $\aleph_2 \subset \overline{CD}_2(S^1)$ is fibered over $S^1/{\mathbb Z}_2$, its fiber over any pair of opposite points $(A, A+\pi) \subset S^1$ consists of all algebras $\gimel(\alpha; A, A+\pi)$, $\alpha \in {\mathbb R}P^1$ (see \S \ref{onnot}). 
It is easy to see that this fiber bundle is non-orientable and thus $\aleph_2$
is homeomorphic to the Klein bottle \ {\unitlength 0.3 mm
\begin{picture}(30,30)
\put(0,0){\vector(1,0){30}}
\put(0,0){\vector(0,1){30}}
\put(0,30){\vector(1,0){30}}
\put(30,30){\vector(0,-1){30}}
\put(1,14){\footnotesize $u$}
\put(14,1){\footnotesize $v$}
\end{picture}} \ . A generator $v$ of the group $H_1(\aleph_2, {\mathbb Z}_2) \simeq {\mathbb Z}_2^2$ is given by the cross-section of this 
fiber bundle, whose value at each point $\{A, A+\pi\}$ is the algebra
$\gimel(1; A, A+\pi)$. 
Consider the constant two-dimensional vector bundle on $\overline{CD}_2(S^1)$, whose fibers
consist of the functions 
(\ref{fou1a}). For each point of the cycle $\aleph_2$,
the {\em tautological homomorphism} of this bundle
to the bundle $\N_2$ maps each such function to its coset modulo the corresponding subalgebra.
 
This homomorphism is an isomorphism
 over the cycle $v$: indeed, no nontrivial function of the form (\ref{fou1a}) can satisfy both conditions $f(A) = f(A+\pi)$ and $f'(A) = f'(A+\pi)$ at the same point $A$. Thus, the basis cohomology class $W \equiv w_1(\N_2) \in H^1\left(\overline{CD}_2(S^1), {\mathbb Z}_2\right)$ takes zero value on the cycle $v$. Therefore, $v$ represents the zero homology class of $\overline{CD}_2(S^1)$.

Another generator $u$ of the group $H_1(\aleph_2, {\mathbb Z}_2)$ is defined by a fiber of the same fiber bundle $\aleph_2 \to S^1/{\mathbb Z}_2$: this fiber
consists of all algebras $\gimel(\alpha; 0, \pi)$ with arbitrary $\alpha \in {\mathbb R}P^1$. According to Lemma \ref{lem86}, the class $w_1(\N_2)$ takes the nonzero value on the cycle $u$.

So, the ring homomorphism 
\begin{equation}
\label{restr}
H^*(\overline{CD}_2(S^1), {\mathbb Z}_2) \to H^*(\aleph_2, {\mathbb Z}_2)\end{equation} defined by the inclusion $\aleph_2 \subset \overline{CD}_2(S^1)$
maps the basis class $W \in H^1(\overline{CD}_2(S^1), {\mathbb Z}_2)$ 
to the 1-cohomology class of the Klein bottle $\aleph_2$, which takes the value 1 on the generator $u$ and the value $0$ on the generator $v$. The square of the latter cohomology class is nontrivial in $H^2(\aleph_2, {\mathbb Z}_2)$. Due to the functoriality of the cup product, this squared class coincides with the image of the class $W^2$ under the map (\ref{restr}), hence $W^2 \neq 0$ in $H^2(\overline{CD}_2(S^1), {\mathbb Z}_2)$. \hfill $\Box_{\ref{mthmmult}(1)}$ \medskip

This is sufficient to prove Corollary \ref{corSW}. Indeed, according to Lemma \ref{le79}, the class $w_2(\N_2)$ is a nontrivial element of the group $H^2(\overline{CD}_2(S^1), {\mathbb Z}_2) \simeq {\mathbb Z}_2$, so it is equal to $W^2$. \hfill $\Box_{\ref{corSW}}$
\medskip

Therefore, the total Stiefel--Whitney class of the third Cartesian power $(\N_2)^3$ of the normal bundle is equal to $(1+W+W^2)^3 \equiv 1+W+W^3$. To prove the second statement of Theorem \ref{mthmmult}, it remains to prove the following lemma.

\begin{lemma}
\label{le56}
$w_3((\N_2)^3)) =0.$
\end{lemma}

\noindent
{\it Proof}. Consider the six-dimensional subspace of the space $C^\infty (S^1, {\mathbb R^3}) \equiv (C^{\infty}(S^1, {\mathbb R}))^3$, which consists of 
maps $(f_1, f_2, f_3)$ whose three components $f_i$ are the functions of the form (\ref{fou1a}). This subspace defines a trivialization of the restriction of the bundle $(\N_3)^3$
to the closure of the cell 
\unitlength 0.7mm
 \begin{picture}(12,10)
\put(6,3){\circle{10}} \put(10.8,3){\circle*{1}}
\bezier{80}(11,3)(8.7,2.5)(6.4,2)
\bezier{80}(5.6,1.9)(3.65,1.5)(1.6,1.1)
\put(6,-1.9){\line(0,1){10}}
\end{picture} , which generates the group $H_3\left(\overline{CD}_2(S^1), {\mathbb Z}_2\right)$. Indeed, no nontrivial function of the form (\ref{fou1a}) can belong to a subalgebra from the closure of the cell 
 \begin{picture}(12,10)
\put(6,3){\circle{10}} \put(10.8,3){\circle*{1}}
\bezier{80}(11,3)(8.7,2.5)(6.4,2)
\bezier{80}(5.6,1.9)(3.65,1.5)(1.6,1.1)
\put(6,-1.9){\line(0,1){10}}
\end{picture} . \hfill $\Box_{\ref{le56}}$ $ \Box_{\ref{mthmmult}(2)}$

\subsection{Proof of Theorem \ref{corm}}
\label{prcorm}

\begin{lemma}
\label{regul}
 The nontrivial element of the group $H_2\left(\overline{CD}_2(S^1), {\mathbb Z}_2\right)$ can be represented by a two-dimensional cycle lying in the subspace $CD_2(S^1) \subset \overline{CD}_2(S^1)$, see Definition \ref{defcd}. 
\end{lemma}

\noindent
{\it Proof.} Let $\varepsilon$ be a small positive number.
For any point $\gimel(\alpha; A, A+\pi) $, $\alpha \in {\mathbb R}P^1$, of the cycle $\aleph_2$ consider the pair of chords 
\begin{equation}
\label{desin}
\left(A +\varepsilon \frac{\alpha}{|\alpha|+1}, A +\pi +\varepsilon \frac{1}{|\alpha|+1}\right) \ \mbox{ and } \ \left(A -\varepsilon \frac{\alpha}{|\alpha|+1}, A+\pi -\varepsilon \frac{1}{|\alpha|+1}\right).
\end{equation} 
These two chords never coincide (although they have a common endpoint when $\alpha=0$ or $\alpha=\infty$) and thus define a point of $CD_2(S^1)$.
The formulas (\ref{desin}) give the same result if we simultaneously replace $A$ with $A+\pi$ and $\alpha$ with $\alpha^{-1}$, so they define a map $\aleph_2 \to CD_2(S^1)$. Tending $\varepsilon$ to 0 defines a homotopy between this map 
and the identical embedding of the cycle $\aleph_2$, whose class 
in $H_2(\overline{CD}_2(S^1))$ is non-zero by Lemma \ref{le79}. \hfill $\Box_{\ref{regul}}$
\medskip

Let $U^2$ be a two-dimensional cycle in $CD_2(S^1)$. If for each its point the 
 third Cartesian power of the corresponding equilevel algebra
intersects the subspace ${\mathcal L}^7 \subset (C^{\infty}(S^1, {\mathbb R}))^3$ transversally, then these intersections form a 1-dimensional vector bundle on $U^2$, whose sum with  the bundle $\N_2^3$ is isomorphic to the constant bundle with fiber ${\mathcal L}^7$. By Theorems \ref{mthmmult} and \ref{SWth}, 
the total Stiefel--Whitney class of this 1-dimensional bundle is
 $(1 + W + W^2)^{-3} = (1+ W)^3 = 1 + W + W^2 $. 
This is impossible if the class $W^2 \in H^*(\overline{CD}_2(S^1), {\mathbb Z}_2)$ takes the non-zero value on the cycle $U^2$. Therefore, every 2-cycle in $CD_2(S^1)$, which is homologous to the cycle constructed in Lemma \ref{regul}, has a nonempty intersection with the set of algebras whose third Cartesian power intersects the space ${\mathcal L}^7$ non-transversally.
 \hfill $\Box_{\ref{corm}}$

\section{Classification of equilevel algebras of codimension three in $C^\infty(S^1, {\mathbb R})$. Cell structure on $\overline{CD}_3(S^1)$}
\label{cells}

\begin{theorem}
\label{thmcell}
The space $\overline{CD}_3(S^1)$ has the structure of a CW-complex with 260 cells. It has 130 cells of the first type $($i.e., with support not containing the distinguished point $\bullet \in S^1):$ 15 cells of dimension 6, 41 cells of dimension 5, 44 cells of dimension 4, 23 cells of dimension 3, six cells of dimension 2, and one cell of dimension 1. It also has the same numbers of cells of the second type, which are of dimensions 5, 4, 3, 2, 1, and 0, respectively. These cells are all described below in \S \ref{cells}. \end{theorem}

The incidence coefficients of these cells that define the boundary operators of the corresponding cell complex mod 2 are given in \S \ref{ic}.

\subsection{On notation} The cells of the first type in $\overline{CD}_3(S^1)$
are denoted by letters with indices like $A_{13},$ $c_2$, etc. The notations of the corresponding cells of the second type are obtained from them by drawing the bar above:
$\bar A_{13},$ $\bar c_2,$ ...
In the figures illustrating the algebras from the cells of the first type, the horizontal interval shows the open part $(0, 2\pi)$ of the source circle of functions $S^1 \to {\mathbb R}$. Correspondingly, instead of the chords in the circles used in the pictures for codimension two algebras we now draw the arcs connecting the points of this interval. A single arc connecting points $A$ and $B$ denotes the condition $f(A)=f(B)$. A tripod with feet at points $A,$ $ B,$ and $C$ denotes the double condition $f(A)=f(B)=f(C)$. A double arc connecting two ordinary (not marked by asterisks) points $A$ and $B$ denotes a family of pairs of conditions
 $f(A)=f(B)$ and $ f'(A) = \alpha f'(B)$, $\alpha \in {\mathbb R}^1 \setminus \{0\}$. The sign $+$ or $-$ over this double arc is the sign of the parameter $\alpha$ of all algebras forming the corresponding cell. The symbols $\ast$, $\ast\ast$, and $\circledast$ mean the same as in \S \ref{onnot}. More complicated notations are explained where they appear.

\unitlength 1.3mm

\begin{figure}[h]
\begin{picture}(24,10)
\put(0,5){\line(1,0){24}}
\put(2,5){\circle*{0.5}}
\put(6,5){\circle*{0.5}}
\put(10,5){\circle*{0.5}}
\put(14,5){\circle*{0.5}}
\put(18,5){\circle*{0.5}}
\put(22,5){\circle*{0.5}}
\put(4,5){\oval(4,3)[t]}
\put(12,5){\oval(4,3)[t]}
\put(20,5){\oval(4,3)[t]}
\put(10,1){\small $A_{11}$}
\end{picture} \qquad
\begin{picture}(24,10)
\put(0,5){\line(1,0){24}}
\put(2,5){\circle*{0.5}}
\put(6,5){\circle*{0.5}}
\put(10,5){\circle*{0.5}}
\put(14,5){\circle*{0.5}}
\put(18,5){\circle*{0.5}}
\put(22,5){\circle*{0.5}}
\put(4,5){\oval(4,3)[t]}
\put(14,5){\oval(8,3)[t]}
\put(18,5){\oval(8,3)[b]}
\put(10,1){\small $A_{12}$}
\end{picture} \qquad
\begin{picture}(24,10)
\put(0,5){\line(1,0){24}}
\put(2,5){\circle*{0.5}}
\put(6,5){\circle*{0.5}}
\put(10,5){\circle*{0.5}}
\put(14,5){\circle*{0.5}}
\put(18,5){\circle*{0.5}}
\put(22,5){\circle*{0.5}}
\put(4,5){\oval(4,3)[t]}
\put(16,5){\oval(12,3)[t]}
\put(16,5){\oval(4,3)[b]}
\put(10,1){\small $A_{13}$}
\end{picture}

\begin{picture}(24,10)
\put(0,5){\line(1,0){24}}
\put(2,5){\circle*{0.5}}
\put(6,5){\circle*{0.5}}
\put(10,5){\circle*{0.5}}
\put(14,5){\circle*{0.5}}
\put(18,5){\circle*{0.5}}
\put(22,5){\circle*{0.5}}
\put(6,5){\oval(8,3)[t]}
\put(10,5){\oval(8,3)[b]}
\put(20,5){\oval(4,3)[t]}
\put(10,1){\small $A_{21}$}
\end{picture} \qquad
\begin{picture}(24,10)
\put(0,5){\line(1,0){24}}
\put(2,5){\circle*{0.5}}
\put(6,5){\circle*{0.5}}
\put(10,5){\circle*{0.5}}
\put(14,5){\circle*{0.5}}
\put(18,5){\circle*{0.5}}
\put(22,5){\circle*{0.5}}
\put(6,5){\oval(8,3)[t]}
\put(12,5){\oval(12,3)[b]}
\put(18,5){\oval(8,3)[t]}
\put(10,1){\small $A_{22}$}
\end{picture} \qquad
\begin{picture}(24,10)
\put(0,5){\line(1,0){24}}
\put(2,5){\circle*{0.5}}
\put(6,5){\circle*{0.5}}
\put(10,5){\circle*{0.5}}
\put(14,5){\circle*{0.5}}
\put(18,5){\circle*{0.5}}
\put(22,5){\circle*{0.5}}
\put(6,5){\oval(8,3)[t]}
\put(14,5){\oval(16,3)[b]}
\put(16,5){\oval(4,3)[t]}
\put(10,1){\small $A_{23}$}
\end{picture} 

\begin{picture}(24,10)
\put(0,5){\line(1,0){24}}
\put(2,5){\circle*{0.5}}
\put(6,5){\circle*{0.5}}
\put(10,5){\circle*{0.5}}
\put(14,5){\circle*{0.5}}
\put(18,5){\circle*{0.5}}
\put(22,5){\circle*{0.5}}
\put(8,5){\oval(12,3)[t]}
\put(8,5){\oval(4,3)[b]}
\put(20,5){\oval(4,3)[t]}
\put(10,1){\small $A_{31}$}
\end{picture} \qquad
\begin{picture}(24,10)
\put(0,5){\line(1,0){24}}
\put(2,5){\circle*{0.5}}
\put(6,5){\circle*{0.5}}
\put(10,5){\circle*{0.5}}
\put(14,5){\circle*{0.5}}
\put(18,5){\circle*{0.5}}
\put(22,5){\circle*{0.5}}
\put(8,5){\oval(12,3)[t]}
\put(12,5){\oval(12,3)[b]}
\put(16,5){\oval(12,4)[t]}
\put(10,1){\small $A_{32}$}
\end{picture} \qquad
\begin{picture}(24,10)
\put(0,5){\line(1,0){24}}
\put(2,5){\circle*{0.5}}
\put(6,5){\circle*{0.5}}
\put(10,5){\circle*{0.5}}
\put(14,5){\circle*{0.5}}
\put(18,5){\circle*{0.5}}
\put(22,5){\circle*{0.5}}
\put(8,5){\oval(12,3)[t]}
\put(14,5){\oval(16,4)[b]}
\put(14,5){\oval(8,3)[b]}
\put(10,1){\small $A_{33}$}
\end{picture}

\begin{picture}(24,10)
\put(0,5){\line(1,0){24}}
\put(2,5){\circle*{0.5}}
\put(6,5){\circle*{0.5}}
\put(10,5){\circle*{0.5}}
\put(14,5){\circle*{0.5}}
\put(18,5){\circle*{0.5}}
\put(22,5){\circle*{0.5}}
\put(10,5){\oval(16,3)[t]}
\put(8,5){\oval(4,3)[b]}
\put(18,5){\oval(8,3)[b]}
\put(10,1){\small $A_{41}$}
\end{picture} \qquad
\begin{picture}(24,10)
\put(0,5){\line(1,0){24}}
\put(2,5){\circle*{0.5}}
\put(6,5){\circle*{0.5}}
\put(10,5){\circle*{0.5}}
\put(14,5){\circle*{0.5}}
\put(18,5){\circle*{0.5}}
\put(22,5){\circle*{0.5}}
\put(10,5){\oval(16,4)[t]}
\put(10,5){\oval(8,3)[t]}
\put(16,5){\oval(12,3)[b]}
\put(10,1){\small $A_{42}$}
\end{picture} \qquad
\begin{picture}(24,10)
\put(0,5){\line(1,0){24}}
\put(2,5){\circle*{0.5}}
\put(6,5){\circle*{0.5}}
\put(10,5){\circle*{0.5}}
\put(14,5){\circle*{0.5}}
\put(18,5){\circle*{0.5}}
\put(22,5){\circle*{0.5}}
\put(10,5){\oval(16,3)[t]}
\put(14,5){\oval(16,4)[b]}
\put(12,5){\oval(4,3)[b]}
\put(10,1){\small $A_{43}$}
\end{picture}

\begin{picture}(24,10)
\put(0,5){\line(1,0){24}}
\put(2,5){\circle*{0.5}}
\put(6,5){\circle*{0.5}}
\put(10,5){\circle*{0.5}}
\put(14,5){\circle*{0.5}}
\put(18,5){\circle*{0.5}}
\put(22,5){\circle*{0.5}}
\put(12,5){\oval(20,3)[t]}
\put(8,5){\oval(4,3)[b]}
\put(16,5){\oval(4,3)[b]}
\put(10,1){\small $A_{51}$}
\end{picture} \qquad
\begin{picture}(24,10)
\put(0,5){\line(1,0){24}}
\put(2,5){\circle*{0.5}}
\put(6,5){\circle*{0.5}}
\put(10,5){\circle*{0.5}}
\put(14,5){\circle*{0.5}}
\put(18,5){\circle*{0.5}}
\put(22,5){\circle*{0.5}}
\put(12,5){\oval(20,4)[t]}
\put(10,5){\oval(8,3)[t]}
\put(14,5){\oval(8,3)[b]}
\put(10,1){\small $A_{52}$}
\end{picture} \qquad
\begin{picture}(24,10)
\put(0,5){\line(1,0){24}}
\put(2,5){\circle*{0.5}}
\put(6,5){\circle*{0.5}}
\put(10,5){\circle*{0.5}}
\put(14,5){\circle*{0.5}}
\put(18,5){\circle*{0.5}}
\put(22,5){\circle*{0.5}}
\put(12,5){\oval(20,4)[t]}
\put(12,5){\oval(12,3)[b]}
\put(12,5){\oval(4,3)[t]}
\put(10,1){\small $A_{53}$}
\end{picture}
\caption{$A$}
\label{A}
\end{figure}

\subsection{Six-dimensional cells of $\overline{CD}_3(S^1)$}
\label{sixdim}
There are exactly 15 such cells called $A_{i j}$, see Fig.~\ref{A}. They are parameterized by the points of the supports of their algebras. The indices $i$ and $j$ in the notation of these cells indicate the partition of these six naturally ordered points of the interval $(0, 2\pi)$ into pairs. For example, $A_{32}$ denotes the picture in which the first point from the left is connected to the third remaining point, and the leftmost point remaining after that is connected to the second of the remaining three points.

\FloatBarrier

\unitlength 1.6mm

\begin{figure}[h]
\begin{picture}(20,9)
\put(0,5){\line(1,0){20}}
\put(2,5){\makebox(0,0)[cc]{$\ast$}}
\put(6,5){\circle*{0.5}}
\put(10,5){\circle*{0.5}}
\put(14,5){\circle*{0.5}}
\put(18,5){\circle*{0.5}}
\put(8,5){\oval(4,3)[t]}
\put(16,5){\oval(4,3)[t]}
\put(9,0){\small $B_1^1$}
\end{picture} \qquad
\begin{picture}(20,8)
\put(0,5){\line(1,0){20}}
\put(2,5){\makebox(0,0)[cc]{$\ast$}}
\put(6,5){\circle*{0.5}}
\put(10,5){\circle*{0.5}}
\put(14,5){\circle*{0.5}}
\put(18,5){\circle*{0.5}}
\put(10,5){\oval(8,3)[t]}
\put(14,5){\oval(8,3)[b]}
\put(9,0){\small $B_1^2$}
\end{picture} \qquad
\begin{picture}(20,8)
\put(0,5){\line(1,0){20}}
\put(2,5){\makebox(0,0)[cc]{$\ast$}}
\put(6,5){\circle*{0.5}}
\put(10,5){\circle*{0.5}}
\put(14,5){\circle*{0.5}}
\put(18,5){\circle*{0.5}}
\put(12,5){\oval(12,3)[t]}
\put(12,5){\oval(4,3)[b]}
\put(9,0){\small $B_1^3$}
\end{picture}

\begin{picture}(20,9)
\put(0,5){\line(1,0){20}}
\put(2,5){\circle*{0.5}}
\put(6,5){\makebox(0,0)[cc]{$\ast$}}
\put(10,5){\circle*{0.5}}
\put(14,5){\circle*{0.5}}
\put(18,5){\circle*{0.5}}
\put(6,5){\oval(8,3)[t]}
\put(16,5){\oval(4,3)[t]}
\put(9,0){\small $B_2^1$}
\end{picture} \qquad
\begin{picture}(20,8)
\put(0,5){\line(1,0){20}}
\put(2,5){\circle*{0.5}}
\put(6,5){\makebox(0,0)[cc]{$\ast$}}
\put(10,5){\circle*{0.5}}
\put(14,5){\circle*{0.5}}
\put(18,5){\circle*{0.5}}
\put(8,5){\oval(12,3)[t]}
\put(14,5){\oval(8,3)[b]}
\put(9,0){\small $B_2^2$}
\end{picture} \qquad
\begin{picture}(20,8)
\put(0,5){\line(1,0){20}}
\put(2,5){\circle*{0.5}}
\put(6,5){\makebox(0,0)[cc]{$\ast$}}
\put(10,5){\circle*{0.5}}
\put(14,5){\circle*{0.5}}
\put(18,5){\circle*{0.5}}
\put(10,5){\oval(16,3)[t]}
\put(12,5){\oval(4,3)[b]}
\put(9,0){\small $B_2^3$}
\end{picture}

\begin{picture}(20,9)
\put(0,5){\line(1,0){20}}
\put(2,5){\circle*{0.5}}
\put(6,5){\circle*{0.5}}
\put(10,5){\makebox(0,0)[cc]{$\ast$}}
\put(14,5){\circle*{0.5}}
\put(18,5){\circle*{0.5}}
\put(4,5){\oval(4,3)[t]}
\put(16,5){\oval(4,3)[t]}
\put(9,0){\small $B_3^1$}
\end{picture} \qquad
\begin{picture}(20,8)
\put(0,5){\line(1,0){20}}
\put(2,5){\circle*{0.5}}
\put(6,5){\circle*{0.5}}
\put(10,5){\makebox(0,0)[cc]{$\ast$}}
\put(14,5){\circle*{0.5}}
\put(18,5){\circle*{0.5}}
\put(8,5){\oval(12,3)[t]}
\put(12,5){\oval(12,3)[b]}
\put(9,0){\small $B_3^2$}
\end{picture} \qquad
\begin{picture}(20,8)
\put(0,5){\line(1,0){20}}
\put(2,5){\circle*{0.5}}
\put(6,5){\circle*{0.5}}
\put(10,5){\makebox(0,0)[cc]{$\ast$}}
\put(14,5){\circle*{0.5}}
\put(18,5){\circle*{0.5}}
\put(10,5){\oval(16,3)[t]}
\put(10,5){\oval(8,3)[b]}
\put(9,0){\small $B_3^3$}
\end{picture}

\begin{picture}(20,9)
\put(0,5){\line(1,0){20}}
\put(2,5){\circle*{0.5}}
\put(6,5){\circle*{0.5}}
\put(10,5){\circle*{0.5}}
\put(14,5){\makebox(0,0)[cc]{$\ast$}}
\put(18,5){\circle*{0.5}}
\put(4,5){\oval(4,3)[t]}
\put(14,5){\oval(8,3)[t]}
\put(9,0){\small $B_4^1$}
\end{picture} \qquad
\begin{picture}(20,8)
\put(0,5){\line(1,0){20}}
\put(2,5){\circle*{0.5}}
\put(6,5){\circle*{0.5}}
\put(10,5){\circle*{0.5}}
\put(14,5){\makebox(0,0)[cc]{$\ast$}}
\put(18,5){\circle*{0.5}}
\put(6,5){\oval(8,3)[t]}
\put(12,5){\oval(12,3)[b]}
\put(9,0){\small $B_4^2$}
\end{picture} \qquad
\begin{picture}(20,8)
\put(0,5){\line(1,0){20}}
\put(2,5){\circle*{0.5}}
\put(6,5){\circle*{0.5}}
\put(10,5){\circle*{0.5}}
\put(14,5){\makebox(0,0)[cc]{$\ast$}}
\put(18,5){\circle*{0.5}}
\put(10,5){\oval(16,3)[t]}
\put(8,5){\oval(4,3)[b]}
\put(9,0){\small $B_4^3$}
\end{picture}

\begin{picture}(20,9)
\put(0,5){\line(1,0){20}}
\put(2,5){\circle*{0.5}}
\put(6,5){\circle*{0.5}}
\put(10,5){\circle*{0.5}}
\put(14,5){\circle*{0.5}}
\put(18,5){\makebox(0,0)[cc]{$\ast$}}
\put(4,5){\oval(4,3)[t]}
\put(12,5){\oval(4,3)[t]}
\put(9,0){\small $B_5^1$}
\end{picture} \qquad
\begin{picture}(20,8)
\put(0,5){\line(1,0){20}}
\put(2,5){\circle*{0.5}}
\put(6,5){\circle*{0.5}}
\put(10,5){\circle*{0.5}}
\put(14,5){\circle*{0.5}}
\put(18,5){\makebox(0,0)[cc]{$\ast$}}
\put(6,5){\oval(8,3)[t]}
\put(10,5){\oval(8,3)[b]}
\put(9,0){\small $B_5^2$}
\end{picture} \qquad
\begin{picture}(20,8)
\put(0,5){\line(1,0){20}}
\put(2,5){\circle*{0.5}}
\put(6,5){\circle*{0.5}}
\put(10,5){\circle*{0.5}}
\put(14,5){\circle*{0.5}}
\put(18,5){\makebox(0,0)[cc]{$\ast$}}
\put(8,5){\oval(12,3)[t]}
\put(8,5){\oval(4,3)[b]}
\put(9,0){\small $B_5^3$}
\end{picture}
\caption{$B$}
\label{B}
\end{figure}

\subsection{Five-dimensional cells}
\label{fivedim}

There are 41 such cells of the first type, see Figs.~\ref{B}, \ref{C}, \ref{D}, and \ref{E}. Additionally, there are 15 \ five-dimensional cells $\bar A_{i j}$ of the second type, which correspond to the cells described in the previous subsection.

\begin{figure}[h]
\begin{picture}(16,9)
\put(0,5){\line(1,0){16}}
\put(2,5){\circle*{0.5}}
\put(5,5){\circle*{0.5}}
\put(8,5){\circle*{0.5}}
\put(11,5){\circle*{0.5}}
\put(14,5){\circle*{0.5}}
\put(3.5,5){\oval(3,3)[b]}
\put(8,8){\line(0,-1){3}}
\put(8,8){\line(1,-1){3}}
\put(8,8){\line(2,-1){6}}
\put(6,0){\small $C_{12}$}
\end{picture} \quad
\begin{picture}(16,8)
\put(0,5){\line(1,0){16}}
\put(2,5){\circle*{0.5}}
\put(5,5){\circle*{0.5}}
\put(8,5){\circle*{0.5}}
\put(11,5){\circle*{0.5}}
\put(14,5){\circle*{0.5}}
\put(5,5){\oval(6,3)[b]}
\put(8,8){\line(-1,-1){3}}
\put(8,8){\line(1,-1){3}}
\put(8,8){\line(2,-1){6}}
\put(6,0){\small $C_{13}$}
\end{picture} \quad
\begin{picture}(16,8)
\put(0,5){\line(1,0){16}}
\put(2,5){\circle*{0.5}}
\put(5,5){\circle*{0.5}}
\put(8,5){\circle*{0.5}}
\put(11,5){\circle*{0.5}}
\put(14,5){\circle*{0.5}}
\put(6.5,5){\oval(9,3)[b]}
\put(8,8){\line(-1,-1){3}}
\put(8,8){\line(0,-1){3}}
\put(8,8){\line(2,-1){6}}
\put(6,0){\small $C_{14}$} 
\end{picture} 
\quad
\begin{picture}(16,8)
\put(0,5){\line(1,0){16}}
\put(2,5){\circle*{0.5}}
\put(5,5){\circle*{0.5}}
\put(8,5){\circle*{0.5}}
\put(11,5){\circle*{0.5}}
\put(14,5){\circle*{0.5}}
\put(8,5){\oval(12,3)[b]}
\put(8,8){\line(-1,-1){3}}
\put(8,8){\line(0,-1){3}}
\put(8,8){\line(1,-1){3}}
\put(6,0){\small $C_{15}$}
\end{picture}

\begin{picture}(16,9)
\put(0,5){\line(1,0){16}}
\put(2,5){\circle*{0.5}}
\put(5,5){\circle*{0.5}}
\put(8,5){\circle*{0.5}}
\put(11,5){\circle*{0.5}}
\put(14,5){\circle*{0.5}}
\put(6.5,5){\oval(3,3)[b]}
\put(8,8){\line(-2,-1){6}}
\put(8,8){\line(1,-1){3}}
\put(8,8){\line(2,-1){6}}
\put(6,0){\small $C_{23}$}
\end{picture} \qquad
\begin{picture}(16,8)
\put(0,5){\line(1,0){16}}
\put(2,5){\circle*{0.5}}
\put(5,5){\circle*{0.5}}
\put(8,5){\circle*{0.5}}
\put(11,5){\circle*{0.5}}
\put(14,5){\circle*{0.5}}
\put(8,5){\oval(6,3)[b]}
\put(8,8){\line(0,-1){3}}
\put(8,8){\line(-2,-1){6}}
\put(8,8){\line(2,-1){6}}
\put(6,0){\small $C_{24}$}
\end{picture} \qquad
\begin{picture}(16,8)
\put(0,5){\line(1,0){16}}
\put(2,5){\circle*{0.5}}
\put(5,5){\circle*{0.5}}
\put(8,5){\circle*{0.5}}
\put(11,5){\circle*{0.5}}
\put(14,5){\circle*{0.5}}
\put(9.5,5){\oval(9,3)[b]}
\put(8,8){\line(0,-1){3}}
\put(8,8){\line(1,-1){3}}
\put(8,8){\line(-2,-1){6}}
\put(6,0){\small $C_{25}$}
\end{picture}

\begin{picture}(16,10)
\put(0,5){\line(1,0){16}}
\put(5,5){\circle*{0.5}}
\put(2,5){\circle*{0.5}}
\put(8,5){\circle*{0.5}}
\put(11,5){\circle*{0.5}}
\put(14,5){\circle*{0.5}}
\put(9.5,5){\oval(3,3)[b]}
\put(8,8){\line(-2,-1){6}}
\put(8,8){\line(-1,-1){3}}
\put(8,8){\line(2,-1){6}}
\put(6,0){\small $C_{34}$}
\end{picture} \qquad
\begin{picture}(16,8)
\put(0,5){\line(1,0){16}}
\put(5,5){\circle*{0.5}}
\put(2,5){\circle*{0.5}}
\put(8,5){\circle*{0.5}}
\put(11,5){\circle*{0.5}}
\put(14,5){\circle*{0.5}}
\put(11,5){\oval(6,3)[b]}
\put(8,8){\line(-2,-1){6}}
\put(8,8){\line(-1,-1){3}}
\put(8,8){\line(1,-1){3}}
\put(6,0){\small $C_{35}$}
\end{picture} \qquad
\begin{picture}(16,8)
\put(0,5){\line(1,0){16}}
\put(2,5){\circle*{0.5}}
\put(5,5){\circle*{0.5}}
\put(8,5){\circle*{0.5}}
\put(11,5){\circle*{0.5}}
\put(14,5){\circle*{0.5}}
\put(12.5,5){\oval(3,3)[b]}
\put(8,8){\line(-2,-1){6}}
\put(8,8){\line(-1,-1){3}}
\put(8,8){\line(0,-1){3}}
\put(6,0){\small $C_{45}$}
\end{picture}
\caption{$C$}
\label{C}
\end{figure}

\begin{figure}[h]
\begin{picture}(16,9)
\put(0,5){\line(1,0){16}}
\put(2,5){\circle*{0.5}}
\put(6,5){\circle*{0.5}}
\put(10,5){\circle*{0.5}}
\put(14,5){\circle*{0.5}}
\put(3,7.5){\small $+$}
\put(4,5){\oval(4.3,3.8)[t]}
\put(4,5){\oval(3.7,3.2)[t]}
\put(12,5){\oval(4,3)[b]}
\put(6,1){\small $D_{12}^+$}
\end{picture} \quad
\begin{picture}(16,8)
\put(0,5){\line(1,0){16}}
\put(2,5){\circle*{0.5}}
\put(6,5){\circle*{0.5}}
\put(10,5){\circle*{0.5}}
\put(14,5){\circle*{0.5}}
\put(3,7){\small $-$}
\put(4,5){\oval(4.3,3.8)[t]}
\put(4,5){\oval(3.7,3.2)[t]}
\put(12,5){\oval(4,3)[b]}
\put(6,1){\small $D_{12}^-$}
\end{picture} \quad
\begin{picture}(16,8)
\put(0,5){\line(1,0){16}}
\put(2,5){\circle*{0.5}}
\put(6,5){\circle*{0.5}}
\put(10,5){\circle*{0.5}}
\put(14,5){\circle*{0.5}}
\put(11,7.5){\small $+$}
\put(12,5){\oval(4.3,3.8)[t]}
\put(12,5){\oval(3.7,3.2)[t]}
\put(4,5){\oval(4,3)[b]}
\put(6,1){\small $D_{34}^+$}
\end{picture}
\label{D34+} \quad
\begin{picture}(16,8)
\put(0,5){\line(1,0){16}}
\put(2,5){\circle*{0.5}}
\put(6,5){\circle*{0.5}}
\put(10,5){\circle*{0.5}}
\put(14,5){\circle*{0.5}}
\put(11,7){\small $-$}
\put(12,5){\oval(4.3,3.8)[t]}
\put(12,5){\oval(3.7,3.2)[t]}
\put(4,5){\oval(4,3)[b]}
\put(6,1){\small $D_{34}^-$}
\end{picture}

\begin{picture}(16,10)
\put(0,5){\line(1,0){16}}
\put(2,5){\circle*{0.5}}
\put(6,5){\circle*{0.5}}
\put(10,5){\circle*{0.5}}
\put(14,5){\circle*{0.5}}
\put(5,7.5){\small $+$}
\put(6,5){\oval(8.3,3.8)[t]}
\put(6,5){\oval(7.7,3.2)[t]}
\put(10,5){\oval(8,3)[b]}
\put(6,1){\small $D_{13}^+$}
\end{picture} \quad
\begin{picture}(16,8)
\put(0,5){\line(1,0){16}}
\put(2,5){\circle*{0.5}}
\put(6,5){\circle*{0.5}}
\put(10,5){\circle*{0.5}}
\put(14,5){\circle*{0.5}}
\put(5,7){\small $-$}
\put(6,5){\oval(8.3,3.8)[t]}
\put(6,5){\oval(7.7,3.2)[t]}
\put(10,5){\oval(8,3)[b]}
\put(6,1){\small $D_{13}^-$}
\end{picture} \quad
\begin{picture}(16,8)
\put(0,5){\line(1,0){16}}
\put(2,5){\circle*{0.5}}
\put(6,5){\circle*{0.5}}
\put(10,5){\circle*{0.5}}
\put(14,5){\circle*{0.5}}
\put(9,7.5){\small $+$}
\put(10,5){\oval(8.3,3.8)[t]}
\put(10,5){\oval(7.7,3.2)[t]}
\put(6,5){\oval(8,3)[b]}
\put(6,1){\small $D_{24}^+$}
\end{picture} \quad
\begin{picture}(16,8)
\put(0,5){\line(1,0){16}}
\put(2,5){\circle*{0.5}}
\put(6,5){\circle*{0.5}}
\put(10,5){\circle*{0.5}}
\put(14,5){\circle*{0.5}}
\put(9,7){\small $-$}
\put(10,5){\oval(8.3,3.8)[t]}
\put(10,5){\oval(7.7,3.2)[t]}
\put(6,5){\oval(8,3)[b]}
\put(6,1){\small $D_{24}^-$}
\end{picture}

\begin{picture}(16,10)
\put(0,5){\line(1,0){16}}
\put(2,5){\circle*{0.5}}
\put(6,5){\circle*{0.5}}
\put(10,5){\circle*{0.5}}
\put(14,5){\circle*{0.5}}
\put(7,7.5){\small $+$}
\put(8,5){\oval(12.3,3.8)[t]}
\put(8,5){\oval(11.7,3.2)[t]}
\put(8,5){\oval(4,3)[b]}
\put(6,1){\small $D_{14}^+$}
\end{picture} \quad
\begin{picture}(16,8)
\put(0,5){\line(1,0){16}}
\put(2,5){\circle*{0.5}}
\put(6,5){\circle*{0.5}}
\put(10,5){\circle*{0.5}}
\put(14,5){\circle*{0.5}}
\put(7,7){\small $-$}
\put(8,5){\oval(12.3,3.8)[t]}
\put(8,5){\oval(11.7,3.2)[t]}
\put(8,5){\oval(4,3)[b]}
\put(6,1){\small $D_{14}^-$}
\end{picture} \quad
\begin{picture}(16,8)
\put(0,5){\line(1,0){16}}
\put(2,5){\circle*{0.5}}
\put(6,5){\circle*{0.5}}
\put(10,5){\circle*{0.5}}
\put(14,5){\circle*{0.5}}
\put(7,7.5){\small $+$}
\put(8,5){\oval(4.3,3.8)[t]}
\put(8,5){\oval(3.7,3.2)[t]}
\put(8,5){\oval(12,3)[b]}
\put(6,1){\small $D_{23}^+$}
\end{picture} \quad
\begin{picture}(16,8)
\put(0,5){\line(1,0){16}}
\put(2,5){\circle*{0.5}}
\put(6,5){\circle*{0.5}}
\put(10,5){\circle*{0.5}}
\put(14,5){\circle*{0.5}}
\put(7,7){\small $-$}
\put(8,5){\oval(4.3,3.8)[t]}
\put(8,5){\oval(3.7,3.2)[t]}
\put(8,5){\oval(12,3)[b]}
\put(6,1){\small $D_{23}^-$}
\end{picture}
\caption{$D$}
\label{D}
\end{figure}

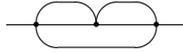
\begin{figure}[h]
\unitlength 2mm
\begin{picture}(12,5)
\put(0,3){\line(1,0){12}}
\put(2,3){\circle*{0.5}}
\put(6,3){\circle*{0.5}}
\put(10,3){\circle*{0.5}}
\put(4,3){\oval(4,3)[t]}
\put(8,3){\oval(4,3)[t]}
\put(6,3){\oval(8,3)[b]}
\end{picture}
\caption{$E_+$, $E_1$, $E_2$, $E_3$}
\label{E}
\end{figure}

The definition of the cells of types $B_i^j$, $C_{i j}$, and $D_{i j}^\pm$ is clear from their pictures, see Figs.~\ref{B}, \ref{C}, and \ref{D}. In particular,
each algebra of class $D$ depends on a parameter $\alpha \neq 0$ and four points $A, B, C, D \in (0, 2\pi)$; it consists of all functions $f$ such that $f(A)=f(B), f(C)=f(D)$, and $f'(A) = \alpha f'(B)$. These algebras form the cells $D^+_{i j}$ and $D^-_{i j}$, where the index $\pm$ is the sign of $\alpha$, and the indices $i<j$ are the numbers of the points $A$ and $B$ among all four points of the support.

Algebras of class $E$ are a special case of the algebras $\divideontimes$ defined in \S \ref{dvd}. Every such algebra
depends on one point $(\alpha: \beta: \gamma) \in {\mathbb R}P^2,$ $\alpha \cdot \beta \cdot \gamma \neq 0$, and on three points $A < B < C \in (0, 2\pi)$. Namely, the algebra $E((\alpha:\beta:\gamma); A, B, C)$
 consists of all functions $f$ such that $f(A) = f (B) = f(C)$ and $\alpha f'(A) + \beta f'(B) + \gamma f'(C)=0$. Such an algebra belongs to the 
cell $E_+$ if all the coefficients $\alpha, \beta,$ and $ \gamma$ have the same sign. It belongs to the cell $E_1$ (respectively, $E_2,$ respectively, $E_3$) if the sign of $\alpha$ (respectively, $\beta$, respectively, $\gamma$) differs from the signs of the other two coefficients.

\FloatBarrier

\subsection{Four-dimensional cells}
\label{fourdim}

There are 44 four-dimensional cells of the first type in $\overline{CD}_3(S^1)$, see Figs.~\ref{a}--\ref{k}. The definitions of the subalgebras of classes $a, b, c, d, e,$ and $ g,$ which form the corresponding cells, are clear from the pictures.

\unitlength 2mm
\begin{figure}[h]
\begin{picture}(16,7)
\put(0,4){\line(1,0){16}}
\put(2,4){\makebox(0,0)[cc]{$\ast$}}
\put(6,4){\makebox(0,0)[cc]{$\ast$}}
\put(10,4){\circle*{0.5}}
\put(14,4){\circle*{0.5}}
\put(12,4){\oval(4,3)[t]}
\put(7,1){$a_{12}$}
\end{picture} \qquad
\begin{picture}(16,6)
\put(0,4){\line(1,0){16}}
\put(2,4){\makebox(0,0)[cc]{$\ast$}}
\put(10,4){\makebox(0,0)[cc]{$\ast$}}
\put(6,4){\circle*{0.5}}
\put(14,4){\circle*{0.5}}
\put(10,4){\oval(8,3)[t]}
\put(7,1){$a_{13}$}
\end{picture} \qquad
\begin{picture}(16,6)
\put(0,4){\line(1,0){16}}
\put(2,4){\makebox(0,0)[cc]{$\ast$}}
\put(14,4){\makebox(0,0)[cc]{$\ast$}}
\put(10,4){\circle*{0.5}}
\put(6,4){\circle*{0.5}}
\put(8,4){\oval(4,3)[t]}
\put(7,1){$a_{14}$}
\end{picture}

\begin{picture}(16,7)
\put(0,4){\line(1,0){16}}
\put(10,4){\makebox(0,0)[cc]{$\ast$}}
\put(6,4){\makebox(0,0)[cc]{$\ast$}}
\put(2,4){\circle*{0.5}}
\put(14,4){\circle*{0.5}}
\put(8,4){\oval(12,3)[t]}
\put(7,1){$a_{23}$}
\end{picture} \qquad
\begin{picture}(16,6)
\put(0,4){\line(1,0){16}}
\put(14,4){\makebox(0,0)[cc]{$\ast$}}
\put(6,4){\makebox(0,0)[cc]{$\ast$}}
\put(10,4){\circle*{0.5}}
\put(2,4){\circle*{0.5}}
\put(6,4){\oval(8,3)[t]}
\put(7,1){$a_{24}$}
\end{picture} \qquad
\begin{picture}(16,6)
\put(0,4){\line(1,0){16}}
\put(10,4){\makebox(0,0)[cc]{$\ast$}}
\put(14,4){\makebox(0,0)[cc]{$\ast$}}
\put(6,4){\circle*{0.5}}
\put(2,4){\circle*{0.5}}
\put(4,4){\oval(4,3)[t]}
\put(7,1){$a_{34}$}
\end{picture}
\caption{$a$}
\label{a}
\end{figure}

\begin{figure}[h]
\begin{picture}(16,8)
\put(0,5){\line(1,0){16}}
\put(2,5){\makebox(0,0)[cc]{$\ast$}}
\put(6,5){\circle*{0.5}}
\put(10,5){\circle*{0.5}}
\put(14,5){\circle*{0.5}}
\put(12,5){\oval(4,3)[t]}
\put(4,5){\oval(4,3)[t]}
\put(6,2){$b_{12}$}
\end{picture} \qquad
\begin{picture}(16,7)
\put(0,5){\line(1,0){16}}
\put(2,5){\makebox(0,0)[cc]{$\ast$}}
\put(6,5){\circle*{0.5}}
\put(10,5){\circle*{0.5}}
\put(14,5){\circle*{0.5}}
\put(10,5){\oval(8,3)[b]}
\put(6,5){\oval(8,3)[t]}
\put(6,2){$b_{13}$}
\end{picture} \qquad
\begin{picture}(16,7)
\put(0,5){\line(1,0){16}}
\put(2,5){\makebox(0,0)[cc]{$\ast$}}
\put(6,5){\circle*{0.5}}
\put(10,5){\circle*{0.5}}
\put(14,5){\circle*{0.5}}
\put(8,5){\oval(12,3)[t]}
\put(8,5){\oval(4,3)[b]}
\put(6,2){$b_{14}$}
\end{picture}

\begin{picture}(16,7)
\put(0,5){\line(1,0){16}}
\put(6,5){\makebox(0,0)[cc]{$\ast$}}
\put(2,5){\circle*{0.5}}
\put(10,5){\circle*{0.5}}
\put(14,5){\circle*{0.5}}
\put(12,5){\oval(4,3)[t]}
\put(4,5){\oval(4,3)[t]}
\put(6,2){$b_{21}$}
\end{picture} \qquad
\begin{picture}(16,7)
\put(0,5){\line(1,0){16}}
\put(6,5){\makebox(0,0)[cc]{$\ast$}}
\put(2,5){\circle*{0.5}}
\put(10,5){\circle*{0.5}}
\put(14,5){\circle*{0.5}}
\put(8,5){\oval(4,3)[t]}
\put(8,5){\oval(12,3)[b]}
\put(6,2){$b_{23}$}
\end{picture} \qquad
\begin{picture}(16,7)
\put(0,5){\line(1,0){16}}
\put(6,5){\makebox(0,0)[cc]{$\ast$}}
\put(2,5){\circle*{0.5}}
\put(10,5){\circle*{0.5}}
\put(14,5){\circle*{0.5}}
\put(10,5){\oval(8,3)[t]}
\put(6,5){\oval(8,3)[b]}
\put(6,2){$b_{24}$}
\end{picture}

\begin{picture}(16,7)
\put(0,5){\line(1,0){16}}
\put(10,5){\makebox(0,0)[cc]{$\ast$}}
\put(6,5){\circle*{0.5}}
\put(2,5){\circle*{0.5}}
\put(14,5){\circle*{0.5}}
\put(6,5){\oval(8,3)[t]}
\put(10,5){\oval(8,3)[b]}
\put(6,2){$b_{31}$}
\end{picture} \qquad
\begin{picture}(16,7)
\put(0,5){\line(1,0){16}}
\put(10,5){\makebox(0,0)[cc]{$\ast$}}
\put(6,5){\circle*{0.5}}
\put(2,5){\circle*{0.5}}
\put(14,5){\circle*{0.5}}
\put(8,5){\oval(4,3)[t]}
\put(8,5){\oval(12,3)[b]}
\put(6,2){$b_{32}$}
\end{picture} \qquad
\begin{picture}(16,7)
\put(0,5){\line(1,0){16}}
\put(10,5){\makebox(0,0)[cc]{$\ast$}}
\put(6,5){\circle*{0.5}}
\put(2,5){\circle*{0.5}}
\put(14,5){\circle*{0.5}}
\put(12,5){\oval(4,3)[t]}
\put(4,5){\oval(4,3)[t]}
\put(6,2){$b_{34}$}
\end{picture}

\begin{picture}(16,7)
\put(0,5){\line(1,0){16}}
\put(14,5){\makebox(0,0)[cc]{$\ast$}}
\put(6,5){\circle*{0.5}}
\put(10,5){\circle*{0.5}}
\put(2,5){\circle*{0.5}}
\put(8,5){\oval(12,3)[t]}
\put(8,5){\oval(4,3)[b]}
\put(6,2){$b_{41}$}
\end{picture} \qquad
\begin{picture}(16,7)
\put(0,5){\line(1,0){16}}
\put(14,5){\makebox(0,0)[cc]{$\ast$}}
\put(6,5){\circle*{0.5}}
\put(10,5){\circle*{0.5}}
\put(2,5){\circle*{0.5}}
\put(10,5){\oval(8,3)[t]}
\put(6,5){\oval(8,3)[b]}
\put(6,2){$b_{42}$}
\end{picture} \qquad
\begin{picture}(16,7)
\put(0,5){\line(1,0){16}}
\put(14,5){\makebox(0,0)[cc]{$\ast$}}
\put(6,5){\circle*{0.5}}
\put(10,5){\circle*{0.5}}
\put(2,5){\circle*{0.5}}
\put(12,5){\oval(4,3)[t]}
\put(4,5){\oval(4,3)[t]}
\put(6,2){$b_{43}$}
\end{picture}
\caption{$b$}
\label{b}
\end{figure}

{\unitlength 1.6mm
\begin{figure}[h]
\begin{picture}(16,8)
\put(0,4){\line(1,0){16}}
\put(2,4){\makebox(0,0)[cc]{$\ast$}}
\put(6,4){\circle*{0.5}}
\put(10,4){\circle*{0.5}}
\put(14,4){\circle*{0.5}}
\put(8,8){\line(-1,-2){2}}
\put(8,8){\line(1,-2){2}}
\put(8,8){\line(3,-2){6}}
\put(7,0){$c_1$}
\end{picture} \quad
\begin{picture}(16,8)
\put(0,4){\line(1,0){16}}
\put(6,4){\makebox(0,0)[cc]{$\ast$}}
\put(14,4){\circle*{0.5}}
\put(10,4){\circle*{0.5}}
\put(2,4){\circle*{0.5}}
\put(8,8){\line(-3,-2){6}}
\put(8,8){\line(1,-2){2}}
\put(8,8){\line(3,-2){6}}
\put(7,0){$c_2$}
\end{picture} \quad
\begin{picture}(16,8)
\put(0,4){\line(1,0){16}}
\put(10,4){\makebox(0,0)[cc]{$\ast$}}
\put(6,4){\circle*{0.5}}
\put(14,4){\circle*{0.5}}
\put(2,4){\circle*{0.5}}
\put(8,8){\line(-3,-2){6}}
\put(8,8){\line(-1,-2){2}}
\put(8,8){\line(3,-2){6}}
\put(7,0){$c_3$}
\end{picture} \quad
\begin{picture}(16,8)
\put(0,4){\line(1,0){16}}
\put(14,4){\makebox(0,0)[cc]{$\ast$}}
\put(6,4){\circle*{0.5}}
\put(10,4){\circle*{0.5}}
\put(2,4){\circle*{0.5}}
\put(8,8){\line(-3,-2){6}}
\put(8,8){\line(-1,-2){2}}
\put(8,8){\line(1,-2){2}}
\put(7,0){$c_4$}
\end{picture}
\caption{$c$}
\label{c}
\end{figure}
}

\begin{figure}[h]
\begin{picture}(16,6)
\put(0,2){\line(1,0){16}}
\put(14,2){\circle*{0.5}}
\put(6,2){\circle*{0.5}}
\put(10,2){\circle*{0.5}}
\put(2,2){\circle*{0.5}}
\put(8,6){\line(-3,-2){6}}
\put(8,6){\line(-1,-2){2}}
\put(8,6){\line(1,-2){2}}
\put(8,6){\line(3,-2){6}}
\put(7,0){$d$}
\end{picture}
\caption{$d$}
\label{d}
\end{figure}

\begin{figure}[h]
\begin{picture}(12,8)
\put(0,4){\line(1,0){12}}
\put(8,6.5){\makebox(0,0)[cc]{\small $+$}}
\put(6,4){\circle*{0.5}}
\put(10,4){\circle*{0.5}}
\put(2,4){\makebox(0,0)[cc]{$\ast$}}
\put(8,4){\oval(4.3,3.2)[t]}
\put(8,4){\oval(3.7,2.6)[t]}
\put(5,1.5){$e_1^+$}
\end{picture} \qquad \quad
\begin{picture}(12,7)
\put(0,4){\line(1,0){12}}
\put(6,6.5){\makebox(0,0)[cc]{\small $+$}}
\put(2,4){\circle*{0.5}}
\put(10,4){\circle*{0.5}}
\put(6,4){\makebox(0,0)[cc]{$\ast$}}
\put(6,4){\oval(8.3,3.2)[t]}
\put(6,4){\oval(7.7,2.6)[t]}
\put(5,1.5){$e_2^+$}
\end{picture} \qquad \quad
\begin{picture}(12,7)
\put(0,4){\line(1,0){12}}
\put(4,6.5){\makebox(0,0)[cc]{\small $+$}}
\put(2,4){\circle*{0.5}}
\put(6,4){\circle*{0.5}}
\put(10,4){\makebox(0,0)[cc]{$\ast$}}
\put(4,4){\oval(4.3,3.2)[t]}
\put(4,4){\oval(3.7,2.6)[t]}
\put(5,1.5){$e_3^+$}
\end{picture}

\begin{picture}(12,8)
\put(0,4){\line(1,0){12}}
\put(8,6.2){\makebox(0,0)[cc]{\small $-$}}
\put(6,4){\circle*{0.5}}
\put(10,4){\circle*{0.5}}
\put(2,4){\makebox(0,0)[cc]{$\ast$}}
\put(8,4){\oval(4.3,3.2)[t]}
\put(8,4){\oval(3.7,2.6)[t]}
\put(5,1.5){$e_1^-$}
\end{picture} \qquad \quad
\begin{picture}(12,7)
\put(0,4){\line(1,0){12}}
\put(6,6.2){\makebox(0,0)[cc]{\small $-$}}
\put(2,4){\circle*{0.5}}
\put(10,4){\circle*{0.5}}
\put(6,4){\makebox(0,0)[cc]{$\ast$}}
\put(6,4){\oval(8.3,3.2)[t]}
\put(6,4){\oval(7.7,2.6)[t]}
\put(5,1.5){$e_2^-$}
\end{picture} \qquad \quad
\begin{picture}(12,7)
\put(0,4){\line(1,0){12}}
\put(4,6.2){\makebox(0,0)[cc]{\small $-$}}
\put(2,4){\circle*{0.5}}
\put(6,4){\circle*{0.5}}
\put(10,4){\makebox(0,0)[cc]{$\ast$}}
\put(4,4){\oval(4.3,3.2)[t]}
\put(4,4){\oval(3.7,2.6)[t]}
\put(5,1.5){$e_3^-$}
\end{picture}
\caption{$e$}
\label{e}
\end{figure}

{\unitlength 1.7mm
\begin{figure}[h]
\begin{picture}(12,7)
\put(0,4){\line(1,0){12}}
\put(2,4){\makebox(0,0)[cc]{$\circledast$}}
\put(6,4){\circle*{0.5}}
\put(10,4){\circle*{0.5}}
\put(8,4){\oval(4,3)[t]}
\put(5,1){$g_1$}
\end{picture} \quad 
\begin{picture}(12,7)
\put(0,4){\line(1,0){12}}
\put(6,4){\makebox(0,0)[cc]{$\circledast$}}
\put(2,4){\circle*{0.5}}
\put(10,4){\circle*{0.5}}
\put(6,4){\oval(8,3)[t]}
\put(5,1){$g_2$}
\end{picture} \quad 
\begin{picture}(12,7)
\put(0,4){\line(1,0){12}}
\put(10,4){\makebox(0,0)[cc]{$\circledast$}}
\put(2,4){\circle*{0.5}}
\put(6,4){\circle*{0.5}}
\put(4,4){\oval(4,3)[t]}
\put(5,1){$g_3$}
\end{picture}
\qquad \quad
\begin{picture}(10,7)
\put(0,4){\line(1,0){10}}
\put(5,7.5){\makebox(0,0)[cc]{\small $+$}}
\put(2,4){\circle*{0.5}}
\put(8,4){\circle*{0.5}}
\put(5,4){\oval(6.6,4.1)[t]}
\put(5,4){\oval(6,3.5)[t]}
\put(5,4){\oval(5.4,2.9)[t]}
\put(5,0){$h^+$}
\end{picture} \quad 
\begin{picture}(10,7)
\put(0,4){\line(1,0){10}}
\put(5,7){\makebox(0,0)[cc]{\small $-$}}
\put(2,4){\circle*{0.5}}
\put(8,4){\circle*{0.5}}
\put(5,4){\oval(6.6,4.1)[t]}
\put(5,4){\oval(6,3.5)[t]}
\put(5,4){\oval(5.4,2.9)[t]}
\put(5,0){$h^-$}
\end{picture}
\caption{$g$ and $h$}
\label{g h}
\end{figure}
}

The cells $h^+$ and $h^-$ are formed by the algebras $\gimel_3(\Delta; A,B)$, see \S \ref{twoser}. They depend on two parameters, $\alpha \in {\mathbb R}\setminus 0$ and $\beta \in {\mathbb R}$, and on two points, $A, B \in S^1$.
The algebra $h (\alpha, \beta; A, B)$ consists of all functions $f$ such that \begin{equation}
\label{eqh}
f(A)=f(B), \quad f'(A) = \alpha f'(B), \quad f''(A) = \beta f'(B) + \alpha^2 f''(B).
\end{equation}
 These algebras satisfy the identity
\begin{equation}
\label{idenh}
h(\alpha, \beta; A, B) \equiv h (\alpha^{-1}, - \beta \alpha^{-3}; B, A).
\end{equation}
They fill in two cells, $h^+$ and $h^-$, which are characterized by the sign of the parameter $\alpha$. 

The algebras of class $i$ depend on two points $A, B \in S^1$ and one point $(\alpha : \beta : \gamma) \in {\mathbb R}P^2$, $\alpha \neq 0 \neq \gamma$.
The algebra $i((\alpha : \beta : \gamma); A, B)$ consists of all functions $f$ such that 
\begin{equation}
\label{ialg}
f(A) = f(B), \quad f'(A)=0, \quad \alpha f'''(A) + \beta f''(A) + \gamma f'(B) =0.
\end{equation}
Each such algebra of the first type belongs to the cell $i_1^\pm$ or $i_2^\pm$, where the 
lower index is $1$ if $A<B$ and is $2$ if $A> B$; the upper index $\pm$ is the sign of the ratio $\alpha / \gamma$. Near the typical singular points of classes $i^+$ (respectively, $i^-$) the singular knots in ${\mathbb R}^3$ look as  
{\unitlength 0.6 mm
\begin{picture}(16,10)
\bezier{100}(8,0)(8,4)(16,8)
\bezier{100}(8,0)(8,4)(0,8)
\put(16,2){\vector(-4,-1){16}}
\put(16,7.7){\vector(4,1){1}}
\end{picture}
} (respectively, 
{\unitlength 0.6 mm
\begin{picture}(16,10)
\bezier{100}(8,0)(8,4)(16,8)
\bezier{100}(8,0)(8,4)(0,8)
\put(0,-2){\vector(4,1){16}}
\put(16,7.7){\vector(4,1){1}}
\end{picture}
}) where the derivative of the smooth branch and the second and third derivatives of the cuspidal branch lie in the same plane.

The algebras of class $k$ depend on one real parameter $\alpha \neq 0$ 
and three ordered points $A, B, C$ of the support. Each such algebra
$k(\alpha; A, B, C)$ consists of functions that satisfy the conditions $f(A)=f(B)=f(C)$ and $f'(A) = \alpha f'(B)$. Obviously, $k(\alpha; A, B, C) = k(\alpha^{-1}; B, A, C)$.
These algebras form six cells $k_i^+$ and $k_i^-$, $i=1, 2,$ or $3$.

Additionally, there are 41 four-dimensional cells $\bar B_i^j$, $\bar C_{i j}$, $\bar D_{i j}^\pm$, $\bar E_\star$ of the second type, corresponding to the five-dimensional cells described in Figs.~\ref{B}--\ref{E}.

\begin{figure}[h] 
\begin{picture}(12,6)
\put(0,3){\line(1,0){12}}
\put(2,3){\makebox(0,0)[cc]{$\circledast$}}
\put(10,3){\circle*{0.5}}
\put(6,3){\oval(7.7,3.2)[t]}
\put(6,3){\oval(8.3, 3.8)[t]}
\put(5,0){$i_1^+$}
\put(6,6){\makebox(0,0)[cc]{\small $+$}}
\end{picture} \qquad 
\begin{picture}(12,6)
\put(0,3){\line(1,0){12}}
\put(2,3){\makebox(0,0)[cc]{$\circledast$}}
\put(10,3){\circle*{0.5}}
\put(6,3){\oval(7.7,3.2)[t]}
\put(6,3){\oval(8.3, 3.8)[t]}
\put(5,0){$i_1^-$}
\put(6,5.5){\makebox(0,0)[cc]{\small $-$}}
\end{picture} \qquad 
\begin{picture}(12,6)
\put(0,3){\line(1,0){12}}
\put(10,3){\makebox(0,0)[cc]{$\circledast$}}
\put(2,3){\circle*{0.5}}
\put(6,3){\oval(7.7,3.2)[t]}
\put(6,3){\oval(8.3, 3.8)[t]}
\put(5,0){$i_2^+$}
\put(6,6){\makebox(0,0)[cc]{\small $+$}}
\end{picture} \qquad
\begin{picture}(12,6)
\put(0,3){\line(1,0){12}}
\put(10,3){\makebox(0,0)[cc]{$\circledast$}}
\put(2,3){\circle*{0.5}}
\put(6,3){\oval(7.7,3.2)[t]}
\put(6,3){\oval(8.3, 3.8)[t]}
\put(5,0){$i_2^-$}
\put(6,5.5){\makebox(0,0)[cc]{\small $-$}}
\end{picture} 
\caption{$i$ }
\label{i}
\end{figure}

\begin{figure}[h]
\begin{picture}(12,10)
\put(0,5){\line(1,0){12}}
\put(9,8){\makebox(0,0)[cc]{\small $+$}}
\put(6,5){\circle*{0.5}}
\put(10,5){\circle*{0.5}}
\put(2,5){\circle*{0.5}}
\put(8,5){\oval(4.3,3.2)[t]}
\put(8,5){\oval(3.7,2.6)[t]}
\bezier{200}(2,5)(7,10)(8,6.6)
\put(4,2){$k_1^+$}
\end{picture} \qquad \quad
\begin{picture}(12,10)
\put(0,5){\line(1,0){12}}
\put(6,8){\makebox(0,0)[cc]{\small $+$}}
\put(2,5){\circle*{0.5}}
\put(10,5){\circle*{0.5}}
\put(6,5){\circle*{0.5}}
\put(6,5){\oval(8.3,3.2)[t]}
\put(6,5){\oval(7.7,2.6)[t]}
\bezier{50}(6,5)(6,5.7)(6,6.3)
\put(4,2){$k_2^+$}
\end{picture} \qquad \quad
\begin{picture}(12,10)
\put(0,5){\line(1,0){12}}
\put(3,8){\makebox(0,0)[cc]{\small $+$}}
\put(2,5){\circle*{0.5}}
\put(6,5){\circle*{0.5}}
\put(10,5){\circle*{0.5}}
\put(4,5){\oval(4.3,3.2)[t]}
\put(4,5){\oval(3.7,2.6)[t]}
\bezier{200}(10,5)(5,10)(4,6.6)
\put(4,2){$k_3^+$}
\end{picture}

\begin{picture}(12,9)
\put(0,5){\line(1,0){12}}
\put(9,7.5){\makebox(0,0)[cc]{\small $-$}}
\put(6,5){\circle*{0.5}}
\put(10,5){\circle*{0.5}}
\put(2,5){\circle*{0.5}}
\put(8,5){\oval(4.3,3.2)[t]}
\put(8,5){\oval(3.7,2.6)[t]}
\bezier{200}(2,5)(7,10)(8,6.6)
\put(4,2){$k_1^-$}
\end{picture} \qquad \quad
\begin{picture}(12,8)
\put(0,5){\line(1,0){12}}
\put(6,7.5){\makebox(0,0)[cc]{\small $-$}}
\put(2,5){\circle*{0.5}}
\put(10,5){\circle*{0.5}}
\put(6,5){\circle*{0.5}}
\put(6,5){\oval(8.3,3.2)[t]}
\put(6,5){\oval(7.7,2.6)[t]}
\bezier{50}(6,5)(6,5.7)(6,6.3)
\put(4,2){$k_2^-$}
\end{picture} \qquad \quad
\begin{picture}(12,8)
\put(0,5){\line(1,0){12}}
\put(3,7.5){\makebox(0,0)[cc]{\small $-$}}
\put(2,5){\circle*{0.5}}
\put(6,5){\circle*{0.5}}
\put(10,5){\circle*{0.5}}
\put(4,5){\oval(4.3,3.2)[t]}
\put(4,5){\oval(3.7,2.6)[t]}
\bezier{200}(10,5)(5,10)(4,6.6)
\put(4,2){$k_3^-$}
\end{picture}
\caption{$k$}
\label{k}
\end{figure}

\FloatBarrier

\subsection{Three-dimensional cells}
\label{threedim}

{
\unitlength 1.6mm

\begin{figure}[h]
\begin{picture}(10,8)
\put(0,5){\line(1,0){10}}
\put(2,5){\makebox(0,0)[cc]{$\ast$}}
\put(8,5){\makebox(0,0)[cc]{$\ast$}}
\put(5,5){\circle*{0.5}}
\put(3.5,5){\oval(3,3)[t]}
\put(4,0){\small $U_{12}$}
\end{picture} \quad
\begin{picture}(10,8)
\put(0,5){\line(1,0){10}}
\put(2,5){\makebox(0,0)[cc]{$\ast$}}
\put(5,5){\makebox(0,0)[cc]{$\ast$}}
\put(8,5){\circle*{0.5}}
\put(5,5){\oval(6,3)[t]}
\put(4,0){\small $U_{13}$}
\end{picture} \quad
\begin{picture}(10,8)
\put(0,5){\line(1,0){10}}
\put(5,5){\makebox(0,0)[cc]{$\ast$}}
\put(2,5){\makebox(0,0)[cc]{$\ast$}}
\put(8,5){\circle*{0.5}}
\put(6.5,5){\oval(3,3)[t]}
\put(4,0){\small $U_{23}$}
\end{picture} \qquad
\begin{picture}(10,8)
\put(0,5){\line(1,0){10}}
\put(5,5){\makebox(0,0)[cc]{$\ast$}}
\put(8,5){\makebox(0,0)[cc]{$\ast$}}
\put(2,5){\circle*{0.5}}
\put(3.5,5){\oval(3,3)[t]}
\put(4,0){\small $U_{21}$}
\end{picture} \quad
\begin{picture}(10,8)
\put(0,5){\line(1,0){10}}
\put(8,5){\makebox(0,0)[cc]{$\ast$}}
\put(5,5){\makebox(0,0)[cc]{$\ast$}}
\put(2,5){\circle*{0.5}}
\put(5,5){\oval(6,3)[t]}
\put(4,0){\small $U_{31}$}
\end{picture} \quad
\begin{picture}(10,8)
\put(0,5){\line(1,0){10}}
\put(8,5){\makebox(0,0)[cc]{$\ast$}}
\put(2,5){\makebox(0,0)[cc]{$\ast$}}
\put(5,5){\circle*{0.5}}
\put(6.5,5){\oval(3,3)[t]}
\put(4,0){\small $U_{32}$}
\end{picture}
\caption{$U$} 
\label{U}
\end{figure}
}

{
\unitlength 1.8mm

\begin{figure}[h]
\begin{picture}(12,8)
\put(0,5){\line(1,0){12}}
\put(2,4.5){\makebox(0,0)[cc]{$\ast$}}
\put(2,5.3){\makebox(0,0)[cc]{$\ast$}}
\put(6,5){\circle*{0.5}}
\put(10,5){\circle*{0.5}}
\put(8,5){\oval(4,3)[t]}
\put(5,2){\small $V_{1}$}
\end{picture}
 \qquad \quad
\begin{picture}(12,8)
\put(0,5){\line(1,0){12}}
\put(6,4.5){\makebox(0,0)[cc]{$\ast$}}
\put(6,5.3){\makebox(0,0)[cc]{$\ast$}}
\put(2,5){\circle*{0.5}}
\put(10,5){\circle*{0.5}}
\put(6,5){\oval(8,3)[t]}
\put(5,2){\small $V_{2}$}
\end{picture} \qquad \quad
\begin{picture}(12,8)
\put(0,5){\line(1,0){12}}
\put(10,4.5){\makebox(0,0)[cc]{$\ast$}}
\put(10,5.3){\makebox(0,0)[cc]{$\ast$}}
\put(2,5){\circle*{0.5}}
\put(6,5){\circle*{0.5}}
\put(4,5){\oval(4,3)[t]}
\put(5,2){\small $V_{3}$}
\end{picture}

\begin{picture}(12,7)
\put(0,4){\line(1,0){12}}
\put(2,4){\makebox(0,0)[cc]{$\ast$}}
\put(6,7){\line(-4,-3){4}}
\put(6,7){\line(0,-1){3}}
\put(6,7){\line(4,-3){4}}
\put(5,1){\small $W_1$}
\end{picture} \qquad
\begin{picture}(12,7)
\put(0,4){\line(1,0){12}}
\put(6,4){\makebox(0,0)[cc]{$\ast$}}
\put(6,7){\line(-4,-3){4}}
\put(6,7){\line(0,-1){3}}
\put(6,7){\line(4,-3){4}}
\put(5,1){\small $W_2$}
\end{picture} \qquad 
\begin{picture}(12,7)
\put(0,4){\line(1,0){12}}
\put(10,4){\makebox(0,0)[cc]{$\ast$}}
\put(6,7){\line(-4,-3){4}}
\put(6,7){\line(0,-1){3}}
\put(6,7){\line(4,-3){4}}
\put(5,1){\small $W_3$}
\end{picture} 
\qquad \quad
\begin{picture}(12,7)
\put(0,4){\line(1,0){12}}
\put(2,4){\makebox(0,0)[cc]{$\ast$}}
\put(6,4){\makebox(0,0)[cc]{$\ast$}}
\put(10,4){\makebox(0,0)[cc]{$\ast$}}
\put(5,1){\small $X$}
\end{picture}
\caption{$V$, $W$ and $X$}
\label{VWX}
\end{figure}
}

There are 23 three-dimensional cells of the first type. They are shown in Figs.~\ref{U}--\ref{TJ}. The definitions of the algebras of classes $U, V, W, X, $ and $Y$ are clear from their pictures.

Every algebra $Z (\alpha; A, B)$ depends on a number $\alpha \in {\mathbb R} \setminus \{0\}$ and two points $A, B \in (0, 2\pi)$. It consists of all functions $f$ such that $f(A) = f(B)$, $f'(A)=0$, and $f''(A) = \alpha f'(B)$. These algebras form four cells $Z^{\pm}_{1,2}$, where the upper index $+$ or $-$ is the sign of $\alpha$, and the lower index is 1 if $A<B$ in $(0, 2\pi)$ and is $2$ if $A>B$.

{
\unitlength 1.7mm
\begin{figure}[h]
\begin{picture}(9,7)
\put(0,4){\line(1,0){9}}
\put(2,4){\makebox(0,0)[cc]{$\circledast$}}
\put(7,4){\makebox(0,0)[cc]{$\ast$}}
\put(4,0){$Y_{1}$}
\end{picture} \quad
\begin{picture}(9,7)
\put(0,4){\line(1,0){9}}
\put(2,4){\makebox(0,0)[cc]{$\ast$}}
\put(7,4){\makebox(0,0)[cc]{$\circledast$}}
\put(4,0){\small $Y_{2}$}
\end{picture} 
\qquad 
\begin{picture}(10,10)
\put(0,4){\line(1,0){10}}
\put(2,4){\makebox(0,0)[cc]{$\ast$}}
\put(5,4){\oval(5.7,2.7)[t]}
\put(5,4){\oval(6.3,3.3)[t]}
\put(8,4){\circle*{0.5}}
\put(4,6.5){\small $+$}
\put(4,0){\small $Z_1^+$}
\end{picture} \quad
\begin{picture}(10,8)
\put(0,4){\line(1,0){10}}
\put(8,4){\makebox(0,0)[cc]{$\ast$}}
\put(5,4){\oval(5.7,2.7)[t]}
\put(5,4){\oval(6.3,3.3)[t]}
\put(2,4){\circle*{0.5}}
\put(4,6.5){\small $+$}
\put(4,0){\small $Z_2^+$}
\end{picture} \quad
\begin{picture}(10,10)
\put(0,4){\line(1,0){10}}
\put(2,4){\makebox(0,0)[cc]{$\ast$}}
\put(5,4){\oval(5.7,2.7)[t]}
\put(5,4){\oval(6.3,3.3)[t]}
\put(8,4){\circle*{0.5}}
\put(4,6.5){\small $-$}
\put(4,0){\small $Z_1^-$}
\end{picture} \quad
\begin{picture}(10,8)
\put(0,4){\line(1,0){10}}
\put(8,4){\makebox(0,0)[cc]{$\ast$}}
\put(5,4){\oval(5.7,2.7)[t]}
\put(5,4){\oval(6.3,3.3)[t]}
\put(2,4){\circle*{0.5}}
\put(4,6.5){\small $-$}
\put(4,0){\small $Z_2^-$}
\end{picture}
\caption{$Y$ and $Z$}
\label{YZ}
\end{figure}

The algebra $T(\alpha; A, B)$ is defined for any $\alpha \in {\mathbb R}$ and any two points $A \neq B $ in $(0, 2\pi)$. It consists of all functions $f$ such that $f(A) = f(B)$, $f'(A) = 0 ,$ and $f'''(A)= \alpha f'' (A)$. These algebras occupy two cells $T_1$ and $T_2$, which consist of all algebras $T(\alpha; A, B)$ with $A <B$ (respectively, $A>B$).

\begin{figure}[h]
\begin{picture}(10,8)
\put(0,4){\line(1,0){10}}
\put(2,4){\makebox(0,0)[cc]{$\circledast$}}
\put(5,4){\oval(6,3)[t]}
\put(8,4){\circle*{0.5}}
\put(4,0){\small $T_1$}
\end{picture} \qquad 
\begin{picture}(10,8)
\put(0,4){\line(1,0){10}}
\put(8,4){\makebox(0,0)[cc]{$\circledast$}}
\put(5,4){\oval(6,3)[t]}
\put(2,4){\circle*{0.5}}
\put(4,0){\small $T_2$}
\end{picture} 
\qquad \qquad
\begin{picture}(8,7)
\put(0,4){\line(1,0){8}}
\put(4,4){\makebox(0,0)[cc]{$\circledast$}}
\put(4,4){\circle{2.8}}
\put(0,0){\small $J$}
\end{picture} \qquad \qquad
\begin{picture}(10,7)
\put(0,4){\line(1,0){10}}
\put(5,3.3){\makebox(0,0)[cc]{$\ast$}}
\put(5,4.5){\makebox(0,0)[cc]{$\ast$}}
\put(5,4){\oval(1.7,2.7)}
\put(5,4){\oval(2.2,3.2)}
\put(0,0){\small $S$}
\end{picture}
\caption{$T$, $J$ and $S$}
\label{TJ}
\end{figure}
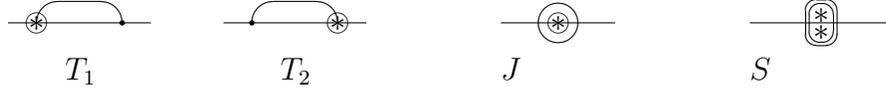
}
Each algebra $J(\alpha, \beta; A)$ depends on two real parameters $\alpha$ and $\beta$ and on its unique support point $A \in S^1$. It 
consists of all functions $f$ such that \begin{equation}
f'(A)=0, \quad f'''(A) = 3 \alpha f''(A), \quad f^{\mbox{\scriptsize V}}(A) = 10 \alpha f^{\mbox{\footnotesize IV}}(A) + \beta f''(A) =0.
\label{jten}
\end{equation}
These algebras of the first type (i.e., with $A \neq \bullet)$ form one cell $J$.

Each algebra $S(\alpha, \beta; A)$ has support at a single point $A$, depends on two parameters $\alpha, \beta$, and is defined by the conditions \begin{equation}
\label{S}
f'(A)=0, \qquad f''(A)=0, \qquad f^{V}(A) = \alpha f^{IV}(A) + \beta f'''(A).\end{equation} 
These algebras of the first type fill in one three-dimensional cell $S$ with the parameters $\alpha, \beta$ and $A$. 

Additionally, there are 44 three-dimensional cells $\bar a_{12}, \dots, \bar k_3^-$ of the second type corresponding to the algebras shown in Figs. \ref{a}--\ref{k}.

\subsection{Two-dimensional cells}
\label{twodim}

There are six two-dimensional cells of the first type, shown in Fig. \ref{Omega}. The cell $\Theta$ consists of
subalgebras $\Theta(\alpha; A)$, each of which has support at a single point $A$, depends on the parameter $\alpha \in {\mathbb R}^1$, and consists of functions $f$ such that $$f'(A)=0, \qquad f''(A)=0, \qquad f^{IV}(A) = \alpha f'''(A).$$ 

Additionally, there are 23 two-dimensional cells of the second type that correspond to the algebras described in the previous subsection.
{\unitlength 1.7mm
\begin{figure}[h]
\begin{picture}(10,8)
\put(0,5){\line(1,0){10}}
\put(2,4.5){\makebox(0,0)[cc]{$*$}}
\put(2,5.3){\makebox(0,0)[cc]{$*$}}
\put(8,5){\makebox(0,0)[cc]{$*$}}
\put(4,0){\small $\Upsilon_1$}
\end{picture} \quad 
\begin{picture}(10,8)
\put(0,5){\line(1,0){10}}
\put(2,5){\makebox(0,0)[cc]{$*$}}
\put(8,4.5){\makebox(0,0)[cc]{$*$}}
\put(8,5.3){\makebox(0,0)[cc]{$*$}}
\put(4,0){\small $\Upsilon_2$}
\end{picture} \quad 
\begin{picture}(9,8)
\put(0,5){\line(1,0){9}}
\put(2,4.5){\makebox(0,0)[cc]{$*$}}
\put(2,5.3){\makebox(0,0)[cc]{$*$}}
\put(4.5,5){\oval(5,3)[t]}
\put(7,5){\circle*{0.5}}
\put(3.5,0){\small $\Omega_1$}
\end{picture} \quad 
\begin{picture}(9,8)
\put(0,5){\line(1,0){9}}
\put(2,5){\circle*{0.5}}
\put(7,4.5){\makebox(0,0)[cc]{$*$}}
\put(7,5.3){\makebox(0,0)[cc]{$*$}}
\put(4.5,5){\oval(5,3)[t]}
\put(3.5,0){\small $\Omega_2$}
\end{picture}
\quad 
\begin{picture}(9,8)
\put(0,5){\line(1,0){9}}
\put(2,5){\makebox(0,0)[cc]{$*$}}
\put(7,5){\makebox(0,0)[cc]{$*$}}
\put(4.5,5){\oval(5,3)[t]}
\put(3.5,0){\small $\Lambda$}
\end{picture} \qquad 
\begin{picture}(6,7)
\put(0,5){\line(1,0){6}}
\put(3,4.4){\makebox(0,0)[cc]{$\ast$}}
\put(3,5.4){\makebox(0,0)[cc]{$\ast$}}
\put(3,5){\oval(1.7,2.7)}
\put(1.5,0){\small $\Theta$}
\end{picture}
\caption{$\Upsilon$, $\Omega$, $\Lambda$, and $\Theta$}
\label{Omega}
\end{figure}
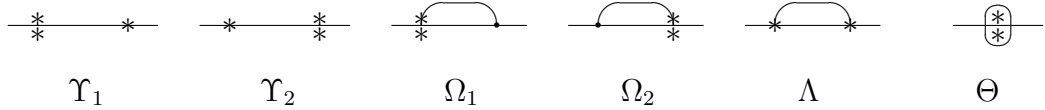
}

\subsection{One- and zero-dimensional cells}
\label{onedim}

There is only one one-dimensional cell of the first type; see Fig. \ref{nabla}. It is swept out by algebras $\nabla(A)$ which are defined by the condition $f'(A)=f''(A)=f'''(A)=0$ for some $A \in S^1 \setminus \{\bullet\}$. 

Additionally, there are six one-dimensional cells $\bar \Upsilon_1,$ $\bar \Upsilon_2$, $\bar \Omega_1$, $\bar \Omega_2$, $\bar \Lambda$, and $\bar \Theta$ of the second type corresponding to the algebras described in the previous subsection.

\begin{figure}[h]
\begin{center}
\begin{picture}(10,7)
\put(0,5){\line(1,0){10}}
\put(5,4.3){\makebox(0,0)[cc]{$*$}}
\put(5,5){\makebox(0,0)[cc]{$*$}}
\put(5,5.7){\makebox(0,0)[cc]{$*$}}
\put(2,0){\small $\nabla$}
\end{picture}
\end{center}
\caption{$\nabla$}
\label{nabla}
\end{figure}

Finally, there is the single zero-dimensional cell $\bar \nabla$, defined by the condition $f'(\bullet) = f''(\bullet) = f'''(\bullet) =0$.

\begin{proposition}
\label{filthprop}
The multiplicities $($see Definition \ref{muldef}$)$ of the equilevel algebras of codimension three described above are as follows: algebras of classes $A,$ $ B,$ $ D, $ $E,$ $ a, $ $e,$ $ g,$ $ h,$ $ i, $ $X,$ $ Y,$ $ J,$ $ S$ are of multiplicity 6; algebras of classes $C,$ $ b,$ $ c, $ $k,$ $ U, $ $V,$ $ Z,$ $ T, \Upsilon,$ $ \Theta$ are of multiplicity 5; algebras of classes $d,$ $ W, $ $\Omega,$ $ \Lambda,$ $ \nabla$ are of multiplicity 4. \hfill $\Box_{\ref{filthprop}}$
\end{proposition}

\subsection{Geometry of some cycles}

By Proposition \ref{trui}, each algebra $\F = S(\alpha, \beta; A)$ is characterized by 
the 2-plane $$(\F \cap \M_A)/\M_A^6 \subset \M_A^3 / \M_A^6 , $$ which is
a point of the Grassmann manifold $G_2(\M^3_A/\M^6_A) \simeq {\mathbb R}P^2$. Each algebra $J(\alpha, \beta; A)$ is characterized by a point of the Grassmann manifold $G_2(\M^2_A/\M^6_A) \simeq G_2({\mathbb R}^4)$. Each algebra $h(\alpha, \beta; A, B)$ is characterized by a point of the manifold $G_2((\M_A \cap \M_B)/(\M^3_A \cap \M_B^3)) \simeq G_2({\mathbb R}^4)$. 
The definitions of these algebras imply the following statements.

\begin{proposition}
\label{Sstruc}
For each $A \in S^1$, the closure of the set of all planes $$(S(\alpha, \beta; A) \cap \M_A)/\M_A^6$$ is the entire space $G_2(\M^3_A/\M^6_A) \simeq {\mathbb R}P^2$. It consists of a two-dimensional family of planes defined by algebras $S(\alpha, \beta; A)$ with parameters $\alpha $ and $\beta$ taking arbitrary real values, a one-dimensional family of planes corresponding to algebras $\Theta(\alpha; A)$, and a plane defined by the algebra $\nabla(A)$. \hfill $\Box_{\ref{Sstruc}}$
\end{proposition}

\begin{proposition}
\label{Jstruc}
For any point $A \in S^1$, the closure of the set of all two-planes $(J(\alpha, \beta; A)\cap \M_A)/ \M_A^6$ is a variety homeomorphic to the two-dimensional sphere with identified poles. It consists of a two-dimensional cell defined by the algebras $J(\alpha, \beta; A)$, a one-dimensional cell defined by the algebras $\Theta(\alpha, A)$, and a point corresponding to the algebra $\nabla(A)$. It is smooth except at the point of class $\nabla$, near which it is locally diffeomorphic to a quadratic cone in ${\mathbb R}^3$. In Pl\"ucker coordinates of the Grassmann manifold $G_2(\M^2_A/\M^6_A)$ related to the coordinates $f''(A),$ $ f'''(A), $ $f^{\tiny \mbox{IV}}(A),$ $ f^{\tiny \mbox{V}}(A) $ of the space $\M^2_A/\M^6_A$, it is given by the equations $$p_{12}=0, \quad 3 p_{14}= 10 p_{23}, \quad
3 p_{13}p_{24}=10 p_{23}^2 \ . \eqno {\Box_{\ref{Jstruc}}} $$
\end{proposition}

\begin{proposition}
\label{proph}
For any pair of points $A$ and $B$, the closure of the union of all two-planes $(h(\alpha, \beta; A, B) \cap \M_A \cap \M_B)/(\M_A^3 \cap \M_B^3)$ is a variety homeomorphic to the projective plane. This variety consists of two two-dimensional cells, whose points correspond to the algebras $h(\alpha, \beta; A, B)$ with $\alpha>0$ or $\alpha<0$ and arbitrary real $\beta$, four intervals of classes $Z_{1,2}^\pm$, and three points of classes $\Omega_1$, $\Omega_2$ and $\Lambda$. It is smooth except at the point of type $\Lambda$. In Pl\"ucker coordinates of the Grassmann manifold $G_2((\M_A \cap \M_B)/(\M_A^3 \cap \M_B^3))$, 
related to the coordinates $f'(A), f''(A), f'(B), $ and $ f''(B)$ in $(\M_A/\M_A^3) \times (\M_B /\M_B^3),$
 it is given by the equations 
\begin{equation}
\label{eq1}
p_{13}=0, \quad p_{14}^2 + p_{23}p_{34}=0, \quad p_{12}p_{34} + p_{14}p_{23}=0, \quad p_{23}^2 -p_{14}p_{12}=0 \ . 
\end{equation}
 In these coordinates, its singular point is given by \ 
\begin{equation}p_{12}=p_{13}=p_{14}=p_{23}=p_{34}=0. 
\label{sipo}
\end{equation} 
The closure of the set of $Z_1$ points in this variety is given by equations $p_{12}=p_{13}=p_{14}=p_{23}=0$. Specifically,  the subset $Z_1^+ $ $($respectively, $Z_1^-$, respectively, $\Omega_1)$ is distinguished by the additional condition $p_{24}>0$ $($respectively, $p_{24}<0$, respectively, $p_{24}=0)$. 
The closure of the set of $Z_2$ points is given by equations $p_{13}=p_{14}=p_{23}=p_{34}=0$; the subset $Z_2^+ $ $($respectively, $Z_2^-$, respectively, $\Omega_2)$ is distinguished in it by the same condition $p_{24}>0$ $($respectively, $p_{24}<0$, respectively, $p_{24}=0)$. 
\end{proposition}

\unitlength 1.3mm
\begin{figure}
\begin{picture}(35,40)
\bezier{600}(20,0)(-10,20)(20,40)
\bezier{600}(20,0)(50,20)(20,40)
\bezier{600}(20,0)(12,20)(20,40)
\bezier{150}(20,0)(28,20)(20,40)
\put(21.5,0){\makebox(0,0)[cc]{\tiny $\Lambda$}}
\put(22,40){\makebox(0,0)[cc]{\tiny $\Lambda$}}
\put(3,20){\makebox(0,0)[cc]{\tiny $\Omega_1$}}
\put(37,20){\makebox(0,0)[cc]{\tiny $\Omega_1$}}
\put(14.5,17){\makebox(0,0)[cc]{\tiny $\Omega_2$}}
\put(26,22.5){\makebox(0,0)[cc]{\tiny $\Omega_2$}}
\put(8,8){\makebox(0,0)[cc]{\tiny $Z_1^+$}}
\put(8.5,32){\makebox(0,0)[cc]{\tiny $Z_1^-$}}
\put(32,8){\makebox(0,0)[cc]{\tiny $Z_1^-$}}
\put(32,32){\makebox(0,0)[cc]{\tiny $Z_1^+$}}
\put(15.5,9){\makebox(0,0)[cc]{\tiny $Z_2^+$}}
\put(15.5,31){\makebox(0,0)[cc]{\tiny $Z_2^-$}}
\put(24.5,9){\makebox(0,0)[cc]{\tiny $Z_2^-$}}
\put(24.5,31){\makebox(0,0)[cc]{\tiny $Z_2^+$}}
\put(5,20){\circle*{0.8}}
\put(35,20){\circle*{0.8}}
\put(16,18){\circle*{0.8}}
\put(24,22){\circle*{0.8}}
\end{picture}
\caption{Preimage of a cycle \{$h(\cdot, \cdot; A, B)$\}}
\label{alph}
\end{figure}

\begin{remark} \rm 
One can consider this variety as the quotient space of a topological two-sphere by the central symmetry, where the point of class $\Lambda$ is obtained from the poles, and the algebras of classes $Z^{\pm}_i$ correspond to the meridian points with longitudes $0, 180, 90, $ and $270$, with equatorial points of these meridians corresponding to points of classes $\Omega_i$, see Fig.~\ref{alph}.
\end{remark}

\noindent
{\it Proof of Proposition \ref{proph}.} The equations (\ref{eq1}) follow directly from the conditions (\ref{eqh}).
The regularity of their solution set outside the point (\ref{sipo}) follows immediately from these equations.
Consider the conic subvariety in ${\mathbb R}^6$ defined by these equations.
 The coordinates $p_{12}, p_{24}$, and $p_{34}$ cannot all vanish simultaneously at a non-zero point of this subvariety. Its projection to the subspace ${\mathbb R}^3$ with these coordinates, defined by forgetting the remaining coordinates, is a homeomorphism which is linear on any line through the origin. Thus, this projection defines a homeomorphism of the projectivization of this subvariety to ${\mathbb R}P^2$. The conditions for $Z_i$ sets follow immediately from their definitions.\hfill $\Box_{\ref{proph}} $

\section{Completeness of the list of equilevel subalgebras of codimension three in $C^\infty(S^1, {\mathbb R})$}
\label{class}

We must prove that all equilevel  subalgebras of codimension three are among the algebras counted in Section \ref{cells}. We will prove this separately for subalgebras whose supports consist of one, two, \dots, six points. Later we will see that all of subalgebras counted in Section \ref{cells} indeed appear in $\overline{CD}_3(S^1)$. 

\subsection{Algebras with one-point support} 

Let $\F$ be an equilevel subalgebra of codimension three in $C^\infty(S^1, {\mathbb R})$ 
with the support at the single point $A \in S^1$.
By Proposition \ref{trui}, $\F$ contains the sixth power $\M_A^6 $ of the maximal ideal in $C^\infty(S^1, {\mathbb R})$ with the support at $A$, as well as all constant functions. All chords of  the configurations $\nu(t)$ participating in the definition of $\F$ contract to the point $A$ when $t$ tends to $0$. Therefore, 
 $f'(A)=0$ for all $f \in \F$. The following possibilities then arise.

\begin{itemize}
\item[1)]
\label{i1}
$f'(A)=f''(A)=f'''(A)=0$ for all elements of the algebra $\F$. 
These are already three conditions and we get a single algebra that depends only on the point $A$. In our notation, it is called $\nabla(A)$.

\item[2)] $f'(A)=f''(A)\equiv 0$ for all $f \in \F$, but there are functions $f \in \F$ with $f'''(A) \neq 0$. The algebra $\F \cap \M_A$ then defines a two-dimensional subspace in the quotient space $\M^3_A / \M^6_A$, so it is given by a linear relation between the coefficients at $x^3$, $x^4$ and $x^5$. In other words, an identity of the form
\begin{equation}
\label{i2}
 \omega f^{\tiny \mbox{V}}(A) + \alpha f^{\tiny \mbox{IV}}(A) + \beta f'''(A) =0
\end{equation}
 for some point $(\omega:\alpha:\beta) \in {\mathbb R}P^2$ should be satisfied for all $f \in \F$.
Any plane defined by such a relation defines a subalgebra of codimension three. The case $\omega=\alpha=0$ was considered in the previous paragraph 1) and contradicts the condition $f'''(A) \not \equiv 0$ which is now presumed. Subalgebras that satisfy the identities (\ref{i2}) with $\omega \neq 0$ depend on two parameters $\alpha/\omega$ and $\beta/\omega$ and the point $A$. These subalgebras form the cells $S$ and $\bar S$, see \S \ref{threedim}. Algebras with $\omega=0,$ $ \alpha \neq 0$ depend on the parameter $\beta/\alpha$ and the point $A$ and form the cells $\Theta$ and $\bar \Theta$, see \S \ref{twodim}. 

\item[3)] 
 $f'(A) \equiv 0$ for all $f \in \F$, but there exists a function $f_0 \in \F$ equal to $x^2+ \alpha x^3 + \omega x^4 + \varkappa x^5$
in the local coordinate $x$ centered at $A$. We can assume that \ $\omega=0$ \ by subtracting $\omega f_0^2$. 
Our subalgebra defines a subspace of codimension two (and thus of dimension two) in the quotient space $\M^2_A/\M^6_A$, so it is linearly spanned by $f_0$ and $f_0^2$ and is completely determined by the coefficients $\alpha$ and $\varkappa$ and the point $A$. All of its elements therefore satisfy the conditions (\ref{jten}) with $\beta = 60 \varkappa$. Thus, the algebras of this kind with $A \in (0,2\pi)$ fill a three-dimensional cell, which we call $J$. The algebras with $A= \bullet$ form the two-dimensional cell $\bar J$.
\end{itemize}

\subsection{Algebras with two-point support}

Let the support of the equilevel subalgebra $\F$ of codimension three consist of two points $A$ and $ B$. Let $\{\nu(t)\}$ be a family of matched six-configurations used in the definition of $\F$, see \S~\ref{maindef}. Each point $\nu(t),$ $t>0$, defines an algebra of type $A_{i j}$ which can be depicted by a 3-chord diagram. As $t$ approaches $0$, these chords 
tend to a graph with two vertices $A, B$ and three edges. This graph can have only one of the five types shown in Fig.~\ref{grrphs} (up to reordering of the vertices).

\unitlength=1.5mm

\begin{figure}
\begin{picture}(14,12)
\put(0,5){\line(1,0){14}}
\put(3,5){\circle*{0.5}}
\bezier{100}(3,5)(0,10)(3,10)
\bezier{100}(3,5)(6,10)(3,10)
\bezier{100}(3,5)(0,0)(3,0)
\bezier{100}(3,5)(6,0)(3,0)
\put(11,5){\circle*{0.5}}
\bezier{100}(11,5)(8,10)(11,10)
\bezier{100}(11,5)(14,10)(11,10)
\put(7,0){\footnotesize 1}
\end{picture} \quad 
\begin{picture}(14,10)
\put(0,5){\line(1,0){14}}
\put(3,5){\circle*{0.5}}
\bezier{100}(3,5)(0,10)(3,10)
\bezier{100}(3,5)(6,10)(3,10)
\bezier{100}(3,5)(0,0)(3,0)
\bezier{100}(3,5)(6,0)(3,0)
\put(11,5){\circle*{0.5}}
\bezier{150}(3,5)(7,8)(11,5)
\put(7,0){\footnotesize 2}
\end{picture} \quad 
\begin{picture}(14,10)
\put(0,5){\line(1,0){14}}
\put(3,5){\circle*{0.5}}
\bezier{100}(3,5)(0,10)(3,10)
\bezier{100}(3,5)(6,10)(3,10)
\put(11,5){\circle*{0.5}}
\bezier{100}(11,5)(8,10)(11,10)
\bezier{100}(11,5)(14,10)(11,10)
\bezier{150}(3,5)(7,2)(11,5)
\put(6,0){\footnotesize 3}
\end{picture} \quad 
\begin{picture}(14,10)
\put(0,5){\line(1,0){14}}
\put(3,5){\circle*{0.5}}
\bezier{100}(3,5)(0,0)(3,0)
\bezier{100}(3,5)(6,0)(3,0)
\put(11,5){\circle*{0.5}}
\bezier{150}(3,5)(7,11)(11,5)
\bezier{150}(3,5)(7,8)(11,5)
\put(7,0){\footnotesize 4}
\end{picture} \quad 
\begin{picture}(14,10)
\put(0,5){\line(1,0){14}}
\put(3,5){\circle*{0.5}}
\put(11,5){\circle*{0.5}}
\bezier{150}(3,5)(7,10)(11,5)
\bezier{150}(3,5)(7,2)(11,5)
\bezier{150}(3,5)(7,8)(11,5)
\put(6,0){\footnotesize 5}
\end{picture}
\caption{Limits of chords: two-point support}
\label{grrphs}
\end{figure}
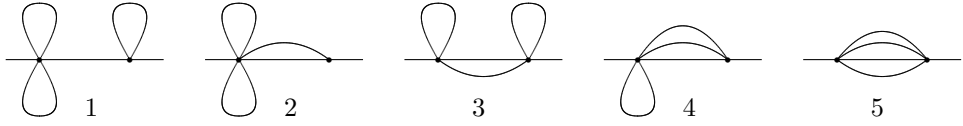
 
Let us analyze the subalgebras that can occur at the limit point $\nu(0).$ 

\begin{itemize}
\item[1.] 
If the endpoints of each chord tend to a single point of the support, then the endpoints of one chord tend to one point, while the endpoints of two other chords tend to the other point. The first point gives us only one vanishing derivative condition. The analysis of the conditions at the second point repeats the consideration of subalgebras of codimension two with one-point support, see \S~\ref{prompro}. Namely, we necessarily have $f' = 0$ at this point.
Our algebra is of type $\Upsilon$ if its third condition at this point is $f''=0$; it is of type $Y$ if this is the condition $f''' = \alpha f''$. Depending on the disposition of the support points in the circle, the algebras of each of these two classes can belong to four different cells whose notations can have lower index 1 or 2 and also to have or not to have the bar above. 
\item[2.] If the limit graph has the vertices $A$ and $B$ of valences 5 and 1, then we immediately obtain two relations $f(A) = f(B)$ and $f'(A)=0$. The third relation should distinguish a two-dimensional subspace in the three-dimensional space $\M^2_A/\M^5_A$. If $f''(A)=0$ for all functions in our algebra, then it is again an algebra of the class $\Upsilon$. Otherwise, $\F$ contains a function $f_0 = x^2 + a x^3 + b x^4$, and our 2-subspace is generated by the classes of functions $f_0$ and $f_0^2$, in particular it is completely characterized by the number $a \equiv \frac{f'''_0(A)}{3 f''_0(A)}$. Thus, we have an algebra of the class $T$, see Fig. \ref{TJ}.
\item[3.] The third graph in Fig.~\ref{grrphs} yields three conditions $f(A)=f(B),$ $f'(A)=0=f'(B)$, and we have an algebra of class $\Lambda$. It depends on only two points $A$ and $B$. 
\item[4.] The fourth graph assumes two conditions $f(A) = f(B) $ and $f'(A) = 0$ for all elements of our algebra. The third condition should distinguish a two-subspace in the three-dimensional space $(\M^2_A \cap \M_B)/(\M^4_A \cap \M^2_B)$,
so it is given by an equation as on the right side of (\ref{ialg}). If in this equation $\alpha=\beta=0$, then we have $f'(B)=0$ for all functions $f \in \F$, and $\F$ is of the class $\Lambda$. If $\alpha = \gamma =0$, then $f''(A) =0$ for all $f \in \F$, and $\F$ is of the class $\Omega$. If $\alpha = 0$ but $\beta \gamma \neq 0$, then the algebra is of class $Z$. If $\gamma=0$ but $\alpha \beta \neq 0$ then it is of class $T$. 
In all other cases, the algebra is of class $i$.
\item[5.] In the case of the fifth graph in Fig.~\ref{grrphs} we have the condition $f(A)=f(B)$, and our algebra defines a two-dimensional subspace in the four-dimensional space $(\M_A \cap \M_B)/(\M^3_A \cap \M^3_B)$. If $f'(A) = f'(B) =0$ for all functions in the algebra, then it is of class $\Lambda$. If the first and the second derivatives of all functions in $\F$ vanish at one endpoint (say, $A$), then $\F$ is an algebra of type $\Omega$. If only the first derivative at $A$ vanishes identically, but the first derivative at $B$ and the second derivative at $A$ generally do not, then we again have a subalgebra of type $Z$. Finally, suppose that there is a function $f_0$ in our algebra with 
\begin{equation}
f_0(A)=f_0(B)=0, \quad f'_0(A) \neq 0 \neq f'_0(B). 
\label{con5}
\end{equation}
Norming \ $f'_0(A)$, \ subtracting \ $\frac{f''_0(A)}{2} f^2_0(A)$ \ and reducing modulo \ $\M^3_A \cap \M^3_B$, we can assume that $f_0 \equiv x$ in the local coordinate $x \equiv \varphi-\varphi(A) $ with the origin at $A$ and $f_0 \equiv \alpha y + \frac{1}{2} \beta y^2$ in the coordinate $y \equiv \varphi - \varphi(B)$ with the origin $B$. These numbers $\alpha$ and $\beta$ completely characterize our subalgebra, which, in this case, is of class $h$, see \S \ref{fourdim}.
\end{itemize}


\subsection{Algebras with three-point support}

There are eight possible types of graphs with three vertices that describe the limit points $\nu(0)$ of matched 3-chord diagrams, see Fig.~\ref{grrphs3}. The corresponding possible algebras of codimension three can be analyzed as in the previous subsection and give us the following classes of algebras. 
1: $X$. 2: $V$ and $g$. 3: $U$. 4: $U$ and $e$. 5 and 6: $W$. 7: $W$ and $k$. 8: $W$, $k$, and $E$.

{\unitlength=1.3mm
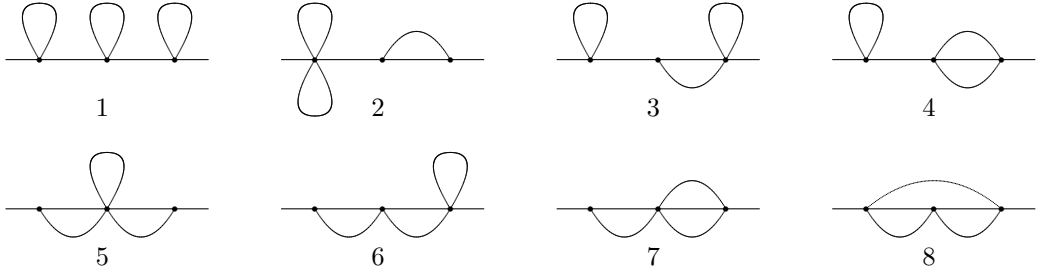
\begin{figure}
\begin{picture}(18,12)
\put(0,5){\line(1,0){18}}
\put(3,5){\circle*{0.5}}
\bezier{100}(3,5)(0,10)(3,10)
\bezier{100}(3,5)(6,10)(3,10)
\bezier{100}(9,5)(6,10)(9,10)
\bezier{100}(9,5)(12,10)(9,10)
\put(9,5){\circle*{0.5}}
\put(15,5){\circle*{0.5}}
\bezier{100}(15,5)(12,10)(15,10)
\bezier{100}(15,5)(18,10)(15,10)
\put(8,0){\footnotesize 1}
\end{picture} \qquad
\begin{picture}(18,10)
\put(0,5){\line(1,0){18}}
\put(3,5){\circle*{0.5}}
\bezier{100}(3,5)(0,10)(3,10)
\bezier{100}(3,5)(6,10)(3,10)
\bezier{100}(3,5)(0,0)(3,0)
\bezier{100}(3,5)(6,0)(3,0)
\put(9,5){\circle*{0.5}}
\put(15,5){\circle*{0.5}}
\bezier{100}(15,5)(12,10)(9,5)
\put(8,0){\footnotesize 2}
\end{picture} \qquad 
\begin{picture}(18,10)
\put(0,5){\line(1,0){18}}
\put(3,5){\circle*{0.5}}
\bezier{100}(3,5)(0,10)(3,10)
\bezier{100}(3,5)(6,10)(3,10)
\bezier{100}(15,5)(12,10)(15,10)
\bezier{100}(15,5)(18,10)(15,10)
\put(9,5){\circle*{0.5}}
\put(15,5){\circle*{0.5}}
\bezier{100}(15,5)(12,0)(9,5)
\put(8,0){\footnotesize 3}
\end{picture} \qquad 
\begin{picture}(18,10)
\put(0,5){\line(1,0){18}}
\put(3,5){\circle*{0.5}}
\bezier{100}(3,5)(0,10)(3,10)
\bezier{100}(3,5)(6,10)(3,10)
\bezier{100}(15,5)(12,10)(9,5)
\put(9,5){\circle*{0.5}}
\put(15,5){\circle*{0.5}}
\bezier{100}(15,5)(12,0)(9,5)
\put(8,0){\footnotesize 4}
\end{picture} 

\begin{picture}(18,13)
\put(0,5){\line(1,0){18}}
\put(3,5){\circle*{0.5}}
\bezier{100}(3,5)(6,0)(9,5)
\bezier{100}(9,5)(6,10)(9,10)
\bezier{100}(9,5)(12,10)(9,10)
\put(9,5){\circle*{0.5}}
\put(15,5){\circle*{0.5}}
\bezier{100}(15,5)(12,0)(9,5)
\put(8,0){\footnotesize 5}
\end{picture} \qquad
\begin{picture}(18,10)
\put(0,5){\line(1,0){18}}
\put(3,5){\circle*{0.5}}
\bezier{100}(3,5)(6,0)(9,5)
\bezier{100}(9,5)(12,0)(15,5)
\put(9,5){\circle*{0.5}}
\put(15,5){\circle*{0.5}}
\bezier{100}(15,5)(12,10)(15,10)
\bezier{100}(15,5)(18,10)(15,10)
\put(8,0){\footnotesize 6}
\end{picture} \qquad
\begin{picture}(18,10)
\put(0,5){\line(1,0){18}}
\put(3,5){\circle*{0.5}}
\bezier{100}(3,5)(6,0)(9,5)
\bezier{100}(9,5)(12,0)(15,5)
\bezier{100}(9,5)(12,10)(15,5)
\put(9,5){\circle*{0.5}}
\put(15,5){\circle*{0.5}}
\put(8,0){\footnotesize 7}
\end{picture} \qquad
\begin{picture}(18,10)
\put(0,5){\line(1,0){18}}
\put(3,5){\circle*{0.5}}
\bezier{100}(3,5)(6,0)(9,5)
\bezier{100}(9,5)(12,0)(15,5)
\bezier{100}(3,5)(9,10)(15,5)
\put(9,5){\circle*{0.5}}
\put(15,5){\circle*{0.5}}
\put(8,0){\footnotesize 8}
\end{picture}
\caption{Limits of chord diagrams: three-point support}
\label{grrphs3}
\end{figure}
}

\subsection{Algebras with four-point support}

There are five types of limit graphs with four vertices, see Fig.~\ref{grrphs4}.
\begin{figure}
\begin{picture}(24,10)
\put(0,5){\line(1,0){24}}
\put(3,5){\circle*{0.5}}
\bezier{100}(3,5)(0,10)(3,10)
\bezier{100}(3,5)(6,10)(3,10)
\bezier{100}(9,5)(6,10)(9,10)
\bezier{100}(9,5)(12,10)(9,10)
\put(9,5){\circle*{0.5}}
\put(15,5){\circle*{0.5}}
\put(21,5){\circle*{0.5}}
\bezier{100}(15,5)(18,10)(21,5)
\put(11,0){\footnotesize 1}
\end{picture} \qquad
\begin{picture}(24,10)
\put(0,5){\line(1,0){24}}
\put(3,5){\circle*{0.5}}
\bezier{100}(3,5)(0,10)(3,10)
\bezier{100}(3,5)(6,10)(3,10)
\bezier{100}(21,5)(18,10)(15,5)
\put(9,5){\circle*{0.5}}
\put(15,5){\circle*{0.5}}
\put(21,5){\circle*{0.5}}
\bezier{100}(9,5)(6,0)(3,5)
\put(11,0){\footnotesize 2}
\end{picture} \qquad 
\begin{picture}(24,10)
\put(0,5){\line(1,0){24}}
\put(3,5){\circle*{0.5}}
\bezier{100}(3,5)(0,10)(3,10)
\bezier{100}(3,5)(6,10)(3,10)
\bezier{100}(9,5)(12,10)(15,5)
\put(9,5){\circle*{0.5}}
\put(15,5){\circle*{0.5}}
\put(21,5){\circle*{0.5}}
\bezier{100}(15,5)(18,10)(21,5)
\put(11,0){\footnotesize 3}
\end{picture} 

\begin{picture}(24,10)
\put(0,5){\line(1,0){24}}
\put(3,5){\circle*{0.5}}
\bezier{100}(3,5)(6,10)(9,5)
\bezier{100}(9,5)(12,10)(15,5)
\put(9,5){\circle*{0.5}}
\put(15,5){\circle*{0.5}}
\put(21,5){\circle*{0.5}}
\bezier{100}(15,5)(18,10)(21,5)
\put(11,0){\footnotesize 4}
\end{picture} \qquad
\begin{picture}(24,10)
\put(0,5){\line(1,0){24}}
\put(3,5){\circle*{0.5}}
\bezier{100}(3,5)(6,10)(9,5)
\bezier{100}(21,5)(18,10)(15,5)
\put(9,5){\circle*{0.5}}
\put(15,5){\circle*{0.5}}
\put(21,5){\circle*{0.5}}
\bezier{100}(9,5)(6,0)(3,5)
\put(11,0){\footnotesize 5}
\end{picture}
\caption{Limits of chord diagrams: four-point support}
\label{grrphs4}
\end{figure}
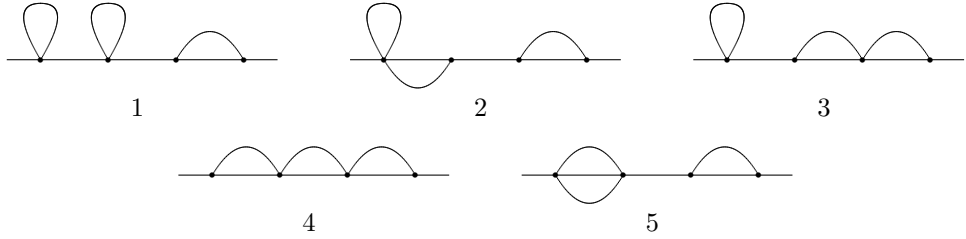
They provide the algebras of the following types: 
1: $a$. 2: $b$. 3: $c$. 4: $d$. 5: $b$ and $D$.
\smallskip

There are also only two types of limit graphs with five vertices and only one type of limit graphs with six vertices, see Fig. \ref{grrphs56}. These graphs correspond to the algebras of classes $B$, $C$, and $A$, respectively.

{\unitlength 1.4mm
\begin{figure}
\begin{picture}(24,10)
\put(0,5){\line(1,0){24}}
\put(2,5){\circle*{0.5}}
\bezier{100}(2,5)(0,9)(2,9)
\bezier{100}(2,5)(4,9)(2,9)
\bezier{100}(7,5)(9.5,8)(12,5)
\put(7,5){\circle*{0.5}}
\put(12,5){\circle*{0.5}}
\put(17,5){\circle*{0.5}}
\put(22,5){\circle*{0.5}}
\bezier{100}(17,5)(19.5,8)(22,5)
\end{picture} \qquad
\begin{picture}(24,10)
\put(0,5){\line(1,0){24}}
\put(2,5){\circle*{0.5}}
\bezier{100}(2,5)(4.5,8)(7,5)
\bezier{100}(7,5)(9.5,8)(12,5)
\put(7,5){\circle*{0.5}}
\put(12,5){\circle*{0.5}}
\put(17,5){\circle*{0.5}}
\put(22,5){\circle*{0.5}}
\bezier{100}(17,5)(19.5,8)(22,5)
\end{picture}
\qquad
\begin{picture}(29,10)
\put(0,5){\line(1,0){29}}
\put(2,5){\circle*{0.5}}
\bezier{100}(2,5)(4.5,8)(7,5)
\bezier{100}(17,5)(14.5,8)(12,5)
\put(7,5){\circle*{0.5}}
\put(12,5){\circle*{0.5}}
\put(17,5){\circle*{0.5}}
\put(22,5){\circle*{0.5}}
\put(27,5){\circle*{0.5}}
\bezier{100}(27,5)(24.5,8)(22,5)
\end{picture}
\caption{Limits of chord diagrams: five- and six-point supports}
\label{grrphs56}
\end{figure}
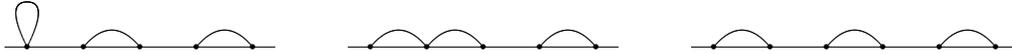
}

The incidence coefficients of the boundary maps, which will be listed in \S \ref{ic}, prove that all these algebras arise from chains of degenerations of 3-chord diagrams. Therefore, they are points of the space $\overline{CD}_3(S^1)$. 
\medskip

\noindent
{\it End of the proof of Theorem \ref{thmcell}.} 
We represented the variety $\overline{CD}_3(S^1)$ as the union of non-intersecting open cells. The proof that these cells define the structure of a CW-complex repeats that for the space $\overline{CD}_2(S^1)$, see Lemma \ref{le21}. The only difference is that $\overline{CD}_2(S^1)$ is replaced by $\overline{CD}_3(S^1)$ and the Grassmann manifold $G_{\dim {\mathcal L}-2}(\mathcal L)$ is replaced by $G_{\dim {\mathcal L}-3}(\mathcal L)$. \hfill $\Box_{\ref{thmcell}}$

\section{Boundary operators mod 2 in the complex $\overline{CD}_3(S^1)$}
\label{ic}

\subsection{}

Consider the cellular chain complex with coefficients in ${\mathbb Z}_2$ defined by our CW-structure of the variety $\overline{CD}_3(S^1)$. 
Most of the incidence coefficients that comprise the formulas for its boundary operators $C_i \to C_{i-1}$ are obvious. In particular, this is the case if the incidence of cells is induced by the collision of only two points of the supports of the subalgebras from the larger cell (see for example the degeneration of cells of type $A_{i j}$ into cells $B_k^l$ or $C_{k l}$), or it is induced by the boundary values of the numerical conditions on the derivatives at the points of a fixed configuration (as in the case of the cells $D^\pm_{i j}$ approaching the cells $b_{i j}$ and $b_{j i}$). The most difficult adjacencies are related with the simultaneous collision of many points of the supports: the asymptotics of these collisions impose conditions on the derivatives of the functions in the algebras of the limit cells. We give here two illustrations.

\subsubsection{Approximation of a five-dimensional cell of class $E$ by the six-dimensional cell $A_{41}$.} 
\label{sss811}
This adjacency is defined by the simultaneous contraction of the segments $[A,B]$, $[CD]$, and $[EF]$ in the lower part of Fig.~\ref{incid1}. 
\unitlength 0.3mm
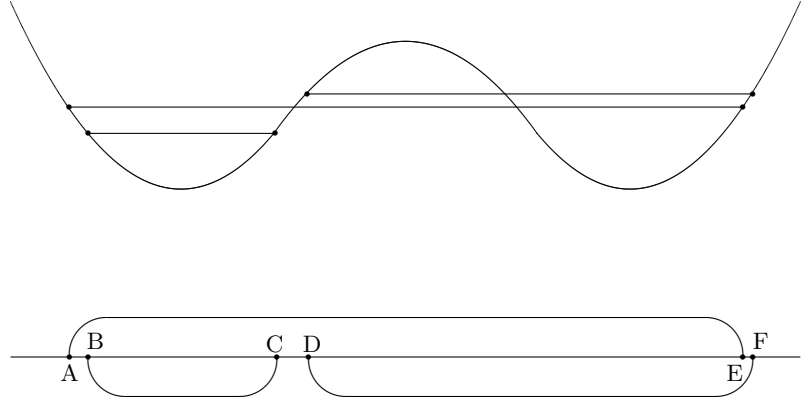
\begin{figure}
\begin{picture}(300, 150)
\bezier{500}(0,145)(50,35)(100,95)
\bezier{500}(100,95)(150,165)(200,95)
\bezier{500}(200,95)(250,35)(300,145)
\put(0,10){\line(1,0){300}}
\put(22.3,105){\circle*{2}}
\put(29.5,95){\circle*{2}}
\put(100.5,95){\circle*{2}}
\put(112.7,110){\circle*{2}}
\put(278.2,105){\circle*{2}}
\put(282,110){\circle*{2}}
\put(22.4,105){\line(1,0){255.8}}
\put(29.5,95){\line(1,0){71.7}}
\put(113.2,110){\line(1,0){168.8}}
\put(22.4,10){\circle*{2}}
\put(29.5,10){\circle*{2}}
\put(101.2,10){\circle*{2}}
\put(113.2,10){\circle*{2}}
\put(278.2,10){\circle*{2}}
\put(282,10){\circle*{2}}
\put(150.3,10){\oval(255.8,30)[t]}
\put(65.35,10){\oval(71.7,30)[b]}
\put(197.6,10){\oval(168.8,30)[b]}
\put(19,1){\scriptsize A}
\put(29,13){\scriptsize B}
\put(97,12){\scriptsize C}
\put(111,12){\scriptsize D}
\put(272,1){\scriptsize E}
\put(282,13){\scriptsize F}
\end{picture}
\caption{Incidence coefficient $[A_{41}, E_3]$}
\label{incid1}
\end{figure}
Let the lengths of these segments be equal to $\alpha t$, $\beta t$, and $\gamma t$, respectively (up to the higher degree terms in $t$), where \ $t \to +0$ \ is the parameter of the contraction. Let $X$, $Y$, and $Z$ be the limit points to which these segments are contracted.
At every moment of the contraction, all functions of the corresponding algebra satisfy the identity $$f(F(t)) - f(E(t)) = \left(f(D(t)) - f(C(t))\right) \ + \ \left(f(B(t))-f(A(t))\right),$$
which asymptotically implies the identity $\gamma f'(Z) = \beta f'(Y) + \alpha f'(X)$ for all functions $f$ of the limit algebra. Here, $\alpha, \beta$, and $\gamma$ are arbitrary positive coefficients, so the limit subalgebra is the algebra $E((\alpha: \beta: -\gamma); X, Y, Z)$ of class $E_3$.

\subsubsection{Incidence coefficient $[a_{23}, S]$.}
Let us fix the generic coefficients $\alpha$ and $\beta$, and investigate the space of parametric algebraic curves $\omega: [0, \varepsilon) \to \overline{CD}_3(S^1)$ such that the point $\omega(0)$ is the algebra $S(\alpha, \beta; A)$ and all the algebras $\omega(t),$ $t \in (0,\varepsilon)$, are of the class $a_{23}$.

A typical function of an algebra of class $a_{23}$ with four support points
approaching the collision of all of them looks as shown in the top part of Fig.~\ref{incid2}. 

\begin{figure}
\unitlength 0.4mm
\begin{picture}(200,152)
\bezier{500}(0,50)(60,180)(100,100)
\bezier{500}(100,100)(140,20)(200,150)
\put(67,130){\circle*{3}}
\put(133,69.5){\circle*{3}}
\put(35.4,110){\circle*{3}}
\put(179.3,110){\circle*{3}}
\put(0,12){\line(1,0){200}}
\put(67,12){\makebox(0,0)[cc]{$\ast$}}
\put(133,12){\makebox(0,0)[cc]{$\ast$}}
\put(35.4,12){\circle*{3}}
\put(179.3,12){\circle*{3}}
\put(35,110){\line(1,0){145}}
\put(107.5, 12){\oval(145,20)[t]}
\put(31.5,1){\small $q$}
\put(65,-1){\small 0}
\put(130,1){\small $t$}
\put(177,1){\small $p$}
\end{picture}
\caption{Incidence coefficient $[a_{23}, S]$}
\label{incid2}
\end{figure}
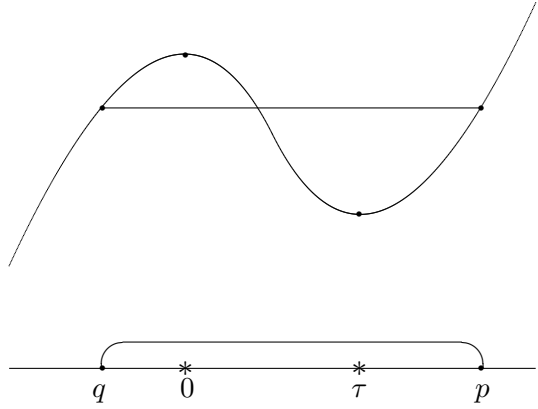

We can assume that the left critical point of the $a_{2 3}$ configuration is fixed at $A$, and that the other three points tend toward it. 
Let $x$ be the local coordinate in $S^1$ centered at $A$.
Then, all intermediate algebras defined by these configurations satisfy the condition $f'(A)\equiv 0$. Take the coordinate $t$ of the other critical point as the parameter of the contraction.
Let the coordinates $p$ and $q$ of the remaining two points of the configuration depend on $t$ as follows:
\begin{eqnarray}
\label{pq}
\left\{
\begin{aligned}
p(t) = \xi t^u + \eta t^v + \zeta t^w+ \dots , \\
q(t) = \lambda t^u + \mu t^v + \varkappa t^w + \dots , \\
\end{aligned}
\right. \ ,
\end{eqnarray} where 
$0< u<v<w< \dots$, the coefficients $\xi $ and $\lambda$ do not vanish simultaneously, the coefficients $\eta$ and $\mu$ also do not vanish  simultaneously, etc. At the end of the contraction, we must obtain the algebra $S(\alpha, \beta; A)$, which is defined by the conditions $f'(A)=0$, $f''(A)=0$, and the linear condition (\ref{i2}) on the values of $f'''(A),$ $ f^{\mbox{\tiny IV}}(A)$ and $f^{\mbox{\tiny V}}(A)$
for all elements $f$ of this algebra. According to Proposition \ref{trui}, this limit algebra is completely characterized by its intersection with the space ${\mathbb R}^4$ of polynomials 
\begin{equation}
a x^2 + b x^3 + c x^4 + d x^5 , 
\label{trans}
\end{equation}
which is transversal to the ideal $\M^6_A $ in the space $\M^2_A$.
Each algebra corresponding to a positive value of the parameter $t$ intersects this space along the set of polynomials (\ref{trans}) that satisfy the two linear conditions $f(p)=f(q)$ and $f'(t)=0$, i.e.,
\begin{multline}
a(p+q)+b(p^2 + p q + q^2) + c(p^3+p^2q+p q^2 + q^3) + \\
d(p^4+ p^3 q + p^2q^2 + p q^3 + p^4) =0 \qquad \qquad
\end{multline}
and
\begin{equation}
2a + 3 b t + 4 c t^2 + 5 d t^3 =0 \ .
\end{equation}
The two hyperplanes in ${\mathbb R}^4$ defined by these conditions tend to the hyperplane $\{a=0\}$ when $t$ tends to $0$. To compute the limit position of their intersection, which is a hyperplane within the latter hyperplane, 
we express the variable \ $a$ \ from one of these equations and substitute it into the other. The resulting condition on the remaining variables is as follows:
\begin{multline}
\label{linbcd}
b\left(2(p^2 + p q + q^2) -3(p+q) t\right) + \\
c\left(2(p^3 + p^2 q + p q^2 + q^3) - 4(p+q) t^2 \right) + \\
 + d\left(2(p^4+ p^3 q + p^2q^2 + p q^3 + p^4) - 5(p+q) t^3 \right) =0 \ .
\end{multline}
Substituting the expressions (\ref{pq})  for $p$ and $q$ into this equation, we obtain a linear equation in the variables $b, c,$ and $d$, whose coefficients are certain power series in $t$. The lowest possible degree of the monomials of these power series at the variable $b$ (respectively, $c$, respectively, $d$) is 2u (respectively, 3u, respectively, 4u). However, some of these monomials should vanish. Indeed, let $\sigma$ be the lowest degree of non-zero monomials of these three power series, so that the monomials in the coefficients of $b, c$ and $d$ are $k t^\sigma, l t^\sigma$, and $m t^\sigma$, respectively. The limit algebra defined by the curve (\ref{pq}) satisfies the characteristic condition (\ref{S}) of the algebra $S(\alpha, \beta; A)$ if and only if 
\begin{equation}
(m: l: k) \equiv (5! : -4! \alpha : -3! \beta)
\label{proto}
\end{equation}
in ${\mathbb R}P^2$. Since $\sigma>4u$, the terms of 
degree below $4u$ of the power series at the variables $b$ and $c$ should vanish. 

If $u<1$ then the degree $2u$ term of the coefficient at $b$ is equal to $\xi^2 + \xi \lambda + \lambda^2$; its vanishing implies the prohibited condition $\xi=0=\lambda$. If $u>1$ then $p(t)<t$ for small $t$, contradicting the shape of the $a_{23}$ configurations. Thus, $u=1$. The vanishing of the degree two monomial in the coefficient at $b$ and of the degree three coefficient in the coefficient  at $c$ imply the system of equations
\begin{equation}
\label{eq101}
2(\xi^2 + \xi \lambda + \lambda^2) -3(\xi+\lambda) = 0, 
\end{equation}
\begin{equation}
\label{eq102}
2(\xi^3 + \xi^2 \lambda + \xi \lambda^2 + \lambda^3) - 4(\xi+\lambda) = 0,
\end{equation}
whose unique solution is
\begin{equation}
\label{solu}
\xi = \frac{\sqrt{3}+1}{2}, \quad
\lambda = \frac{1-\sqrt{3}}{2}.
\end{equation}
These values imply that the lowest-degree monomial in the coefficient at $d$ in (\ref{linbcd}) is equal to $\frac{1}{2}t^4$.

The next terms of the asymptotic expansions of the coefficients at $b$ and $c$ are equal respectively to  
\begin{equation}
\label{eq103}
\left(2(2 \xi \eta + \xi \mu + \lambda \eta + 2 \lambda \mu) -3(\eta+\mu)\right) t^{1+v} \qquad \quad \mbox{and} 
\end{equation} 
\begin{equation}
\label{eq104}
\left(2(3 \xi^2\eta + 2 \xi \eta \lambda + \xi^2 \mu + 2\xi \lambda \mu + \eta \lambda^2 + 3 \lambda^2 \mu) -4 (\eta + \mu)\right) t^{2+v}.
\end{equation}
If $v<2$ then both these terms should vanish. Together with the conditions (\ref{solu}) this gives a contradictory system of equations on $\eta$ and $\mu$.
If $v>2$, then the degrees  in $t$ of all nonzero monomials in the coefficient at $c$ are strictly greater than four. Thus, $v=2$ and the coefficient in (\ref{eq103}) should be equal to zero. Together with (\ref{solu}), this gives  $\eta=\mu.$ The principal (degree four) monomial in the coefficient at $c$ in (\ref{linbcd}) is then equal to $4\eta t^2$. 

If the exponent $w$ in (\ref{pq}) is greater than three, then the degree-four term in the coefficient at $b$ is equal to $6 \eta^2t^4$. Therefore,  it is dependent on the similar term in the coefficient at $c$, and the curves of this form cannot approach general algebras $S(\alpha, \beta; A)$. 

Now, suppose that $w=3$. In this case, the degree four term in the coefficient of the equation  (\ref{linbcd}) at $b$ is equal to $(\sqrt{3}(\zeta - \varkappa) + 6 \eta^2)t^4$. The equation (\ref{proto}) then has the form $(5!:-4!\alpha:-3!\beta) 
\equiv (\frac{1}{2}: 4 \eta: \sqrt{3}(\zeta-\varkappa)-6\eta^2)$. 
This gives $\eta = -\frac{1}{40}\alpha$ and $\zeta-\varkappa=\frac{\sqrt{3}\alpha^2}{800}-\frac{\beta}{40\sqrt{3}}$. For certainty, we can take $\varkappa=0$, and all higher terms in the series (\ref{pq}) to be trivial, so that the desired curve (\ref{pq}) has the form 
\begin{eqnarray}
\label{pqqa}
\left\{
\begin{aligned}
p(t) =  \frac{\sqrt{3}+1}{2} t - \frac{\alpha}{40} t^2 + 
\frac{\sqrt{3}\alpha^2}{800}-\frac{\beta}{40\sqrt{3}} t^3,  \\
q(t) = \frac{1-\sqrt{3}}{2} t - \frac{\alpha}{40} t^2  \ . \qquad \qquad \qquad \\
\end{aligned}
\right. \ .
\end{eqnarray}

Let us prove that any other algebraic curve in the space of $a_{23}$-configurations, that also approaches the same algebra $S(\alpha, \beta; A)$, is homotopic to this  curve (\ref{pqqa}). Let the starting segment of the power expansion of such a curve, containing all its terms of degree $\leq 3$, be 
\begin{eqnarray}
\label{pqq}
\left\{
\begin{aligned}
p(t) =  \frac{\sqrt{3}+1}{2} t - \frac{\alpha}{40} t^2 + \zeta_{w_1} t^{w_1}+  \ldots + \zeta_{w_N} t^{w_N} + \zeta_3 t^3, \quad  \\
\ \ q(t) = \frac{1-\sqrt{3}}{2} t - \frac{\alpha}{40} t^2 + \varkappa_{w_1} t^{w_1} + \ldots + \varkappa_{w_N} t^{w_N} + \varkappa_3 t^3 \ ,  \\
\end{aligned}
\right. \ 
\end{eqnarray}
where $2 < w_1 < \dots < w_N < 3$. 

Regardless of the coefficients $\zeta_{w_i}$ and $\varkappa_{w_i}$, it provides the leading (of degree 4) terms in the coefficients of the linear equation (\ref{linbcd}) at the variables $d$ and $c$ equal to $\frac{1}{2}t^4$ and $\frac{-\alpha}{10}t^4$ respectively. Then $\zeta_{w_1}-\varkappa_{w_1} =0$, because otherwise the leading term of the coefficient at $b$ is equal to $\sqrt{3}(\zeta_{w_1}-\varkappa_{w_1}) t^{1+w_1}$ and is stronger than all terms at $c$ and $d$. In the same way, the coefficients $\zeta_{w_i}$ and $\varkappa_{w_i}$  should coincide for any $i$, and the coefficients $\zeta_3$ and $\varkappa_3$ should satisfy the condition $\zeta_3-\varkappa_3=
\frac{\sqrt{3}\alpha^2}{800}-\frac{\beta}{40\sqrt{3}}.$

For any fixed set of exponents $w_i$, the space of polynomials of the form (\ref{pqq})  that satisfy these conditions is a contractible affine space, and we can deform the curve (\ref{pqq}) to the sample curve (\ref{pqqa}) inside this space. Also, we can remove all terms of degree greater than three, the space of which is also contractible and which do not participate in the restrictions on the curves approaching the cell $S$. Thus, the mod 2 incidence coefficient $[a_{2 3}, S]$ is equal to 1.

\begin{remark} \rm
The same computation can also be applied to the 
cells $a_{1 3}$ and $a_{2 4}$ instead of $a_{2 3}$. In these cases it is convenient to choose the origin respectively at the left and the right critical point. 
This computation yields the same values (\ref{solu}) for $\xi$ and $\lambda$. They contradict the shape of configurations of these classes. On the other hand, the contraction of supports of the algebras of the classes $a_{1 2}$, $a_{1 4}$, and $a_{3 4}$ necessarily gives us 
 the condition $f'(A)=f''(A)=f'''(A)=0$ on all functions of the limit algebras. This condition is not valid for almost all functions of an algebra of the class $S$. Thus, the incidence coefficients $[a_{i j}, S]$ are equal to zero for all $(i, j) \neq (2, 3)$.
\end{remark}

Similar (or simpler) considerations prove all formulas for the
boundary operators $\partial_i: C_i \to C_{i-1}$, $i = 2, 3, 4, 5, 6$, of the cell complex. These formulas are listed in the following
 subsections \ref{2to1}, \ref{3to2}, \ref{4to3}, \ref{5to4}, and \ref{6to5}, respectively. 
The matrices of all these operators are given in the arXiv version of this paper, arxiv: 2503.19343.

{\normalsize

\subsection{\bf $\partial_1: C_1 \to C_0$}
\label{1to0} 

The complex $\overline{CD}_3(S^1)$ contains only one 0-dimensional cell, therefore the boundary map $C_1 \to C_0$ is trivial.

\subsection{\bf $\partial_2: C_2 \to C_1$}
\label{2to1}

\begin{proposition}
\label{2to1prop}
The boundary map $C_2 \to C_1$ of the cell complex of the variety $\overline{CD}_3(S^1)$ with coefficients in ${\mathbb Z}_2$ defined by the cell decomposition described above is given by the following formulas:
\begin{align}
 \partial (\Upsilon_1) & = \nabla + \bar \Upsilon_1 + \bar \Upsilon_2 & 
 \partial (\Upsilon_2) & = \nabla + \bar \Upsilon_2 + \bar \Upsilon_1 \nonumber \\ 
 \partial (\Omega_1) & = \nabla + \bar \Omega_1 + \bar \Omega_2 & 
 \partial (\Omega_2) & = \nabla + \bar \Omega_2 + \bar \Omega_1 \nonumber \\ 
 \partial (\Lambda) & = \nabla & 
 \partial (\Theta) & = 0 \nonumber \\
 \partial (\bar U_{12}) & = \bar \Upsilon_1 + \bar \Lambda + \bar \Omega_2 + \bar \Theta & 
 \partial (\bar U_{13}) & = \bar \Omega_1 + \bar \Lambda + \bar \Upsilon_2 + \bar \Theta \nonumber \\ 
 \partial (\bar U_{21}) & = \bar \Upsilon_1 + \bar \Omega_2 + \bar \Lambda + \bar \Theta & 
 \partial (\bar U_{23}) & = \bar \Omega_1 + \bar \Upsilon_2 + \bar \Lambda + \bar \Theta \nonumber \\ 
 \partial (\bar U_{31}) & = \bar \Lambda + \bar \Omega_2 + \bar \Upsilon_2 + \bar \Theta & 
 \partial (\bar U_{32}) & = \bar \Lambda + \bar \Upsilon_2 + \bar \Omega_2 + \bar \Theta \nonumber \\
 \partial (\bar V_1) & = \bar \Omega_1 + \bar \Upsilon_1 + \bar \Omega_2 & 
 \partial (\bar V_2) & = \bar \Omega_1 + \bar \Omega_2 + \Upsilon_1 \nonumber \\ 
 \partial (\bar V_3) & = \bar \Upsilon_2 \nonumber \\
 \partial (\bar W_1) & = \bar \Omega_1 + \bar \Lambda + \bar \Omega_2 & 
 \partial (\bar W_2) & = \bar \Omega_1 + \bar \Omega_2 + \bar \Lambda \nonumber \\ 
 \partial (\bar W_3) & = \bar \Lambda &
 \partial (\bar X) & = \bar \Upsilon_1 \nonumber \\
 \partial (\bar Y_1) & = 0 &
 \partial (\bar Y_2) & = 0 \nonumber \\
 \partial (\bar Z_1^+) & = \bar \Omega_1 + \bar \Lambda + \bar \Theta & 
 \partial (\bar Z_2^+) & = \bar \Omega_2 + \bar \Lambda + \bar \Theta \nonumber \\
 \partial (\bar Z_1^-) & = \bar \Omega_1 + \bar \Lambda + \bar \Theta & 
 \partial (\bar Z_2^-) & = \bar \Omega_2 + \bar \Lambda + \bar \Theta \nonumber \\ 
 \partial (\bar T_1) & = 0 & 
 \partial (\bar T_2) & = 0 \nonumber \\ 
 \partial (\bar J) & = 0 & 
\partial (\bar S) & = 0 \nonumber
\end{align} 
\end{proposition}

\begin{corollary}
The kernel of the operator $\partial_2: C_2 \to C_1$ is generated by 23 linearly independent cycles $($\ref{52}$)$--$($\ref{74}$)$:
\begin{align}
 & \Upsilon_1+\Upsilon_2 \label{52} \\
& \Omega_1+\Omega_2 \label{53} \\ 
& \Theta \label{54} \\
& \bar U_{12}+ \bar U_{21} \label{55} \phantom{-----------} \\
& \bar U_{13}+\bar U_{23} \label{56} \\
& \bar U_{31} + \bar U_{32} \label{57} \\
& \bar V_1 + \bar V_2 \label{58} \\
& \bar W_1 + \bar W_2 \label{59} 
\end{align}
\begin{align}
& \bar Y_1 \label{60}\\
& \bar Y_2 \label{61} \\
& \bar Z_1^+ + \bar Z_1^- \label{62} \\
& \bar Z_2^+ + \bar Z_2^- \label{63} \\
& \bar T_1 \label{64} \\
& \bar T_2 \label{65} \\ 
& \bar J \label{66} \\
& \bar S \label{67} \\
& \bar U_{12} + \bar X + \bar Z_2^+ \label{68} \\
& \bar U_{31} + \bar V_3 + \bar Z_2^+ \label{69} \\
& \bar U_{13} + \bar V_3 + \bar Z_1^+ \label{70} \\
& \Upsilon_1 + \bar X + \bar V_3 + \Lambda \label{71} \\
& \bar V_1 + \Omega_1 + \Lambda + \bar X \label{72} \\
 & \bar Z_1^+ + \bar Z_2^+ + \bar W_1 + \bar W_3 \label{73} \\
& \bar \Omega_1 + \bar \Lambda + \bar Z_1^+ + \bar Z_2^+ \phantom{------}
\label{74} \end{align}
\end{corollary}

\subsection{\bf $\partial_3: C_3 \to C_2$}
\label{3to2}

\begin{proposition}
\label{3to2prop}
The boundary map $C_3 \to C_2$ of the cell complex of the variety $\overline{CD}_3(S^1)$ with coefficients in ${\mathbb Z}_2$ is given $($in terms of both the generators of $C_2$ described in \S \ref{twodim} and the generators $($\ref{52}$)$--$($\ref{74}$)$ of the subgroup $\ker \partial_2 \subset C_2)$ by the following formulas:
\begin{eqnarray}
 \partial (J) & = & 0 \nonumber \\ 
 \partial(S) & = & 0 \nonumber\\
 \partial (T_1) & = & \Theta + \bar T_1 + \bar T_2 = (\ref{54})+(\ref{64})+(\ref{65}) \nonumber \\ 
 \partial (T_2) & = & \Theta + \bar T_2 + \bar T_1 = (\ref{54})+(\ref{64})+(\ref{65}) \nonumber \\
 \partial (U_{12}) & = & \Upsilon_1 + \Lambda + \bar U_{12} + \bar U_{31} = (\ref{68})+(\ref{69})+(\ref{71}) \nonumber \\
 \partial (U_{13}) & = & \Omega_1 + \Lambda + \Theta + \bar U_{13} + \bar U_{32} = (\ref{54})+(\ref{57})+(\ref{69})+(\ref{70})+(\ref{74}) \nonumber \\
 \partial (U_{21}) & = & \Upsilon_1 + \Omega_2 + \bar U_{21} + \bar U_{13} = (\ref{53})+(\ref{55})+(\ref{68})+(\ref{70})+(\ref{71})+(\ref{74}) \nonumber \\
 \partial (U_{23}) & = & \Omega_1 + \Upsilon_2 + \bar U_{23} + \bar U_{12} = (\ref{52})+(\ref{56})+(\ref{68})+(\ref{70})+(\ref{71})+(\ref{74}) \nonumber 
\end{eqnarray}
\begin{eqnarray}
 \partial (U_{31}) & = & \Lambda + \Omega_2 + \Theta + \bar U_{31} + \bar U_{23}= \nonumber\\
\nonumber & = & (\ref{53})+(\ref{54})+(\ref{56})+ (\ref{69})+(\ref{70})+(\ref{74}) \nonumber \\
 \partial (U_{32}) & = & \Lambda + \Upsilon_2 + \bar U_{32} + \bar U_{21} = (\ref{52})+(\ref{55})+(\ref{57})+(\ref{68})+(\ref{69})+(\ref{71}) \nonumber \\
 \partial (V_1) & = & \Omega_1 + \Upsilon_1 + \bar V_1 + \bar V_3 = (\ref{71})+(\ref{72}) \nonumber \\
 \partial (V_2) & = & \Omega_1 + \Omega_2 + \bar V_2 + \bar V_1 = (\ref{53})+(\ref{58}) \nonumber\\ 
 \partial (V_3) & = & \Upsilon_2 + \Omega_2 + \bar V_3 + \bar V_2 = (\ref{52}) + (\ref{53}) + (\ref{58}) + (\ref{71}) + (\ref{72}) \nonumber\\ 
 \partial (W_1) & = & \Omega_1 + \Lambda + \bar W_1 + \bar W_3 =(\ref{73})+(\ref{74}) \nonumber\\ 
 \partial (W_2) & = & \Omega_1 + \Omega_2 + \bar W_2 + \bar W_1 = (\ref{53})+(\ref{59}) \nonumber\\ 
 \partial (W_3) & = & \Lambda + \Omega_2 + \bar W_3 + \bar W_2 = (\ref{53})+(\ref{59})+(\ref{73})+(\ref{74}) \nonumber\\
 \partial (X) & = & \Upsilon_1 + \Upsilon_2 = (\ref{52})\nonumber\\ 
 \partial (Y_1) & = & \Theta + \bar Y_1 + \bar Y_2 = (\ref{54}) + (\ref{60})+(\ref{61})\nonumber\\ 
 \partial (Y_2) & = & \Theta + \bar Y_2 + \bar Y_1 = (\ref{54})+(\ref{60})+(\ref{61}) \nonumber\\ 
 \partial (Z_1^+) & = & \Omega_1 + \Lambda + \bar Z_1^+ + \bar Z_2^+ = (\ref{74}) \nonumber\\ 
 \partial (Z_2^+) & = & \Omega_2 + \Lambda + \Theta + \bar Z_2^+ + \bar Z_1^+= (\ref{53})+(\ref{54})+(\ref{74})\nonumber\\ 
 \partial (Z_1^-) & = & \Omega_1 + \Lambda + \Theta + \bar Z_1^- + \bar Z_2^-= (\ref{54})+(\ref{62})+(\ref{63}) + (\ref{74})\nonumber\\ 
 \partial (Z_2^-) & = & \Omega_2 + \Lambda + \bar Z_2^- + \bar Z_1^- = (\ref{53})+(\ref{62})+(\ref{63})+(\ref{74})\nonumber \\
 \partial (\bar a_{12}) & = & \bar V_1 + \bar U_{23} + \bar X + \bar U_{32} + \bar S = \nonumber\\
\nonumber & = & (\ref{56})+(\ref{57})+(\ref{67})+(\ref{69})+(\ref{70})+(\ref{72})+(\ref{74}) \nonumber\\ 
 \partial (\bar a_{13}) & = & \bar U_{13} + \bar U_{23} + \bar U_{32} + \bar U_{31} = (\ref{56})+(\ref{57}) \nonumber\\ 
 \partial (\bar a_{14}) & = & \bar U_{12} + \bar X + \bar U_{32} + \bar V_3 + \bar S = (\ref{57})+(\ref{67})+(\ref{68})+(\ref{69}) \nonumber \\
 \partial (\bar a_{23}) & = & \bar U_{13}+ \bar V_2 + \bar U_{31} + \bar X + \bar S = \nonumber\\
\nonumber & = & (\ref{58})+(\ref{67})+(\ref{69}) + (\ref{70})+(\ref{72})+(\ref{74}) \nonumber\\
 \partial (\bar a_{24}) & = & \bar U_{12}+ \bar Y_1+ \bar U_{21}+ \bar U_{31} + \bar U_{32} + \bar Y_2 = (\ref{55})+(\ref{57})+(\ref{60})+(\ref{61})\nonumber \\
 \partial (\bar a_{34}) & = & \bar X+ \bar U_{21}+ \bar V_3 + \bar U_{31} + \bar S = (\ref{55}) + (\ref{67})+(\ref{68})+(\ref{69}) \nonumber\\
 \partial (\bar b_{12} ) & = & \bar V_1 + \bar W_1 + \bar U_{12} + \bar W_3 + \bar Z_2^+ = (\ref{68})+(\ref{72})+(\ref{73})+(\ref{74}) \nonumber \\
 \partial (\bar b_{13}) & = & \bar W_1 + \bar Z_1^+ + \bar W_3 + \bar Z_2^- + \bar T_2 = (\ref{63})+(\ref{65})+(\ref{73}) \nonumber \\
 \partial (\bar b_{14}) & = & \bar U_{13}+ \bar Z_1^- + \bar V_3 = (\ref{62}) + (\ref{70}) \nonumber\\ 
 \partial (\bar b_{21}) & = & \bar V_1 + \bar W_2 + \bar U_{21} + \bar W_1 + \bar Z_1^- = \nonumber \\
\nonumber & = &(\ref{55})+(\ref{59})+(\ref{62})+(\ref{68})+(\ref{72})+(\ref{74})\nonumber \\
 \partial (\bar b_{23}) & = & \bar W_1 + \bar V_2 + \bar W_2 + \bar Z_1^+ + \bar U_{12} = (\ref{58}) + (\ref{59}) + (\ref{68})+(\ref{72})+(\ref{74}) \nonumber\\ 
 \partial (\bar b_{24}) & = & \bar Z_1^- + \bar T_1 + \bar Z_1^+ = (\ref{62})+(\ref{64}) \nonumber\\
 \partial (\bar b_{31}) & = & \bar W_3 + \bar Z_2^+ + \bar T_2 + \bar W_2 + \bar Z_1^- = (\ref{59})+(\ref{62})+(\ref{65})+(\ref{73}) \nonumber \\
 \partial (\bar b_{32}) & = & \bar W_2 + \bar V_2 + \bar W_3 + \bar Z_2^- + \bar U_{21} = \nonumber\\
\nonumber & = &(\ref{55})+(\ref{58})+(\ref{59})+(\ref{63})+(\ref{68})+(\ref{72})+(\ref{73})+(\ref{74}) \nonumber \\
 \partial (\bar b_{34}) & = & \bar U_{23} + \bar V_3 + \bar Z_1^+ = (\ref{56})+(\ref{70}) \nonumber \\
 \partial (\bar b_{41}) & = & \bar U_{31} + \bar Z_2^+ + \bar V_3 =(\ref{69}) \nonumber 
\end{eqnarray}
\begin{eqnarray}
 \partial (\bar b_{42}) & = & \bar Z_2^- + \bar Z_2^+ + \bar T_2 = (\ref{63})+(\ref{65}) \nonumber \\
 \partial (\bar b_{43}) & = & \bar U_{32} + \bar V_3 + \bar Z_2^- = (\ref{57})+(\ref{63})+(\ref{69}) \nonumber \\ 
 \partial (\bar c_1) & = & \bar W_1 + \bar U_{23} + \bar U_{32} + \bar W_3 + \bar T_2 = \nonumber\\
\nonumber & = &(\ref{56})+(\ref{57})+(\ref{65})+(\ref{69})+(\ref{70})+(\ref{73}) \nonumber\\ 
 \partial (\bar c_2) & = & \bar W_1 + \bar W_2 + \bar U_{31} + \bar T_1 + \bar U_{32} = (\ref{57})+(\ref{59})+(\ref{64}) \nonumber\\ 
 \partial (\bar c_3) & = & \bar U_{13} + \bar W_2 + \bar W_3 + \bar T_2 + \bar U_{31} = (\ref{59})+(\ref{65})+(\ref{69})+(\ref{70})+(\ref{73})\nonumber\\ 
 \partial (\bar c_4) & = & \bar U_{12} + \bar U_{21} + \bar T_2 = (\ref{55})+(\ref{65}) \nonumber\\ 
 \partial (\bar d) & = & \bar W_1 + \bar W_2 = (\ref{59}) \nonumber\\
 \partial (\bar e_1^+) & = & \bar Z_1^+ + \bar U_{23} + \bar U_{32} + \bar Z_2^- = (\ref{56})+(\ref{57})+(\ref{63})+(\ref{69})+(\ref{70}) \nonumber\\
 \partial (\bar e_2^+) & = & \bar Z_1^- + \bar Z_2^+ + \bar U_{13} + \bar U_{31} = (\ref{62})+(\ref{69})+(\ref{70}) \nonumber\\ 
 \partial (\bar e_3^+) & = & \bar Z_2^- + \bar U_{12} + \bar U_{21} + \bar Z_2^+ = (\ref{55})+(\ref{63}) \nonumber\\ 
 \partial (\bar e_1^-) & = & \bar Z_1^- + \bar U_{23} + \bar U_{32} + \bar Y_2 + \bar Z_2^+ + \bar J = \nonumber\\
\nonumber & = &(\ref{56})+(\ref{57})+(\ref{61})+(\ref{62})+(\ref{66})+(\ref{69})+(\ref{70}) \nonumber \\
 \partial (\bar e_2^-) & = & \bar Z_1^+ + \bar Z_2^- + \bar U_{13} + \bar U_{31} + \bar J + \bar Y_2 = \nonumber\\
\nonumber & = & (\ref{61})+(\ref{63}) + (\ref{66})+(\ref{69}) +(\ref{70}) \nonumber\\ 
 \partial (\bar e_3^-) & = & \bar Z_2^+ + \bar U_{12} + \bar U_{21} + \bar Y_1 + \bar Z_2^- + \bar J = (\ref{55})+(\ref{60})+(\ref{63})+(\ref{66}) \nonumber\\
 \partial (\bar k_1^+) & = & \bar Z_1^+ + \bar W_2 + \bar W_3 + \bar Z_2^- = (\ref{59})+(\ref{63})+(\ref{73}) \nonumber\\
 \partial (\bar k_2^+) & = & \bar Z_1^- + \bar Z_2^+ + \bar W_1 + \bar W_3 = (\ref{62})+(\ref{73}) \nonumber\\ 
 \partial (\bar k_3^+) & = & \bar Z_2^- + \bar W_1 + \bar W_2 + \bar Z_2^+ = (\ref{59})+(\ref{63}) \nonumber\\
 \partial (\bar k_1^-) & = & \bar Z_1^- + \bar W_2 + \bar W_3 + \bar T_2 + \bar Z_2^+ = (\ref{59})+(\ref{62})+(\ref{65})+(\ref{73})\nonumber\\ 
 \partial (\bar k_2^-) & = & \bar Z_1^+ + \bar Z_2^- + \bar W_1 + \bar W_3 + \bar T_2 = (\ref{63})+(\ref{65})+(\ref{73}) \nonumber\\ 
 \partial (\bar k_3^-) & = & \bar Z_2^+ + \bar W_1 + \bar W_2 + \bar T_1 + \bar Z_2^- = (\ref{59})+ (\ref{63})+(\ref{64}) \nonumber \\
 \partial (\bar g_1) & = & \bar T_1 + \bar Y_1 + \bar T_2 + \bar J = (\ref{60})+(\ref{64})+(\ref{65})+(\ref{66}) \nonumber\\ 
 \partial (\bar g_2) & = & \bar T_1 + \bar T_2 + \bar J + \bar Y_1= (\ref{60}) + (\ref{64}) + (\ref{65})+(\ref{66}) \nonumber\\ 
 \partial (\bar g_3) & = & \bar Y_2 + \bar J = (\ref{61})+(\ref{66})\nonumber\\ 
 \partial (\bar h^+) & = & \bar Z_1^+ + \bar Z_1^- + \bar Z_2^+ + \bar Z_2^- = (\ref{62})+(\ref{63}) \nonumber\\ 
 \partial (\bar h^-) & = & \bar Z_1^+ + \bar Z_1^- + \bar Z_2^+ + \bar Z_2^- = (\ref{62})+(\ref{63}) \nonumber\\ 
 \partial (\bar i_1^+) & = & \bar Z_1^+ + \bar Z_1^- + \bar T_1 = (\ref{62})+(\ref{64}) \nonumber\\ 
 \partial (\bar i_1^-) & = & \bar Z_1^+ + \bar Z_1^- + \bar T_1 = (\ref{62})+(\ref{64}) \nonumber\\ 
 \partial (\bar i_2^+) & = & \bar Z_2^+ + \bar Z_2^- + \bar T_2 = (\ref{63})+(\ref{65})\nonumber\\ 
 \partial (\bar i_2^-) & = & \bar Z_2^+ + \bar Z_2^- + \bar T_2 = (\ref{63})+(\ref{65}) \nonumber
\end{eqnarray}
\end{proposition} 

\begin{corollary}
The kernel of the map $\partial: C_3 \to C_2$ is generated by 46 linearly independent cycles $($\ref{6}$)$--$($\ref{51}$):$
\begin{align}
 & { J} \label{6} \\
& { S} \phantom{----------------------} \label{7} \\
& { T_1 + T_2} \label{8} \\
& { T_2 + Z_1^+ + Z_1^- + \bar h^+ + \bar g_1 + \bar e_3^- + \bar e_3^+} \label{9} \\
& { U_{12} + U_{32} + X + \bar e_3^+ +\bar c_2 + \bar k_3^- } \label{10} \\
& { U_{13}+ Z_1^- + \bar h^+ + \bar d +\bar c_1 + \bar c_3 + \bar b_{41} + \bar b_{34} } \label{11}\\
& { U_{21} + V_3 + X  + \bar a_{23} + \bar a_{34}  } \label{12}\\
& U_{23} + W_2 + V_3 + \bar b_{23} + \bar b_{34} \label{13}\\
& U_{31} + W_2 +  Z_1^- +  \bar b_{41} + \bar b_{34} +  \bar d +  \bar h^+ \label{14} \\
& U_{32}+ V_3 + W_2 + Z_1^+ + \bar b_{32} + \bar b_{43} + \bar k_1^+ + \bar k_3^+ \label{15}\\
& V_1 + W_2 + V_3 + X + \bar a_{23} + \bar a_{34} + \bar b_{14} + \bar b_{21} \label{16} \\
& V_2 +  W_2 +  \bar a_{23} + \bar a_{34} +  \bar b_{14} + \bar b_{21} \label{17} \\
& { W_1+W_2+W_3} \label{18} \\
& W_1 + Z_1^+ + \bar k_1^+ + \bar k_3^+ \label{19} \\
& { W_2 + Z_2^+  + Z_1^-  + \bar d  + \bar h^+ } \label{20}\\
& { Y_1 + Y_2} \label{21}\\
& { Y_2 + Z_1^- + Z_1^+ + \bar h^+ + \bar g_3 + \bar e_3^- + \bar e_3^+ } \label{22} \\
& { Z_1^+ + Z_1^- + Z_2^+ + Z_2^- } \label{23} \\
& { \bar a_{12} + \bar b_{24} + \bar k_3^- + \bar b_{21} + \bar b_{34} + \bar a_{14} + \bar e_3^+ } \label{24} \\
& { \bar a_{13} + \bar e_1^+ + \bar e_2^+ + \bar h^+ } \label{25}\\
& { \bar a_{14} + \bar e_3^+ + \bar a_{34} + \bar c_2 + \bar k_3^- } \label{26} \\
& { \bar a_{23} + \bar b_{14} + \bar b_{23} + \bar b_{24} + \bar a_{34} + \bar e_3^+ + \bar k_3^- } \label{27}\\
& \bar a_{24} + \bar g_3 + \bar e_3^- + \bar c_2 + \bar k_3^- \label{28} \\
& { \bar b_{12} + \bar b_{21} + \bar b_{24} + \bar b_{34} + \bar b_{41} + \bar c_2 + \bar e_1^+ + \bar e_3^+ + \bar k_1^+ + \bar k_3^+ } \label{29}\\
& { \bar b_{13} + \bar k_2^-} \label{30}\\
& { \bar b_{14} + \bar c_2 + \bar c_1 + \bar c_3 +
\bar b_{24} + \bar b_{34} } \label{31} \\
& \bar b_{23} + \bar b_{32} + \bar k_1^+ + \bar e_3^+ + \bar k_3^+ \label{32} \\
& { \bar b_{24} + \bar i_1^+} \label{33}\\
& { \bar b_{31} + \bar k_1^- } \label{34}\\
& { \bar b_{34} + \bar e_1^+ + \bar b_{43} }\label{35}
\end{align}
\begin{align}
& { \bar b_{41} + \bar b_{43} + \bar c_2 + \bar k_3^- } \label{36}\\ 
& { \bar b_{42} + \bar i_2^+ } \label{37} \\
& { \bar c_1 + \bar k_2^- + \bar e_1^+ }\label{38}\\
& { \bar c_3 + \bar k_2^- + \bar e_2^+ + \bar h^+ + \bar d} \label{39} \\
& \bar c_4 + \bar e_3^+ + \bar k_1^- + \bar k_3^+ + \bar d + \bar k_1^+ + \bar h^+ \label{40} \\
& { \bar d + \bar i_1^+ + \bar k_3^- + \bar h^+ }\label{41} \\
& { \bar e_1^+ + \bar e_1^- + \bar g_3 + \bar h^+ }\label{42} \\
& { \bar e_2^+ + \bar e_2^- + \bar g_3 + \bar h^+} \label{43} \\
& \bar e_3^+ + \bar e_3^- + \bar i_1^+ + \bar h^+ +
\bar i_2^+ + \bar g_1 \label{44}\\
& { \bar k_1^+ + \bar k_2^- + \bar k_3^+ +\bar i_2^+ } \label{45}\\
& { \bar k_1^- + \bar k_2^+ + \bar k_3^+ + \bar i_2^+ } \label{46}\\
& { \bar k_2^- + \bar k_2^- + \bar k_3^+ + \bar k_3^- + \bar i_1^+ + \bar i_2^+} \label{47}\\
& { \bar g_1 + \bar g_2} \label{48}\\
& { \bar h^+ + \bar h^-} \label{49}\\
& { \bar i_1^+ + \bar i_1^-} \label{50}\\
& { \bar i_2^+ + \bar i_2^-} \label{51}
\end{align}
\end{corollary}

\subsection{$\partial_4: C_4 \to C_3$}
\label{4to3}

\begin{proposition}
\label{4to3prop}
The boundary operator $\partial_4: C_4 \to C_3 $ is defined by the following formulas $($both in terms of the standard generators of the group $C_3$ and the generators $($\ref{6}$)$---$($\ref{51}$)$ of its subgroup $\ker \partial_3)$:
\begin{eqnarray}
 \partial (a_{12}) & = & V_1 + U_{23} + X + \bar a_{12} + \bar a_{14} = (\ref{13})+(\ref{16})+(\ref{24})+(\ref{27}) \nonumber \\ 
 \partial (a_{13}) & = & U_{13}+U_{23}+U_{32}+Y_2 +\bar a_{13} + \bar a_{24} = (\ref{11})+(\ref{13})+ \nonumber \\
& & + (\ref{15})+(\ref{22})+(\ref{25})+(\ref{28})+(\ref{32}) +(\ref{36})+(\ref{38})+(\ref{39}) \nonumber\\ 
 \partial ( a_{14}) & = & U_{12}+X+ U_{32} + \bar a_{14} + \bar a_{34}= (\ref{10})+(\ref{26})\nonumber\\ 
 \partial ( a_{23}) & = & U_{13}+V_2 + U_{31} + S + \bar a_{23} + \bar a_{12}=\nonumber\\
\nonumber & = & (\ref{7})+(\ref{11})+(\ref{14})+(\ref{17})+(\ref{24})+(\ref{26})+(\ref{31})\nonumber\\ 
 \partial ( a_{24}) & = & U_{12}+Y_1+U_{21}+U_{31} + \bar a_{24} + \bar a_{13}=\nonumber \\
\nonumber & = &(\ref{10})+(\ref{12})+(\ref{14})+(\ref{15})+(\ref{21})+(\ref{22})+ (\ref{25})+\nonumber\\
\nonumber & & +(\ref{27})+(\ref{28})+(\ref{31})+ (\ref{32})+(\ref{36})+(\ref{38})+(\ref{39})\nonumber \\
 \partial ( a_{34}) & = & X+ U_{21}+ V_3 + \bar a_{34} + \bar a_{23} = (\ref{12}) \nonumber 
\end{eqnarray}
\begin{eqnarray} 
 \partial ( b_{12} ) & = & V_1 + W_1 + U_{12} + \bar b_{12} + \bar b_{41}=\nonumber \\
\nonumber & = & (\ref{10})+(\ref{15})+(\ref{16})+(\ref{19})+(\ref{27})+ (\ref{29})+(\ref{32})+(\ref{35}) \nonumber\\ 
 \partial ( b_{13}) & = & W_1 + Z_1^+ + \bar b_{13} + \bar b_{42} = 
(\ref{19})+(\ref{30})+(\ref{37})+(\ref{45}) \nonumber\\ 
 \partial ( b_{14}) & = & U_{13}+ Z_1^- + \bar b_{14} + \bar b_{43} = (\ref{11})+(\ref{31})+(\ref{33})+(\ref{36})+(\ref{41}) \nonumber\\
\partial ( b_{21}) & = & V_1 + W_2 + U_{21} + \bar b_{21} + \bar b_{14} = (\ref{12})+(\ref{16}) \nonumber \\
 \partial (b_{23}) & = & W_1 + V_2 + W_2 + Z_1^+ + \bar b_{23} + \bar b_{12} =\nonumber \\
\nonumber & = & (\ref{17})+(\ref{19})+(\ref{27})+(\ref{29})+(\ref{35})+(\ref{36}) \nonumber \\
 \partial ( b_{24}) & = & W_1 + Z_1^- + T_1 + \bar b_{24} + \bar b_{13} =\nonumber\\
\nonumber & = & (\ref{8})+(\ref{9})+(\ref{19})+(\ref{30})+(\ref{33})+(\ref{44})+(\ref{45}) \nonumber \\
 \partial ( b_{31}) & = & W_3 + Z_2^+ + T_2 + \bar b_{31} + \bar b_{24} 
= (\ref{9})+(\ref{18})+(\ref{19})+ \nonumber \\
& & + (\ref{20})+(\ref{33})+(\ref{34}) + (\ref{41})+(\ref{44})+(\ref{45}) +(\ref{46})+(\ref{47}) \nonumber \\
 \partial ( b_{32}) & = & W_2 + V_2 + W_3 +Z_2^- + \bar b_{32} + \bar b_{21} =\nonumber\\ 
\nonumber &=& (\ref{17})+(\ref{18})+(\ref{19})+(\ref{20})+(\ref{23})+ (\ref{27})+(\ref{32})+(\ref{33})+(\ref{41}) \nonumber \\
 \partial ( b_{34}) & = & U_{23} + W_2 + V_3 + \bar b_{34} + \bar b_{23} = (\ref{13}) \nonumber\\
 \partial ( b_{41}) & = & U_{31} + Z_2^+ +\bar b_{41} +\bar b_{34}= (\ref{14})+(\ref{20}) \nonumber \\
 \partial ( b_{42}) & = & W_3 + Z_2^- + \bar b_{42} + \bar b_{31}= 
 (\ref{18})+(\ref{19})+(\ref{20})+ \nonumber \\
& & + (\ref{23})+(\ref{34})+ (\ref{37})+(\ref{41})+(\ref{45})+(\ref{46})+(\ref{47}) \nonumber\\
\partial ( b_{43}) & = & U_{32} + W_3 + V_3 + \bar b_{43} + \bar b_{32} = (\ref{15})+(\ref{18})+(\ref{19}) \nonumber \\
 \partial (c_1) & = & W_1 + U_{23} + U_{32} + \bar c_1 + \bar c_4 
= (\ref{13})+(\ref{15})+(\ref{19})+ \nonumber \\
& & + (\ref{32})+(\ref{35})+ (\ref{38})+(\ref{40})+ (\ref{41})+(\ref{46})+(\ref{47}) \nonumber \\
 \partial (c_2) & = & W_1 + W_2 + U_{31} + T_1 + \bar c_2 + \bar c_1= 
 (\ref{8})+(\ref{9})+(\ref{14})+ \nonumber \\
& & + (\ref{19})+(\ref{35})+(\ref{36})+(\ref{38}) + (\ref{41})+(\ref{44})+(\ref{45}) \nonumber\\
 \partial (c_3) & = & U_{13} + W_2 + W_3 + T_2 + \bar c_3 + \bar c_2 =
 (\ref{9})+(\ref{11})+(\ref{18})+\nonumber \\
& & + (\ref{19})+(\ref{35})+(\ref{36})+(\ref{38}) + (\ref{41})+(\ref{44})+(\ref{45}) \nonumber\\
 \partial (c_4) & = & U_{12} + U_{21} + W_3 + \bar c_4 + \bar c_3 =\nonumber\\
\nonumber &=&(\ref{10})+(\ref{12})+(\ref{15})+(\ref{18})+(\ref{19})+(\ref{27})+(\ref{31})+\nonumber\\
\nonumber & & +(\ref{32})+(\ref{35})+(\ref{38})+ (\ref{40})+(\ref{41})+(\ref{46})+(\ref{47}) \nonumber\\
 \partial (d) & = & W_1 + W_2 + W_3 = (\ref{18}) \nonumber\\
 \partial (e_1^+) & = & Z_1^+ + U_{23} + U_{32} + \bar e_1^+ + \bar e_3^+ = (\ref{13})+(\ref{15})+(\ref{32})+(\ref{35}) \nonumber\\
 \partial (e_2^+) & = & Z_1^- + Z_2^+ + U_{13} + U_{31} + \bar e_2^+ + \bar e_1^+ =\nonumber\\
\nonumber &=& (\ref{11})+(\ref{14})+(\ref{20})+(\ref{38})+(\ref{39}) \nonumber\\
 \partial (e_3^+) & = & Z_2^- + U_{12} + U_{21} + \bar e_3^+ + \bar e_2^+ =
 (\ref{10})+(\ref{12})+(\ref{15})+\nonumber \\
& & + (\ref{20})+(\ref{23})+(\ref{27})+ (\ref{31})+(\ref{32})+(\ref{35})+(\ref{38})+(\ref{39}) \nonumber \\
 \partial (e_1^-) & = & Z_1^- + U_{23} + U_{32} + Y_2 + \bar e_1^- + \bar e_3^- =\nonumber\\
\nonumber &=& (\ref{13})+(\ref{15})+(\ref{22})+(\ref{32})+(\ref{35})+(\ref{42}) \nonumber
\end{eqnarray}
\begin{eqnarray} 
 \partial (e_2^-) & = & Z_1^+ + Z_2^- + U_{13} + U_{31} + J + \bar e_2^- + \bar e_1^- =\nonumber \\
\nonumber &=& (\ref{6})+(\ref{11})+(\ref{14})+(\ref{20})+(\ref{23})+(\ref{38})+ (\ref{39})+(\ref{42})+(\ref{43}) \nonumber \\
 \partial (e_3^-) & = & Z_2^+ + U_{12} + U_{21} + Y_1 + \bar e_3^- + \bar e_2^- =\nonumber\\
\nonumber &=& (\ref{10})+(\ref{12})+(\ref{15})+(\ref{20})+(\ref{21})+(\ref{22})+(\ref{27})+\nonumber \\
\nonumber & & + (\ref{31})+(\ref{32})+(\ref{35})+(\ref{38})+(\ref{39})+(\ref{43}) \nonumber\\
 \partial (k_1^+) & = & Z_1^+ + W_2 + W_3 + \bar k_1^+ + \bar k_3^+ = (\ref{18}) +(\ref{19}) \nonumber\\
 \partial (k_2^+) & = & Z_1^- + Z_2^+ + W_1 + W_3 + \bar k_2^+ + \bar k_1^+ = \nonumber\\
\nonumber &=& (\ref{18}) +(\ref{20}) +(\ref{41}) +(\ref{45}) + (\ref{47}) \nonumber \\
 \partial (k_3^+) & = & Z_2^- + W_1 + W_2 + \bar k_3^+ + \bar k_2^+ = \nonumber \\
\nonumber &=& (\ref{19}) +(\ref{20}) +(\ref{23})+ (\ref{41}) +(\ref{45})+ (\ref{47}) \nonumber \\
 \partial (k_1^-) & = & Z_1^- + W_2 + W_3 + T_2 + \bar k_1^- + \bar k_3^- =\nonumber\\ 
\nonumber &=& (\ref{9}) + (\ref{18}) +(\ref{19}) +(\ref{44})+ (\ref{45}) +(\ref{46}) +(\ref{47}) \nonumber\\
 \partial (k_2^-) & = & Z_1^+ + Z_2^- + W_1 + W_3 + \bar k_2^- + \bar k_1^- =\nonumber\\ 
\nonumber&=& (\ref{18})+(\ref{20})+(\ref{23}) +(\ref{41})+(\ref{46}) +(\ref{47})\nonumber \\
 \partial (k_3^-) & = & Z_2^+ + W_1 + W_2 + T_1 + \bar k_3^- + \bar k_2^- = \nonumber\\
\nonumber &=& (\ref{8})+(\ref{9}) +(\ref{19})+(\ref{20})+(\ref{41}) +(\ref{44}) +(\ref{45}) \nonumber\\
 \partial (g_1) & = &T_1 + Y_1 + \bar g_1 + \bar g_3 = (\ref{8})+(\ref{9})+(\ref{21})+(\ref{22}) \nonumber\\
 \partial (g_2) & = & T_1 + T_2 + J + \bar g_2 + \bar g_1 = (\ref{6})+(\ref{8})+(\ref{48}) \nonumber\\
 \partial (g_3) & = & Y_2 + T_2 + \bar g_3 + \bar g_2 = (\ref{9})+(\ref{22})+(\ref{48}) \nonumber\\
 \partial (h^+) & = & Z_1^+ + Z_1^- + Z_2^+ + Z_2^- + S = (\ref{7})+(\ref{23}) \nonumber\\
 \partial (h^-) & = & Z_1^+ + Z_1^- + Z_2^+ + Z_2^- + J = (\ref{6})+(\ref{23}) \nonumber\\
 \partial (i_1^+) & = & Z_1^+ + Z_1^- + T_1 + S + \bar i_1^+ + \bar i_2^+ = (\ref{7})+(\ref{8})+(\ref{9})+(\ref{44}) \nonumber\\
 \partial (i_1^-) & = & Z_1^+ + Z_1^- + T_1 + \bar i_1^- + \bar i_2^- =\nonumber\\ \nonumber &=& (\ref{8})+(\ref{9})+(\ref{44})+(\ref{50})+(\ref{51}) \nonumber\\
 \partial (i_2^+) & = & Z_2^+ + Z_2^- + T_2 + S + \bar i_2^+ + \bar i_1^+ = (\ref{7})+(\ref{9})+(\ref{23})+(\ref{44}) \nonumber\\
 \partial (i_2^-) & = & Z_2^+ + Z_2^- + T_2 + \bar i_2^- + \bar i_1^- =\nonumber\\ \nonumber &=& (\ref{9})+(\ref{23})+(\ref{44})+(\ref{50})+(\ref{51}) \nonumber\\
 \partial ( \bar B_1^1) & = & \bar b_{12} + \bar a_{12} + \bar c_1 + \bar a_{14} + \bar b_{43} + \bar i_2^+ =\nonumber\\ 
\nonumber &=& (\ref{24})+(\ref{29})+(\ref{36})+(\ref{38})+(\ref{45}) \nonumber\\
 \partial ( \bar B_1^2 ) & = & \bar b_{13} + \bar c_1 + \bar e_1^+ + \bar i_2^- + \bar b_{42} = (\ref{30})+(\ref{37})+(\ref{38})+(\ref{51}) \nonumber\\
\nonumber \partial ( \bar B_1^3) & = & \bar b_{14} + \bar a_{13} + \bar e_1^- + \bar b_{41} + \bar g_3 =\nonumber \\
&=& (\ref{25})+(\ref{31})+(\ref{33})+(\ref{35})+(\ref{36})+ (\ref{38})+(\ref{39})+(\ref{41})+(\ref{42}) \nonumber \\
 \partial ( \bar B_2^1 ) & = & \bar b_{12} + \bar b_{21} + \bar c_2 + \bar a_{24} + \bar g_1 + \bar c_1 + \bar e_1^- = \nonumber\\ 
\nonumber &=& (\ref{28})+(\ref{29})+(\ref{33})+(\ref{35})+(\ref{36})+(\ref{38})+ (\ref{42})+(\ref{44})+(\ref{45})\nonumber 
\end{eqnarray}
\begin{eqnarray}
 \partial ( \bar B_2^2 ) & = & \bar b_{13} + \bar b_{24} + \bar c_1 + \bar e_1^+ + \bar i_1^- = (\ref{30})+(\ref{33})+(\ref{38})+(\ref{50}) \nonumber\\
 \partial ( \bar B_2^3 ) & = & \bar b_{14} + \bar b_{23} + \bar a_{23} + \bar c_2 + \bar a_{14} + \bar i_1^+ = (\ref{26})+(\ref{27})+(\ref{33}) \nonumber \\
 \partial ( \bar B_3^1 ) & = & \bar a_{12} + \bar b_{21} + \bar b_{34} + \bar a_{34} + \bar c_2 + \bar i_1^+ = (\ref{24})+(\ref{26})+(\ref{33}) \nonumber \\
 \partial ( \bar B_3^2 ) & = & \bar b_{24} + \bar b_{31} + \bar c_3 + \bar e_2^+ + \bar i_1^- =\nonumber\\ 
\nonumber &=& (\ref{33})+(\ref{34})+(\ref{39})+(\ref{41})+(\ref{46})+(\ref{47})+(\ref{50}) \nonumber\\
 \partial ( \bar B_3^3 ) & = & \bar c_2 + \bar b_{23} + \bar b_{32} + \bar c_3 + \bar e_2^- + \bar g_2 + \bar a_{24} =\nonumber\\ 
\nonumber &=& (\ref{28})+(\ref{32})+(\ref{39})+(\ref{41})+(\ref{43})+ (\ref{44})+(\ref{45})+(\ref{48}) \nonumber\\
 \partial ( \bar B_4^1 ) & = & \bar a_{13} + \bar b_{34} + \bar b_{43} + \bar g_3 + \bar e_2^- = (\ref{25})+(\ref{35})+(\ref{43}) \nonumber\\ 
 \partial ( \bar B_4^2 ) & = & \bar b_{31} + \bar b_{42} + \bar c_3 + \bar e_2^+ + \bar i_2^- =\nonumber \\
\nonumber &=& (\ref{34})+(\ref{37})+(\ref{39})+(\ref{41})+(\ref{46})+(\ref{47})+(\ref{51}) \nonumber \\
\partial ( \bar B_4^3 ) & = & \bar c_3 + \bar a_{23} + \bar b_{32} + \bar b_{41} + \bar a_{34} + \bar i_2^+ =\nonumber\\ 
\nonumber &=& (\ref{27})+(\ref{31})+(\ref{32})+(\ref{35})+(\ref{36})+(\ref{38})+(\ref{45}) \nonumber \\
 \partial ( \bar B_5^1 ) & = & \bar a_{14} + \bar c_4 + \bar a_{34} + \bar b_{43} + \bar b_{41} + \bar i_2^+ =\nonumber\\
\nonumber &=& (\ref{26})+(\ref{36})+(\ref{40})+(\ref{41})+(\ref{45})+(\ref{46})+(\ref{47}) \nonumber \\
 \partial ( \bar B_5^2 ) & = & \bar c_4 + \bar e_3^+ + \bar i_2^- = (\ref{40})+(\ref{41})+(\ref{45})+(\ref{46})+(\ref{47}) +(\ref{51})\nonumber\\
 \partial ( \bar B_5^3 ) & = & \bar a_{24} + \bar b_{41} + \bar e_3^- + \bar b_{43} + \bar g_3 = (\ref{28})+(\ref{36}) \nonumber\\
 \partial ( \bar C_{12}) & = & \bar c_1 + \bar b_{34} + \bar b_{43} + \bar k_2^- = (\ref{35})+(\ref{38}) \nonumber\\
\partial ( \bar C_{13}) & = & \bar b_{42} + \bar k_3^+ + \bar k_2^+ + \bar k_1^- = (\ref{37})+(\ref{46}) \nonumber\\
 \partial ( \bar C_{14}) & = & \bar b_{24} + \bar k_2^+ + \bar k_3^- + \bar k_1^+ = (\ref{33})+(\ref{45})+(\ref{47}) \nonumber\\
 \partial ( \bar C_{15}) & = & \bar b_{23} + \bar b_{32} + \bar k_2^- + \bar c_4 = (\ref{32})+(\ref{40})+(\ref{41})+(\ref{46})+(\ref{47}) \nonumber\\
 \partial ( \bar C_{23}) & = & \bar c_2 + \bar b_{41} + \bar k_3^- + \bar b_{43} = (\ref{36}) \nonumber\\
 \partial ( \bar C_{24}) & = & \bar k_1^+ + \bar k_2^- + \bar k_3^+ + \bar b_{42} = (\ref{37})+(\ref{45}) \nonumber\\
 \partial ( \bar C_{25}) & = & \bar b_{31} + \bar k_1^- = (\ref{34}) \nonumber\\
 \partial ( \bar C_{34}) & = & \bar b_{14} + \bar c_3 + \bar k_1^- + \bar b_{41}= \nonumber\\ 
\nonumber &=& (\ref{31})+(\ref{33})+(\ref{35})+(\ref{36})+(\ref{38})+(\ref{46})+(\ref{47}) \nonumber\\
 \partial ( \bar C_{35}) & = & \bar b_{13} + k_2^- = (\ref{30}) \nonumber\\ 
\nonumber 
\partial ( \bar C_{45}) & = & \bar b_{12} + \bar b_{21} + \bar c_4 + \bar k_1^- = \nonumber \\
& = & (\ref{29})+(\ref{33})+(\ref{35})+(\ref{36})+(\ref{40})+(\ref{41}) \nonumber\\
 \partial ( \bar D_{12}^+) & = & \bar k_3^+ + \bar e_3^+ + \bar b_{12} + \bar b_{21} + \bar k_2^+ =\nonumber\\ 
\nonumber &=& (\ref{29})+(\ref{33})+(\ref{35})+(\ref{36})+(\ref{45})+(\ref{47}) \nonumber\\
 \partial ( \bar D_{12}^-) & = & \bar k_3^- + \bar e_3^- + \bar g_1 + \bar b_{12} + \bar b_{21} + \bar k_2^- + \bar h^- =\nonumber\\ 
\nonumber &=& (\ref{29})+(\ref{33})+(\ref{35})+(\ref{36})+(\ref{44})+(\ref{45})+(\ref{49}) \nonumber \\
 \partial ( \bar D_{13}^+) & = & \bar k_2^+ + \bar b_{13} + \bar b_{31} + \bar k_1^+ = (\ref{30})+(\ref{34})+(\ref{45})+(\ref{46}) \nonumber\\
 \partial ( \bar D_{13}^-) & = & \bar k_2^- + \bar b_{13} + \bar b_{31} + \bar k_1^- = (\ref{30})+(\ref{34}) \nonumber 
\end{eqnarray}
\begin{eqnarray}
\partial ( \bar D_{14}^+) & = & \bar e_2^+ + \bar b_{14} + \bar b_{41} =\nonumber \\
\nonumber &=& (\ref{31})+(\ref{33})+(\ref{35})+(\ref{36})+(\ref{38})+(\ref{39})+(\ref{41}) \nonumber \\
 \partial ( \bar D_{14}^-) & = & \bar e_2^- + \bar h^- + \bar b_{14} + \bar b_{41} + \bar g_3 = (\ref{31})+(\ref{33})+(\ref{35})+ \nonumber \\
& & + (\ref{36})+(\ref{38})+(\ref{39}) + (\ref{41})+(\ref{43})+(\ref{49}) \nonumber\\
 \partial ( \bar D_{23}^+) & = & \bar k_1^+ + \bar k_3^+ + \bar b_{23} + \bar b_{32} + \bar e_3^+ = (\ref{32}) \nonumber\\
\partial ( \bar D_{23}^-) & = & \bar k_1^- + \bar k_3^- + \bar g_2 + \bar h^- + \bar b_{23} + \bar b_{32} + \bar e_3^- =\nonumber \\ 
\nonumber &=& (\ref{32})+(\ref{44})+(\ref{45})+(\ref{46})+(\ref{47})+(\ref{48})+(\ref{49}) \nonumber \\
\partial ( \bar D_{24}^+) & = & \bar b_{24} + \bar b_{42} + \bar i_1^+ + \bar i_2^+ = (\ref{33})+(\ref{37}) \nonumber\\
 \partial ( \bar D_{24}^-) & = & \bar b_{24} + \bar b_{42} + \bar i_1^- + \bar i_2^- = (\ref{33})+(\ref{37})+(\ref{50})+(\ref{51}) \nonumber \\
\partial ( \bar D_{34}^+) & = & \bar e_1^+ + \bar b_{34} + \bar b_{43} = (\ref{35}) \nonumber\\
 \partial ( \bar D_{34}^-) & = & \bar e_1^- + \bar g_3 + \bar b_{34} + \bar b_{43} + \bar h^- = (\ref{35})+(\ref{42})+(\ref{49})\nonumber \\
 \partial (\bar E_+) & = & \bar k_1^- + \bar k_2^- + \bar k_3^- + \bar i_1^+ = (\ref{46})+(\ref{47})\nonumber \\
 \partial (\bar E_1) & = & \bar k_1^- + \bar k_2^+ + \bar k_3^+ + \bar i_2^- = (\ref{46})+(\ref{51}) \nonumber\\
 \partial (\bar E_2) & = & \bar k_1^+ + \bar k_2^- + \bar k_3^+ + \bar i_2^- = (\ref{45})+(\ref{51}) \nonumber\\ 
 \partial (\bar E_3) & = & \bar k_1^+ + \bar k_2^+ + \bar k_3^- + \bar i_1^- = (\ref{45})+(\ref{47})+(\ref{50}) \nonumber
\end{eqnarray} 
\end{proposition}

\subsection{\bf $\partial_5: C_5 \to C_4$}
\label{5to4}

\begin{proposition}
\label{5to4prop}
The boundary operator $\partial_5: C_5 \to C_4 $ is described by following formulas: 
\begin{eqnarray}
 \partial (B_1^1) & = & b_{12} + a_{12} + c_1 + a_{14} + \bar B_1^1 + \bar B_5^1 \nonumber \\ 
\partial (B_1^2 ) & = & b_{13} + c_1 + e_1^+ + \bar B_1^2 + \bar B_5^2\nonumber\\ 
\partial (B_1^3) & = & b_{14} + a_{13} + e_1^- + \bar B_1^3 + \bar B_5^3\nonumber\\ 
\partial (B_2^1 ) & = & b_{12} + b_{21} + c_2 + a_{24} + g_1 + \bar B_2^1 + \bar B_1^3\nonumber\\ 
\partial (B_2^2 ) & = & b_{13} + b_{24} + i_1^- + \bar B_2^2 + \bar B_1^2\nonumber\\ 
\partial (B_2^3 ) & = & b_{14} + b_{23} + a_{23} + c_2 + i_1^+ + \bar B_2^3 + \bar B_1^1 \nonumber \\
\partial (B_3^1 ) & = & a_{12} + b_{21} + b_{34} + a_{34} + \bar B_{3}^1 + \bar B_2^3\nonumber\\ 
\partial (B_3^2 ) & = & c_2 + b_{24} + b_{31} + c_3 + e_2^+ + \bar B_3^2 + \bar B_2^2\nonumber \\
\partial (B_3^3 ) & = & c_2 + b_{23} + b_{32} + c_3 + e_2^- + g_2 + \bar B_3^3 + \bar B_2^1\nonumber\\ 
\partial (B_4^1 ) & = & a_{13} + c_3 + b_{34} + b_{43} + g_3 + \bar B_4^1 + \bar B_3^3\nonumber \\
\partial (B_4^2 ) & = & b_{31} + b_{42} + i_2^- + \bar B_4^2 + \bar B_3^2 \nonumber\\ 
\partial (B_4^3 ) & = & c_3 + a_{23} + b_{32} + b_{41} + i_2^+ + \bar B_4^3 + \bar B_3^1\nonumber \\
\partial (B_5^1 ) & = & a_{14} + c_4 + a_{34} + b_{43} + \bar B_5^1 + \bar B_4^3\nonumber\\ 
\partial (B_5^2 ) & = & c_4 + b_{42} + e_3^+ + \bar B_5^2 + \bar B_4^2 \nonumber\\ 
\partial (B_5^3 ) & = & a_{24} + b_{41} + e_3^- + \bar B_5^3 + \bar B_4^1 \nonumber 
\end{eqnarray}
\begin{eqnarray} 
\partial ( C_{12}) & = & c_1 + d + b_{34} + b_{43} + \bar C_{12} + \bar C_{15} \nonumber\\ 
\partial ( C_{13}) & = & d + b_{42} + k_3^+ + \bar C_{13} + \bar C_{25}\nonumber\\ 
\partial ( C_{14}) & = & d + b_{24} + k_2^+ + k_3^- + \bar C_{14} + \bar C_{35}\nonumber\\ 
\partial ( C_{15}) & = & b_{23} + b_{32} + k_2^- + \bar C_{15} + \bar C_{45}\nonumber\\ 
\partial ( C_{23}) & = & c_2 + b_{41} + k_3^- + \bar C_{23} + \bar C_{12}\nonumber \\
\partial ( C_{24}) & = & k_1^+ + k_2^- + k_3^+ + \bar C_{24} + \bar C_{13}\nonumber \\
\partial ( C_{25}) & = & b_{31} + d + k_2^+ + k_1^- + \bar C_{25} + \bar C_{14}\nonumber \\
\partial ( C_{34}) & = & b_{14} + c_3 + k_1^- + \bar C_{34} + \bar C_{23}\nonumber\\ 
\partial ( C_{35}) & = & b_{13} + d + k_1^+ + \bar C_{35} + \bar C_{24}\nonumber\\ 
\partial ( C_{45}) & = & b_{12} + b_{21} + d + c_4 + \bar C_{45} + \bar C_{34}\nonumber \\
\partial (D_{12}^+) & = & k_3^+ + e_3^+ + b_{12} + b_{21} + \bar D_{12}^+ + \bar D_{14}^+\nonumber \\
 \partial (D_{12}^-) & = & k_3^- + e_3^- + g_1 + b_{12} + b_{21} + \bar D_{12}^- + \bar D_{14}^-\nonumber \\
\partial (D_{13}^+) & = & k_2^+ + h^+ + b_{13} + b_{31} + i_2^+ + \bar D_{13}^+ + \bar D_{24}^+\nonumber \\
\partial (D_{13}^-) & = & k_2^- + b_{13} + b_{31} + i_2^- + \bar D_{13}^- + \bar D_{24}^-\nonumber\\ 
\partial (D_{14}^+) & = & e_2^+ + b_{14} + b_{41} + \bar D_{14}^+ + \bar D_{34}^+\nonumber\\ 
\partial (D_{14}^-) & = & e_2^- + h^- + b_{14} + b_{41} + \bar D_{14}^- + \bar D_{34}^- \nonumber\\
\partial (D_{23}^+) & = & k_1^+ + k_3^+ + b_{23} + b_{32} + \bar D_{23}^+ + \bar D_{12}^+\nonumber\\ 
\partial (D_{23}^-) & = & k_1^- + k_3^- + g_2 + h^- + b_{23} + b_{32} + \bar D_{23}^- + \bar D_{12}^-\nonumber\\ 
\partial (D_{24}^+) & = & k_2^+ + h^+ + b_{24} + b_{42} + i_1^+ + \bar D_{24}^+ + \bar D_{13}^+\nonumber\\ 
\partial (D_{24}^-) & = & k_2^- + b_{24} + b_{42} + i_1^- + \bar D_{24}^- + \bar D_{13}^-\nonumber\\ 
\partial (D_{34}^+) & = & e_1^+ + k_1^+ + b_{34} + b_{43} + \bar D_{34}^+ + \bar D_{23}^+\nonumber\\ 
\partial (D_{34}^-) & = & e_1^- + k_1^- + g_3 + b_{34} + b_{43} + \bar D_{34}^- + \bar D_{23}^- \nonumber\\
\partial (E_+) & = & k_1^- + k_2^- + k_3^- + i_1^+ + i_2^+ \nonumber\\ 
\partial (E_1) & = & k_1^- + k_2^+ + k_3^+ + i_2^- + \bar E_1 + \bar E_3 \nonumber\\
\partial (E_2) & = & k_1^+ + k_2^- + k_3^+ + \bar E_2 + \bar E_1\nonumber \\
\partial (E_3) & = & k_1^+ + k_2^+ + k_3^- + i_1^- + \bar E_3 + \bar E_2 \nonumber \\
\partial(\bar A_{11}) & = & \bar B_1^1 + \bar C_{45} + \bar B_3^1 + \bar C_{12} + \bar B_5^1 + \bar C_{23} + \bar E_+\nonumber\\ 
\partial(\bar A_{12}) & = & \bar B_1^2 + \bar C_{35} + \bar C_{12} + \bar D_{34}^+ + \bar C_{24} + \bar E_2\nonumber \\
\partial(\bar A_{13}) & = & \bar B_1^3 + \bar B_4^1 + \bar D_{34}^- + \bar D_{14}^-\nonumber\\ 
\partial(\bar A_{21}) & = & \bar C_{45} + \bar B_5^2 + \bar D_{12}^+ + \bar E_1\nonumber \\
\partial(\bar A_{22}) & = & \bar C_{24} + \bar E_2 + \bar C_{14} + \bar D_{24}^- + \bar E_3\nonumber\\ 
\partial(\bar A_{23}) & = & \bar C_{25} + \bar B_4^2 + \bar C_{13} + \bar E_1 + \bar C_{34} + \bar D_{14}^+\nonumber\\ 
\partial(\bar A_{31}) & = & \bar B_2^1 + \bar C_{23} + \bar B_5^3 + \bar D_{12}^- + \bar C_{12} + \bar D_{34}^-\nonumber\\ 
\partial(\bar A_{32}) & = & \bar C_{25} + \bar C_{13} + \bar C_{24} + \bar C_{35} + \bar D_{13}^+\nonumber\\ 
\partial(\bar A_{33}) & = & \bar C_{13} + \bar C_{14} + \bar D_{24}^- + \bar E_3 + \bar E_1\nonumber
\end{eqnarray}
\begin{eqnarray}
\partial(\bar A_{41}) & = & \bar C_{35} + \bar B_2^2 + \bar C_{14} + \bar E_3 + \bar C_{12} + \bar D_{34}^+\nonumber\\ 
\partial(\bar A_{42}) & = & \bar C_{25} + \bar C_{35} + \bar C_{24} + \bar D_{13}^- + \bar E_1 + \bar C_{13} + \bar E_2\nonumber\\ 
\partial(\bar A_{43}) & = & \bar C_{14} + \bar B_3^2 + \bar C_{25} + \bar C_{34} + \bar D_{14}^+ + \bar E_3\nonumber\\ 
\partial(\bar A_{51}) & = & \bar C_{34} + \bar B_2^3 + \bar C_{15} + \bar B_4^3 + \bar C_{23} + \bar E_+ + \bar B_5^1\nonumber\\ 
\partial(\bar A_{52}) & = & \bar C_{15} + \bar D_{23}^+ + \bar E_2 + \bar B_5^2\nonumber\\ 
 \partial(\bar A_{53}) & = & \bar C_{23} + \bar B_3^3 + \bar C_{34} + \bar D_{23}^- + \bar D_{14}^- + \bar B_5^3 \nonumber
\end{eqnarray}
\end{proposition}

\begin{corollary}
The kernel of this operator is generated by the right parts of 15 formulas from the next Proposition \ref{6to5prop}. \hfill $\Box$
\end{corollary}

\subsection{\bf $\partial: C_6 \to C_5$}
\label{6to5}
\begin{proposition}
\label{6to5prop} The boundary operator $C_6 \to C_5$ is given by the following formulas:
\end{proposition}
\begin{eqnarray}
\partial(A_{11}) & = & B_1^1 + C_{45} + B_3^1 + C_{12} + B_5^1 + \bar A_{11} + \bar A_{51} \nonumber \\
\partial(A_{12}) & = & B_1^2 + C_{35} + C_{12} + D_{34}^+ + \bar A_{12} + \bar A_{52} \nonumber\\
\partial(A_{13}) & = & B_1^3 + C_{34} + B_4^1 + D_{34}^- + \bar A_{13} + \bar A_{53}\nonumber\\ 
 \partial(A_{21}) & = & C_{45} + C_{13} + B_5^2 + D_{12}^+ + \bar A_{21} + \bar A_{23}\nonumber\\ 
\partial(A_{22}) & = & C_{24} + E_2 + \bar A_{22} + \bar A_{33}\nonumber\\ 
\partial(A_{23}) & = & C_{25} + B_4^2 + C_{13} + E_1 + \bar A_{23} + \bar A_{43}\nonumber\\ 
\partial(A_{31}) & = & B_2^1 + C_{23} + B_5^3 + D_{12}^- + \bar A_{31} + \bar A_{13}\nonumber\\ 
\partial(A_{32}) & = & C_{25} + C_{13} + C_{24} + C_{35} + C_{14} + D_{13}^+ + D_{24}^+ + E_+\nonumber\\ 
\partial(A_{33}) & = & C_{24} + C_{13} + C_{14} + D_{24}^- + E_3 + \bar A_{33} + \bar A_{42}\nonumber\\ 
\partial(A_{41}) & = & C_{35} + B_2^2 + C_{14} + E_3 + \bar A_{41} +\bar A_{12}\nonumber\\ 
\partial(A_{42}) & = & C_{25} +C_{35} + C_{24} + D_{13}^- + E_1 + \bar A_{42} + \bar A_{22}\nonumber\\ 
\partial(A_{43}) & = & C_{23} + C_{14} + B_3^2 + C_{25} + C_{34} + D_{14}^+ + \bar A_{43} + \bar A_{41}\nonumber\\ 
\partial(A_{51}) & = & C_{34} + B_2^3 + C_{15} + B_4^3 + C_{23} + E_+ + \bar A_{51} + \bar A_{11}\nonumber\\ 
\partial(A_{52}) & = & C_{15} + D_{23}^+ + E_2 + \bar A_{52} + \bar A_{21}\nonumber\\ 
\partial(A_{53}) & = & C_{23} + B_3^3 + C_{34} + D_{23}^- + D_{14}^- + \bar A_{53} + \bar A_{31} \nonumber
\end{eqnarray}

}

\section{Cohomology ring of the complex $\overline{CD}_3(S^1)$. Stiefel--Whitney classes of the bundle $\N_3$}
\label{homapp}

\subsection{}Theorem \ref{thmhom} now has the following concretization.

\begin{proposition}
\label{hpro3}
The group $H_1\left(\overline{CD}_3(S^1), {\mathbb Z}_2\right)$ is isomorphic to ${\mathbb Z}_2$ and is generated by any one of the three cells \ $\bar \Omega_1,$ $\bar \Omega_2$ \ or \ $\bar \Theta$. 

The group $H_2\left(\overline{CD}_3(S^1), {\mathbb Z}_2\right)$ is isomorphic to ${\mathbb Z}_2^2 $ and is generated by the classes of 1$)$ the cell \ $\bar S$ \ and \ 2$)$ each of the cells \ $\bar Y_1$, $\bar Y_2$, \ and \ $\bar J$. 

The group $H_3\left(\overline{CD}_3(S^1), {\mathbb Z}_2\right)$ is isomorphic to ${\mathbb Z}_2^2$ and is generated by the cycles $ \bar h^+ + \bar h^-$ and 
\begin{equation}
\label{2nd3} \bar d + \bar k_1^- + \bar k_2^- + \bar h^+ .
\end{equation}

All groups $H_i\left(\overline{CD}_3(S^1), {\mathbb Z}_2\right)$ with $ i \geq 4$ are trivial. \end{proposition}

\noindent
{\it Proof.} These statements all follow directly from the formulas for the boundary operators given in the Propositions \ref{2to1prop}, \ref{3to2prop}, \ref{4to3prop}, \ref{5to4prop}, and \ref{6to5prop}. Calculating the homology groups of the complex defined by these operators can be considerably simplified using the filtration defined in \S~\ref{filth}. Namely, let $\Psi_0$ (respectively $ \Psi_1 $, $\Psi_2$) be the subspace of $\overline{CD}_3(S^1)$ consisting of all equilevel algebras of codimension three with multiplicity $4$ (respectively, $\leq 5$, $\leq 6$), see Proposition \ref{filthprop}; in particular $\Psi_2 \equiv \overline{CD}_3(S^1)$. 

\unitlength 0.6mm
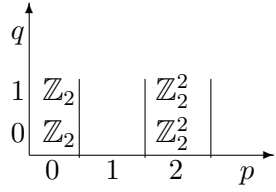
\begin{figure}
\begin{center}
\begin{picture}(100,63)
\put(7,9){\vector(1,0){84}}
\put(7,9){\vector(0,1){55}}
\put(16,0){\small 0}
\put(36,0){\small 1}
\put(56,0){\small 2}
\put(83,0){\small $p$}
\put(0,54){\small $q$}
\put(0,15){\small 0}
\put(0,31){\small 1}
\put(27,9){\line(0,1){33}}
\put(47,9){\line(0,1){33}}
\put(67,9){\line(0,1){33}}
\put(13,14){${\mathbb Z}_2$}
\put(13,30){${\mathbb Z}_2$}
\put(53,14){${\mathbb Z}_2^2$}
\put(53,30){${\mathbb Z}_2^2$}
\end{picture}
\end{center}
\caption{Spectral sequence for $H_*(\overline{CD}_3(S^1), {\mathbb Z}_2)$}
\label{spse}
\end{figure}

\begin{lemma}
\label{prohh}
All nonzero groups $E^1_{p,q}$ of the 
term $E^1$ of the spectral sequence that computes the group \ $H_*(\overline{CD}_3(S^1), {\mathbb Z}_2)$ and is defined by the filtration $\{\Psi_i\}$ are given in Fig.~\ref{spse}. The group $E^1_{0,1} \equiv H_1(\Psi_0, {\mathbb Z}_2)$ is generated by the class of either the cell $\bar \Omega_1$ or the cell $\bar \Omega_2$.
The group $E^1_{2,0} \equiv H_2(\Psi_2, \Psi_1; {\mathbb Z}_2)$ is generated by the relative homology classes of a$)$ either of the cells $\bar S$ or $\bar X$, and b$)$ either of the cells $\bar Y_1$, $\bar Y_2$ or $\bar J$. The group $E^1_{2,1} \equiv H_3(\Psi_2, \Psi_1; {\mathbb Z}_2)$ is generated by the classes of the cells \ $\bar h^+$ \ and \ $\bar h^-$.
\end{lemma}

\noindent
{\it On the proof of Lemma \ref{prohh}.} The groups $H_*(\Psi_0)$, $H_*(\Psi_1, \Psi_0)$, and $H_*(\Psi_2, \Psi_1)$, which form the term $E^1$, can be calculated using the exact sequence defined by the additional two-term filtration, whose first term consists of all algebras of the second type. The complexes of relative cycles of these exact sequences (generated by the cells of the first type) are acyclic in all three cases, and the subcomplexes of cycles of the second type provide the entire homology groups shown in Fig.~\ref{spse} and represented by the cycles given in Lemma \ref{prohh}. \hfill $\Box_{\ref{prohh}}$
\medskip

Thus, the groups $H_i(\overline{CD}_3(S^1), {\mathbb Z}_2)$ are no greater than those indicated in Proposition \ref{hpro3}. The formulas of the previous section show that the chains given in this proposition are absolute cycles in $\overline{CD}_3(S^1)$. Their classes in the groups $E^1_{p,q}$ are independent, so they define independent classes in $H_*(\overline{CD}_3(S^1), {\mathbb Z}_2)$. \hfill $\Box_{\ref{hpro3}}$ $\Box_{\ref{thmhom}}$
\medskip

The first two statements of Theorem \ref{thmsw} now follow from Lemmas \ref{lem86} and \ref{lem87}. Indeed, by Lemma \ref{lem87}, the classes $w_1(\N_3)$ and $w_1^2(\N_3)$ are nontrivial. By Lemmas \ref{lem86} and \ref{lem87}, the class $w_2(\N_3)$ is also nontrivial, but it is different from $w_1^2(\N_3)$.

Below, we give alternative proofs of these two statements of Theorem \ref{thmsw}, which are based on the topology of cycles $S$ and $J$. Namely, these statements also follow from the following Propositions \ref{prosw1} and \ref{propJ}.

\begin{proposition}
\label{prosw1}
1. The first Stiefel--Whitney class of the vector bundle $\N_3$ takes the nonzero value on the cycle defined by the closure of the cell $\bar \Theta$.

2. The square of the first Stiefel--Whitney class of the bundle $\N_3$ takes the nonzero value on the generator $\bar S$ of the group $H_2\left(\overline{CD}_3(S^1), {\mathbb Z}_2\right)$, while the second Stiefel--Whitney class of the same bundle takes the zero value on this generator.
\end{proposition}

\noindent
{\it Proof}. Restricted to the closure $\{\bar S\}$ of the cell $\bar S$, the bundle $\N_3$ is the sum of a two-dimensional constant bundle, whose fiber is the quotient space of $C^{\infty}(S^1, {\mathbb R})$ by the subspace of functions with $f'(\bullet)=f''(\bullet)=0$, and the bundle of one-dimensional normals to the fibers of the tautological bundle 
on the space of two-dimensional subspaces in 
$\M^3_\bullet /\M^6_\bullet \simeq {\mathbb R}^3$, see Proposition \ref{Sstruc}. Thus, the values of the characteristic classes $w_1^2(\N_3)$ and $w_2(\N_3)$ on the cycle $\{\bar S\}$ are equal to the values of the analogous characteristic classes of the tautological one-dimensional bundle on ${\mathbb R}P^2$. The closure of the cell
$\bar \Theta$ defines a nontrivial homology class in this projective plane. Therefore, the class $w_1(\N_3)$ takes the non-zero value on it. 
\hfill {$\Box_{\ref{prosw1}} \ \Box_{\ref{thmsw}(1)}$}

\begin{proposition}
\label{propJ}
The class $w_2(\N_3)$ takes the nonzero value on the class $\{\bar J\subset H_2\left(\overline{CD}_3(S^1), {\mathbb Z}_2\right)$ defined by the closure of the cell $\bar J$, while the class $w_1^2(\N_3)$ takes the zero value on it.
\end{proposition}

\noindent
{\it Proof.} The cycle $\{\bar J\}$ can be identified with a two-dimensional subvariety of the Grassmann manifold of two-dimensional subspaces of the space $\M_\bullet^2/\M_\bullet^6 \simeq {\mathbb R}^4$, see Proposition \ref{Jstruc}: each algebra ${\mathcal F}$ in this class is associated with the subspace $({\mathcal F} \cap \M_\bullet)/\M^6_{\bullet} \subset \M_\bullet^2/\M_\bullet^6$. The restriction of the bundle $\N_3$ to this variety is equal to the sum of the one-dimensional constant bundle $\M_\bullet/\M_\bullet^2$ and the two-dimensional canonical {\em normal} bundle of this Grassmann manifold.

The group $H_2(G_2({\mathbb R}^4), {\mathbb Z}_2) \simeq {\mathbb Z}_2^2$ is generated by the classes of two Schubert cells $X(1,4)$ and $X(2,3)$, see e.g. \cite{MS}.

\begin{lemma}
\label{le51}
The intersection indices of the cycle $\{\bar J\} \subset G_2(\M_\bullet^2/\M_\bullet^6)$ with the classes of the Schubert cells $X(1,4)$ and $X(2,3)$ are nontrivial.
\end{lemma}

\noindent
{\it Proof.} The homology class of the Schubert cell $X(1,4)$ in the Grassmann manifold $G_2({\mathbb R}^4)$ can be realized by the set of all two-dimensional subspaces that contain a fixed one-dimensional subspace of ${\mathbb R}^4$. Take for this subspace in $\M^2_\bullet/\M^6_\bullet$ the line spanned by a polynomial of the form $x^2 + a x^3 + b x^4 + c x^5,$ where $x$ is a local coordinate in $S^1$ with the origin at the point $\bullet$. Each polynomial $f$ of this form belongs to a single equilevel subalgebra of the class $J$, consisting of the linear combinations of the functions $1$, $f$, $f^2$, 
and the ideal $\M_\bullet^6$. 

The homology class of the Schubert cell $X(2,3)$ can be realized by the variety of all two-dimensional subspaces contained in a fixed three-dimensional subspace in ${\mathbb R}^4$. For this three-dimensional subspace, take the space of polynomials $a x^2 + b x^3 + c x^4 + d x^5$ whose coefficients satisfy the condition $ a + \alpha b + \beta c + \gamma d = 0$ for a generic triple of coefficients $(\alpha, \beta, \gamma)$. Each such subspace with $\gamma \neq 0$ contains a unique subspace from the cycle $\{\bar J\}$ spanned by the polynomials $$f_0 = x^2 + \frac{\beta}{2\gamma} x^3 + \frac{\alpha \beta - 2 \gamma}{2 \gamma^2} x^5$$ and $f_0^2 (\mbox{mod }\M_\bullet^6)$. \hfill $\Box_{\ref{le51}}$ 

\begin{corollary}
\label{cor84}
The subvariety \ $ \{\bar J\} \subset G_2(\M_\bullet^2/\M_\bullet^6)$ is homologous to the sum of the Schubert cells $X(1,4)$ and $X(2,3)$.
\end{corollary}

\noindent
{\it Proof.} The mod 2 intersection indices of the two-di\-men\-si\-onal Schubert cells in $G_2({\mathbb R}^4)$ are as follows: $\langle X(1,4), X(1,4) \rangle =$ $\langle X(2,3), X(2,3) \rangle =1$, $\langle X(1,4), X(2,3) \rangle =0$. Thus, the linear functions on the group $H_2(G_2({\mathbb R}^4), {\mathbb Z}_2)$ defined by the intersection indices with cycles $\{\bar J\}$ and $X(1,4) +X(2,3)$ are the same. By Poincar\'e duality, this implies that these cycles are homological. \hfill $\Box_{\ref{cor84}}$ 
\medskip

The pairings of the Schubert cells with Stiefel--Whitney classes of the tautological bundle over $G_2({\mathbb R}^4)$ are as follows: $w_1^2(X(1,4))=1, $ $w_2(X(1,4))=0$, $w_1^2(X(2,3))=1, $ and $w_2(X(2,3))=1$. The Stiefel--Whitney classes $\bar w_1$ and $\bar w_2$ of the canonical normal bundle of the Grassmann manifold $G_2({\mathbb R}^4)$ are equal to $w_1$ and $w_2+w_1^2$, respectively. Therefore, $\bar w_2 (X(1,4)+X(2,3)) = 1$, $\bar w_1^2 (X(1,4)+X(2,3))=0$. 
\hfill $\Box_{\ref{propJ}}$ $\Box_{\ref{thmsw}(2)}$ 

\begin{proposition} 
\label{pro18}
The cycle
$\bar h^+ + \bar h^-$
 is homologous to the cycle $\aleph_3$, see \S~\ref{gener}. 
\end{proposition}

\noindent
{\it Proof.} The closure of the set of all algebras $h(\alpha, \beta; A, B)$ is fibered over the space $\overline{B(S^1,2)}$ of pairs of (possibly coinciding) points $A, B \in S^1$ (with exceptional fibers over boundary points of the form $(A, A)).$ The closures of the subvarieties $\aleph_3$ and $\bar h^+ + \bar h^-$ of codimension 1 in this closure are the preimages of two loops in $\overline{B(S^1,2)}$: one consisting of pairs of opposite points and the other consisting of pairs containing the point $\bullet$. These loops are homologous to each other, so the homology between their preimages is defined by the preimage of the area in $\overline{B(S^1,2)}$ realizing the homology between these loops. \hfill $\Box_{\ref{pro18}}$ 

\begin{proposition}
\label{protri2}
The cosets of functions 
\begin{equation}
\label{fou1}
 a \sin \varphi + b \cos \varphi
\end{equation}
with arbitrary $a$ and $b$ 
form a trivial two-dimensional subbundle of the restriction of the vector bundle $\N_3$ to the closure of 
the cycle $($\ref{2nd3}$)$. 
\end{proposition}

\noindent
{\it Proof.} No non-zero function of the form (\ref{fou1}) belongs to any subalgebra from this cycle. Indeed, every function from an algebra of the class $h^+$ has a pair of points with the same value and {\em the same sign of the derivative} at these points.
All functions from all algebras which belong to cells $d$ and $k_i^{\pm}$ or to their boundaries or to the boundary of the cell $h^+$ have values of multiplicity $\geq 3$.
 \hfill $\Box_{\ref{protri2}}$

\begin{corollary}
\label{coro36}
The Stiefel--Whitney classes $w_3(\N_3)$ and $w_1(\N_3) \smile w_2(\N_3)$ take the zero value on the homology class of the cycle $($\ref{2nd3}$)$. \hfill $\Box_{\ref{coro36}}$
\end{corollary}

In the notation of \S~\ref{cells}, the cycle $\diamondsuit_3$ (see \S~\ref{dvd}) is the closure of the union of all algebras $E(A, A+2\pi/3, A+ 4\pi/3)$ over all $A \in [0,2\pi/3]$.

\begin{proposition}
\label{probet}
The cycle $\diamondsuit_3$ is homologous to the cycle $($\ref{2nd3}$)$ in $\overline{CD}_3(S^1)$.
\end{proposition}

This proposition will be proved in \S~\ref{prsw0}. 

\subsection{Topology of the cycle $\aleph_3$.} 

The cycle $\aleph_3$ is the space of a fiber bundle 
\begin{equation}
\label{fb2}
\aleph_3 \to S^1/{\mathbb Z}_2
\end{equation}
over the space of diameters of $S^1$. 
 
According to Proposition \ref{proph}, its fiber over any point $\{A, A+\pi\}$ is homeomorphic to the projective plane and is smooth everywhere except at the point of class $\Lambda$. In particular, $\aleph_3$ is a topological manifold.

The natural {\it equatorial circle} $\circ(A)\equiv \circ(A+\pi)$ of this fiber is the closure of the set of all algebras $h(\alpha, \beta; A, A+\pi)$
defined by the conditions (\ref{eqh}) with $\beta=0$ and arbitrary $\alpha \in {\mathbb R} P^1$. This fiber also has two distinguished {\it meridians} $M_1(A)$ and $M_{-1}(A)$  which are intervals with the ends at the singular point: they consist of all algebras $h(1, \beta; A, A+\pi)$ and $h(-1, \beta; A, A+\pi)$
 with arbitrary $\beta \in {\mathbb R}$. 

The fiber bundle (\ref{fb2}) has a canonical local trivialization: the corresponding flat connection maps any algebra $h(\alpha, \beta; A, A+\pi)$ to $h(\alpha, \beta; \tilde A, \tilde A+\pi)$ for any two neighboring points $A$ and $\tilde A$.

\begin{proposition}
\label{monodr}
 The transportation over the entire loop $S^1/{\mathbb Z}_2$ by this local trivialization

a$)$ moves the equator $\circ(A)$ to itself, changing its orientation;

b$)$ moves the distinguished meridian $M_{1}(A)$ to itself, changing its orientation; 

c$)$ moves the meridian $M_{-1}(A)$ to itself identically.
\end{proposition}

\noindent
{\it Proof.} a) This transportation moves any algebra $h(\alpha, 0; 0, \pi)$ to $h(\alpha, 0; \pi, 0) \equiv $ $ h(\alpha^{-1}, 0; 0, \pi)$. 

b) It also moves any algebra $h(1, \beta; 0,\pi)$ to the algebra $h(1, \beta; \pi, 0),$ which, by identity (\ref{idenh}), equals $h(1, - \beta; 0, \pi )$. 

c) Similarly, it moves any algebra $h(-1, \beta; 0,\pi)$ to the algebra 
$h(-1, \beta; \pi, 0) \equiv h(-1, \beta; 0, \pi).$ \hfill $\Box_{\ref{monodr}}$

\begin{corollary}
\label{cormonodr}
There is a group isomorphism $$H_*(\aleph_3, {\mathbb Z}_2) \simeq H_*(S^1, {\mathbb Z}_2) \otimes H_*({\mathbb R} P^2, {\mathbb Z}_2).$$
\end{corollary}

\noindent
{\it Proof.} By Proposition \ref{monodr}, all differentials of the mod 2 Wang's exact sequence related with the fiber bundle (\ref{fb2}) are trivial.
\hfill $\Box_{\ref{cormonodr}}$
\medskip

There are four important 1-cycles in $\aleph_3$: 

1) the union of the singular points of all fibers, 

2) an arbitrary non-contractible cycle in the fiber over the point $\{0,\pi\} \in S^1/{\mathbb Z}_2$,

3) the union of the points $h(1, 0; A, A+\pi)$ over all $A$; 

4) the union of the points $h(-1, 0; A, A+\pi)$ over all $A$. 

\begin{lemma}
\label{le38}
All these cycles define nontrivial elements of the group $H_1(\aleph_3, {\mathbb Z}_2) \simeq {\mathbb Z}_2^2$. The transportations along the cycles 1$)$ and 4$)$ preserve the orientation of the  bundle $\N_3$, while the transportation along the cycles 2$)$ and 3$)$ change it.
\end{lemma}

\noindent
{\it Proof.} The projections of cycles 1), 3), and 4) to the base of the fiber bundle (\ref{fb2}) are maps of degree 1. The entire group $H_1(\aleph_3, {\mathbb Z}_2) \simeq {\mathbb Z}_2^2$ is generated by either of these cycles and the nontrivial homology class of the fiber, so the cycle 2) representing the latter class is also not zero-homologous in $\aleph_3$.

The points of cycle 1) are the algebras defined by the three linear equations $f(A) - f(A+\pi)=0$, $f'(A)=0$, $f'(A+\pi)=0$. Moving $A$ along the interval $(0, \pi)$ multiplies the gradient of the left part of the first equation by $-1$ and permutes the gradients of the other two conditions. Thus, the frame of the normal bundle formed by these gradients preserves its orientation.

The algebras of cycle 3) are distinguished by the three linear equations $f(A) - f(A+\pi)=0$, $f'(A) - f'(A+\pi)=0$, $f''(A) - f''(A+\pi)=0$. The gradients of all their left parts are multiplied by $-1$ by the transportation over the cycle $S^1/{\mathbb Z}_2$, so this cycle changes the orientation of $\N_3$. 

The case of cycle 4) is very similar, but the second equation has in this case the form
$f'(A) + f'(A+\pi)=0$, so its gradient is preserved by the transportation.

Cycles 1) and 2) generate the group $H_1(\aleph_3, {\mathbb Z}_2)$, so at least one of them must change the orientation of the bundle $\N_3$. We know that this is not cycle 1). \hfill $\Box_{\ref{le38}}$

\begin{corollary}
\label{cor91}
Cycles 1$)$ and 4$)$ are homologous to each other. The sum of cycles 1$)$, 2$),$ and 3$)$ is homologous to zero.
\end{corollary}

\noindent
{\it Proof.} There is only one nontrivial element in the group $H_1(\aleph_3, {\mathbb Z}_2) \simeq {\mathbb Z}_2^2$ on which the nontrivial linear form $w_1(\N_3)$ takes the value zero. By Lemma \ref{le38}, both cycles 1) and 4) satisfy this condition. 

Cycles 2) and 3) are not homologous to each other, since they have different degrees of projection onto the base of the fiber bundle. Thus, their difference is a nontrivial cycle that does not change the orientation of $\N_3$, so it is homologous to the cycle 1). \hfill $\Box_{\ref{cor91}}$
\medskip

By Lemma \ref{le79} and Corollary \ref{cormonodr}, the group $H^3(\aleph_3, {\mathbb Z}_2)$ is generated by the third Stiefel--Whitney class of the normal bundle $\N_3$.

\begin{proposition}
\label{prohh2}
1$)$ The class $w_1^3(\N_3)$ takes the nonzero value on the generator of the group $H_3(\aleph_3, {\mathbb Z}_2)$, 

2$)$ the class $w_1(\N_3) \smile w_2(\N_3)$ takes zero value on this generator,

3$)$ the classes $w_1^2(\N_3)$ and $w_2(\N_3)$ generate the group $H^2(\aleph_3, {\mathbb Z}_2)$,

4$)$ both of these classes $w_1^2(\N_3)$ and $w_2(\N_3)$ take the nonzero value on the fibers of the fiber bundle $($\ref{fb2}$)$.
\end{proposition}

\noindent
{\it Proof.} 
Define the surface ${\mathcal E} \subset \aleph_3$ as the union of the equatorial circles $\circ(A)$ of all fibers over the points $\{A, A+\pi\} \in S^1/{\mathbb Z}_2$.

\begin{lemma}
\label{cor72}
The first Stiefel--Whitney class of the vector bundle $\N_3$ is equal to the intersection index with the submanifold ${\mathcal E}$ in the topological manifold $\aleph_3$. 
\end{lemma}

\noindent
{\it Proof.} By Lemma \ref{le38}, the orientation of this vector bundle changes over a one-dimensional cycle in $\aleph_3$ if and only if this cycle has an odd intersection number with ${\mathcal E}$. \hfill $\Box_{\ref{cor72}}$
\medskip

 The cohomology class of the square of this Stiefel-Whitney class on $ \aleph_3$ is therefore Poincar\'e dual to the intersection set of this surface ${\mathcal E}$ with any other cycle that is homologous to it and intersects it transversally. Let us construct such a cycle. 

Fix a small number $\varepsilon > 0$. Consider the points $h(1, \varepsilon; 0, \pi)$ and $h(1,-\varepsilon; 0, \pi)$ and connect them in the fiber of the bundle (\ref{fb2}) over the point $\{0, \pi\} \in S^1/{\mathbb Z}_2$ by the curve that goes at the distance $\varepsilon$ from the equatorial circle $\circ(0)$ from the first of these points in the direction of the increase of $\alpha$. Then, connect the endpoints of this curve by the segment consisting of algebras $h(1, \beta; 0, \pi)$ with $\beta \in [-\varepsilon, \varepsilon]$. The obtained closed curve intersects the equator of our fiber at the single point $h(1, 0; 0, \pi)$. Using the standard local trivialization of the bundle (\ref{fb2}), transport this curve to the fibers over all other points $\{A, A+\pi\}$, $A \in [0, \pi)$. By Proposition \ref{monodr}, these curves sweep out a closed surface in $\aleph_3$. This surface is obviously homologous to the surface ${\mathcal E}$ and intersects it transversally along the curve 3) from Lemma \ref{le38}.
So the class $w_1^2(\N_3)$ in $\aleph_3$ can be defined by the intersection indices with this curve. 

By Lemma \ref{le38}, this curve changes the orientation of the bundle $\N_3$, so $w_1^3(\N_3)$ takes the non-zero value on the cycle $\aleph_3$. \hfill $\Box_{\ref{prohh2}(1)}$
\medskip

To realize the class $w_2$ of the normal bundle, consider the two-dimensional subspace in $C^\infty(S^1, {\mathbb R})$ consisting of the functions (\ref{fou1}). For every $A \in [0,\pi]$ there is exactly one algebra $h(\alpha, \beta; A, A+\pi)$ having a nontrivial intersection with this subspace: it is the algebra $h(-1, 0; A, A+\pi)$. The union of these algebras over all $A \in [0,\pi]$ is the curve 4) from Lemma \ref{le38}, so the restriction of the class $w_2(\N_3)$ to the manifold $\aleph_3$ is Poincar\'e dual to the class of this curve. By the same lemma, the class $w_1(\N_3)$ takes the zero value on this curve, so $w_1(\N_3) \smile w_2(\N_3) =0$ in $H^3(\aleph_3, {\mathbb Z}_2)$. \hfill $\Box_{\ref{prohh2}(2)}$
\medskip

By Lemma \ref{le38}, the 1-cycles 3) and 4) from this lemma generate the group $H_1(\aleph_3, {\mathbb Z}_2)$, so their Poincar\'e dual classes $w_1^2(\N_3)$ and $w_2(\N_3)$ generate the group $H^2(\aleph_3, {\mathbb Z}_2)$. \hfill $\Box_{\ref{prohh2}(3)}$
\medskip

Each of these two cycles has a single transversal intersection with each fiber of the bundle (\ref{fb2}). \hfill $\Box_{\ref{prohh2}(4)}$

\subsection{The topology of the cycle $\diamondsuit_3$.}

Consider the fiber bundle
\begin{equation}
\label{fbt}
\diamondsuit_3 \to S^1/{\mathbb Z}_3 ,
\end{equation}
which is a special case of the fiber bundle (\ref{fb9}) for $k=3$. Its fiber
over an arbitrary point $\{A, A+ 2\pi/3, A+ 4\pi/3\} \in S^1/{\mathbb Z}_3$ is diffeomorphic to ${\mathbb R} P^2$. This fiber
consists of all algebras defined by the conditions 
\begin{equation}
\label{eqdiam}
\begin{split}
& f(A)=f(A+2\pi/3)=f(A+4\pi/3), \\ & \alpha f'(A)+\beta f'(A+2\pi/3)+ \gamma f'(A+4\pi/3)=0
\end{split}
\end{equation} 
for an arbitrary point $(\alpha:\beta:\gamma)\in {\mathbb R}P^2$. The mod 2 Wang exact sequence of this fiber bundle yields the group isomorphism 
\begin{equation}
\label{kun}
H^*(\diamondsuit_3, {\mathbb Z}_2) \simeq H^*(S^1, {\mathbb Z}_2) \otimes H^*({\mathbb R} P^2, {\mathbb Z}_2) .
\end{equation}

By Lemma \ref{lem87}, the class $w_1^2(\N_3)$ takes the non-zero value on any fiber of the fiber bundle $($\ref{fbt}$)$. By Lemma \ref{prt23}, the Stiefel--Whitney classes $w_2$ and $w_3$ of the restriction of the bundle $\N_3$ to the manifold $\diamondsuit_3$ are trivial. 

\begin{proposition}
\label{pro22}
The class $w_1^3(\N_3)$ takes zero value on the cycle $\diamondsuit_3$.
\end{proposition}

\noindent
{\it Proof.}
The group $H_1(\diamondsuit_3, {\mathbb Z}_2) \simeq {\mathbb Z}_2^2$ is generated by

 1) an arbitrary 1-cycle in the fiber of the fiber bundle (\ref{fbt}), which generates the 1-homology group of this fiber, and 

2) the cross-section of this fiber bundle that consists of all algebras defined by the conditions (\ref{eqdiam}) with $\alpha=\beta=\gamma.$ 

The first cycle changes the orientation of $\N_3$ while the second cycle preserves it.

For any $A \in [0, 2\pi/3]$, take the curve in the corresponding fiber of the bundle (\ref{fbt}), which consists of all algebras defined by the equations (\ref{eqdiam}) with \begin{equation}
\label{qu}
\alpha^3+\beta^3+\gamma^3=0.
\end{equation}
 The union of these curves over all fibers is a closed surface in $\diamondsuit_3$. The first basic cycle of $H_1(\diamondsuit_3, {\mathbb Z}_2)$ intersects it once, and the second cycle does not intersect it. Therefore, the first Stiefel-Whitney class of the restriction of $\N_3$ to $\diamondsuit_3$ is equal to the intersection index with this surface. Thus, the class $w_1^2(\N_3)$ in $\diamondsuit_3$ is Poincar\'e dual to the self-intersection homology class of this surface. For another realization of the homology class of this surface, we can take the union of the curves in the fibers given by the equation $\alpha^3+\beta^3+\gamma^3 + \varepsilon \cdot \alpha \beta \gamma =0,$ 
where $\varepsilon$ is a small positive number. The intersection of these two surfaces is presented in each fiber by a triple of points that is invariant under all permutations of the coordinates $\alpha, \beta, $ and $\gamma$. Consider a single fiber and shift these three points in it slightly off the curve (\ref{qu}) in a way that leaves the resulting triple invariant under the action on ${\mathbb R} P^2$ of cyclic permutations of the coordinates $\alpha \to \beta \to \gamma \to \alpha$. Repeating this shift in all the other fibers yields a family of triples of points, which sweeps out a closed curve in $\diamondsuit_3$ that also defines the cocycle $w_1^2(\N_3)$. This curve does not intersect our surface formed by curves (\ref{qu}), so $w_1^3(\N_3)$ is Poincar\'e dual to the empty cycle. \hfill $\Box_{\ref{pro22}}$ $\Box_{\ref{th41}}$

\subsection{Proof of Proposition \ref{probet}.}
\label{prsw0}
Consider two two-dimensional cycles in $\overline{CD}_3(S^1)$: 

a) the fiber of the fiber bundle (\ref{fb2}) over the point 
$\{\pi/2, 3\pi/2\} \in S^1/{\mathbb Z}_2$, 

b) the fiber of the fiber bundle
(\ref{fbt}) over the point $\{0, 2\pi/3, 4\pi/3\} \in S^1/{\mathbb Z}_3$.

\begin{lemma}
\label{lemma72}
These two cycles define independent elements of the group
$H_2(\overline{CD}_3(S^1), {\mathbb Z}_2)$.
\end{lemma}

\noindent
{\it Proof.} By Proposition \ref{prohh2}(4) and Lemmas \ref{lem87}(2) and \ref{prt23}, the cohomology classes $w_1^2(\N_3)$ and $w_2(\N_3)$ take the values 1 and 1 on the first cycle and the values 1 and 0 on the second cycle. \hfill $\Box_{\ref{lemma72}}$
\medskip

By Proposition \ref{prohh}, $H_2(\Psi_1, {\mathbb Z}_2)=0$. Therefore, by the exact sequence of the pair $(\Psi_2, \Psi_1)$, the two cycles considered in Lemma \ref{lemma72} define 
independent elements of the group 
$$H_2(\overline{CD}_3(S^1), \Psi_1; {\mathbb Z}_2) \equiv H^{l f}_2(\Psi_2\setminus \Psi_1, {\mathbb Z}_2) \simeq {\mathbb Z}_2^2 ,$$ where $H^{l f}$ means the Borel--Moore homology. 

Denote by $\varkappa$ the fiber over the point $0 \in S^1$ of the fiber bundle 
 \begin{equation}
\label{fibrez}
\Psi_2\setminus \Psi_1 \to S^1
\end{equation}
 introduced in \S \ref{filth}. The intersections of our two cycles with the set $\Psi_2 \setminus \Psi_1$ lie in $\varkappa$, so they define 
two linearly independent elements of the group $H^{l f}_2(\varkappa; {\mathbb Z}_2)$.

Transportations of these cycles by the covering homotopy over the base $S^1$ of the bundle (\ref{fibrez}) can be realized by the intersections of the fibers of this bundle with the Borel--Moore cycles $\aleph_3 \cap (\Psi_2 \setminus \Psi_1)$ and $\diamondsuit_3 \cap (\Psi_2 \setminus \Psi_1)$. Therefore, by Wang's exact sequence of this fiber bundle, the cycles $\aleph_3$ and $\diamondsuit_3$
define independent elements of the group $H^{l f}_3(\Psi_3 \setminus \Psi_2, {\mathbb Z}_2) \equiv H_3(\Psi_3, \Psi_2; {\mathbb Z}_2)$, and thus, by Proposition \ref{prohh}, also of $H_3\left(\overline{CD}_3(S^1), {\mathbb Z}_2\right)$. 

By Proposition \ref{pro22}, the cocycle $w_3(\N_3)$ takes zero value on $\diamondsuit_3$. By Lemma \ref{le79} and Corollary \ref{coro36}, the cycle (\ref{2nd3}) defines the only nontrivial element of the group $H_3\left(\overline{CD}_3(S^1), {\mathbb Z}_2\right)$ with this property. Therefore,  this cycle (\ref{2nd3}) is homologous to the cycle $\diamondsuit_3$. \hfill $\Box_{\ref{probet}}$

\subsection{Proof of last two statements of Theorem \ref{thmsw}.}
\label{prsw}

By Lemma \ref{le79} and Proposition \ref{prohh2}(1), the classes $w_3(\N_3)$ and $w_1^3(\N_3) $ take equal values on the cycle $\aleph_3$. By 
Corollary \ref{coro36} and Proposition \ref{pro22}, they also take equal values on the cycle $\diamondsuit_3$. \hfill $\Box_{\ref{thmsw}(4)}$
\medskip

By Proposition \ref{prohh2}(2) and Corollary \ref{coro36},
the class $w_1(\N_3) \smile w_2(\N_3)$ takes the zero value on each of these two homology classes. \hfill $\Box_{\ref{thmsw}(5)}$

}
\newpage

{
\section*{Appendix: matrices of boundary operators}
\vspace{0.5cm}
 
\begin{center}
\begin{tabular}{| l || c | c | c | c|c|c|c|}
\hline
$\vphantom{\overline{(}}$ & $\bar \Upsilon_1 \vphantom{\overline{(}}$ & $\bar \Upsilon_2$ & $\bar \Omega_1$ & $\bar \Omega_2$ & $\bar \Lambda$ & $\bar \Theta$ & $\nabla$ \\
\hline 
\hline
$\Upsilon_1 $ & + & + & & & & & + \\
\hline
$\Upsilon_2$ & + & + & & & & & + \\
\hline
$\Omega_1$ & & & + & + & & & + \\
\hline
$\Omega_2$ & & & + &+ & & & + \\
\hline
$\Lambda$ & & & & & & & + \\
\hline
$\Theta$ & & & & & & & \\
\hline
$\bar U_{12} \vphantom{\overline{(}}$ & + & & & + & + & + & \\
\hline
$\bar U_{13} \vphantom{\overline{(}}$ & & + & + & & + & + & \\
\hline
$\bar U_{21} \vphantom{\overline{(}}$ & + & & & + & + & + & \\
\hline
$\bar U_{23} \vphantom{\overline{(}}$ & & + & + & & + & + & \\
\hline
$\bar U_{31} \vphantom{\overline{(}}$ & & + & & + & + & + & \\
\hline
$\bar U_{32} \vphantom{\overline{(}}$ & & + & & + & + & + & \\
\hline
$\bar V_1 \vphantom{\overline{(}}$ & + & & + & + & & & \\
\hline
$\bar V_2 \vphantom{\overline{(}}$ & + & & + & + & & & \\
\hline
$\bar V_3 \vphantom{\overline{(}}$ & & + & & & & & \\
\hline
$\bar W_1 \vphantom{\overline{(}}$ & & & + & + & + & & \\
\hline
$\bar W_2 \vphantom{\overline{(}}$ & & & + & + &+ & & \\
\hline
$\bar W_3 \vphantom{\overline{(}}$ & & & & & + & & \\
\hline
$\bar X \vphantom{\overline{(}}$ & + & & & & & & \\
\hline
$\bar Y_1 \vphantom{\overline{(}}$ & & & & & & & \\
\hline
$\bar Y_2 \vphantom{\overline{(}}$ & & & & & & & \\
\hline
$\bar Z_1^+ \vphantom{\overline{(}}$ & & & + & & + & + & \\
\hline
$\bar Z_2^+ \vphantom{\overline{(}}$ & & & & + & + & + & \\
\hline
$\bar Z_1^- \vphantom{\overline{(}}$ & & & + & & + & + & \\
\hline
$\bar Z_2^- \vphantom{\overline{(}}$ & & & & + & + & + & \\
\hline
$\bar T_1 \vphantom{\overline{(}}$ & & & & & & & \\
\hline
$\bar T_2 \vphantom{\overline{(}}$ & & & & & & & \\
\hline
$\bar J \vphantom{\overline{(}}$ & & & & & & & \\
\hline
$\bar S \vphantom{\overline{(}}$ & & & & & & & \\
\hline
\end{tabular}
\label{page2to1}

{\bf Boundary operator $\partial_2: C_2 \to C_1$.}
\end{center}
}

{
\scriptsize

\noindent
\begin{tabular}{| l || c | c | c | c|c|c|c|c|c|c|c|c|c|c|c|c|c|c|c|}
\hline
 & $\Theta \vphantom{\overline{(}}$ & $\Upsilon_1$ & $\Upsilon_2$ & $\Omega_1$ & $\Omega_2$ & $\Lambda$ & $\bar U_{12}$ & $\bar U_{13}$ & $\bar U_{21}$ & $\bar U_{23}$ & $\bar U_{31}$ & $\bar U_{32}$ & $\bar V_1$ & $\bar V_2$ & $\bar V_3$ & $\bar W_1$ & $\bar W_2$ & $\bar W_3$ & $ \bar X$ \\
\hline
$U_{12}$ & & + & & & & + & + & & & & + & & & & & & & & \\
\hline
$U_{13}$ & + & & & + & & + & & + & & & & + & & & & & & & \\
\hline
$U_{21}$ & & + & & & + & & & + & + & & & & & & & & & & \\
\hline
$U_{23}$ & & & + & + & & & + & & & + & & & & & & & & &\\
\hline
$U_{31}$ & + & & & & + & + & & & & + & + & & & & & & & &\\
\hline
$U_{32}$ & & & + & & & + & & & + & & & + & & & & & & &\\
\hline
$V_1$ & & + & & + & & & & & & & & & + & & + & & & &\\
\hline
\end{tabular}
\label{page3to2}\\
{\bf Boundary map $\partial_3: C_3 \to C_2$, left top}
\newpage

\begin{tabular}{| l || c | c | c | c|c|c|c|c|c|c|}
\hline
& $\bar Y_1 \vphantom{\overline{(}}$ & $\bar Y_2$ & $\bar Z_1^+$ & $\bar Z_2^+$ & $\bar Z_1^-$ & $\bar Z_2^-$ & $\bar T_1$ & $\bar T_2$ & $\bar J$ & $\bar S$ \\
\hline 
$U_{12}$ & & & & & & & & & & \\
\hline
$U_{13}$ & & & & & & & & & & \\
\hline
$U_{21}$ & & & & & & & & & & \\
\hline
$U_{23}$ & & & & & & & & & & \\
\hline
$U_{31}$ & & & & & & & & & & \\
\hline
$U_{32}$ & & & & & & & & & & \\
\hline
$V_1$ & & & & & & & & & & \\
\hline
\end{tabular}
\nopagebreak

{\bf Boundary map $\partial_3: C_3 \to C_2$, right top}
\newpage

{\scriptsize
\begin{tabular}{| l || c | c | c | c|c|c|c|c|c|c|c|c|c|c|c|c|c|c|c|c|}
\hline
 & $\Theta \vphantom{\overline{(}}$ & $\Upsilon_1$ & $\Upsilon_2$ & $\Omega_1$ & $\Omega_2$ & $\Lambda$ & $\bar U_{12}$ & $\bar U_{13}$ & $\bar U_{21}$ & $\bar U_{23}$ & $\bar U_{31}$ & $\bar U_{32}$ & $\bar V_1$ & $\bar V_2$ & $\bar V_3$ & $\bar W_1$ & $\bar W_2$ & $\bar W_3$ & $ \bar X$ \\
\hline
$\bar e_1^+$ & & & & & & & & & & + & & + & & & & & & &\\
\hline
$\bar e_2^+$ & & & & & & & & + & & & + & & & & & & & &\\ 
\hline
$\bar e_3^+$ & & & & & & & + & & + & & & & & & & & & & \\
\hline
$\bar e_1^-$ & & & & & & & & & & + & & + & & & & & & &\\
\hline
$\bar e_2^-$ & & & & & & & & + & & & + & & & & & & & &\\
\hline
$\bar e_3^-$ & & & & & & & + & & + & & & & & & & & & &\\
\hline
$\bar k_1^+ \vphantom{\overline{(}}$ & & & & & & & & & & & & & & & & & + & + &\\
\hline
$\bar k _2^+ \vphantom{\overline{(}}$ & & & & & & & & & & & & & & & & + & & + & \\ 
\hline
$\bar k_3^+ \vphantom{\overline{(}}$ & & & & & & & & & & & & & & & & + & + & & \\
\hline
$\bar k_1^- \vphantom{\overline{(}}$ & & & & & & & & & & & & & & & & & + & + & \\
\hline
$\bar k_2^- \vphantom{\overline{(}}$ & & & & & & & & & & & & & & & & + & & + &\\
\hline
$\bar k_3^- \vphantom{\overline{(}}$ & & & & & & & & & & & & & & & & + & + & &\\
\hline
$\bar g_1$ & & & & & & & & & & & & & & & & & & &\\
\hline
$\bar g_2$ & & & & & & & & & & & & & & & & & & &\\
\hline
$\bar g_3$ & & & & & & & & & & & & & & & & & & & \\
\hline
$\bar h^+ \vphantom{\overline{(}}$ & & & & & & & & & & & & & & & & & & &\\
\hline
$\bar h^- \vphantom{\overline{(}}$ & & & & & & & & & & & & & & & & & & &\\
\hline
$\bar i_1^+ \vphantom{\overline{(}}$ & & & & & & & & & & & & & & & & & & &\\
\hline
$\bar i_2^+ \vphantom{\overline{(}}$ & & & & & & & & & & & & & & & & & & &\\ 
\hline
$\bar i_1^- \vphantom{\overline{(}}$ & & & & & & & & & & & & & & & & & & &\\
\hline
$\bar i_2^- \vphantom{\overline{(}}$ & & & & & & & & & & & & & & & & & & &\\ 
\hline
\end{tabular}

{\bf Boundary map $\partial_3: C_3 \to C_2$, left bottom}
}
\medskip

{\scriptsize 
\begin{tabular}{| l || c | c | c | c|c|c|c|c|c|c|c|}
\hline
& $\bar Y_1 \vphantom{\overline{(}}$ & $\bar Y_2$ & $\bar Z_1^+$ & $\bar Z_2^+$ & $\bar Z_1^-$ & $\bar Z_2^-$ & $\bar T_1$ & $\bar T_2$ & $\bar J$ & $\bar S$ \\
\hline
$\bar e_1^+$ & & & + & & & + & & & & \\
\hline
$\bar e _2^+$ & & & & + & + & & & & & \\
\hline
$\bar e_3^+$ & & & & + & & + & & & & \\
\hline
$\bar e_1^-$ & & + & & + & + & & & & + & \\
\hline
$\bar e_2^-$ & & + & + & & & + & & & + & \\
\hline
$\bar e_3^-$ & + & & & + & & + & & & + & \\
\hline
$\bar k_1^+ \vphantom{\overline{(}}$ & & & + & & & + & & & & \\
\hline
$\bar k _2^+ \vphantom{\overline{(}}$ & & & & + & + & & & & & \\ 
\hline
$\bar k_3^+ \vphantom{\overline{(}}$ & & & & + & & + & & & & \\
\hline
$\bar k_1^- \vphantom{\overline{(}}$ & & & & + & + & & & + & & \\
\hline
$\bar k_2^- \vphantom{\overline{(}}$ & & & + & & & + & & + & & \\
\hline
$\bar k_3^- \vphantom{\overline{(}}$ & & & & + & & + & + & & & \\
\hline
$\bar g_1$ & + & & & & & & + & + & + & \\
\hline
$\bar g_2$ & + & & & & & & + & + & + & \\
\hline
$\bar g_3$ & & + & & & & & & & + & \\
\hline
$\bar h^+ \vphantom{\overline{(}}$ & & & + & + & + & + & & & & \\
\hline
$\bar h^- \vphantom{\overline{(}}$ & & & + & + & + & + & & & & \\
\hline
$\bar i_1^+ \vphantom{\overline{(}}$ & & & + & & + & & + & & & \\
\hline
$\bar i_2^+ \vphantom{\overline{(}}$ & & & & + & & + & & + & & \\ 
\hline
$\bar i_1^- \vphantom{\overline{(}}$ & & & + & & + & & + & & & \\
\hline
$\bar i_2^- \vphantom{\overline{(}}$ & & & & + & & + & & + & & \\ 
\hline
\end{tabular}
\vspace{2mm}
\nopagebreak

{\bf Boundary map $\partial_3: C_3 \to C_2$, right bottom}
\label{page3to2e}
}
}

{
\tiny

\noindent
\begin{tabular}{| l || c | c | c | c|c|c|c|c|c|c|c|c|c|c|c|c|c|c|c|c|c|c|c|}
\hline
& $U_{12}$ & $U_{13}$ & $U_{21}$ & $U_{23}$ & $U_{31}$ & $U_{32}$ & $V_1$ & $V_2$ & $V_3$ & $W_1$ & $W_2$ & $W_3$ & $X$ & $Y_1$ & $Y_2$ & $Z_1^+$ & $Z_1^-$ & $Z_2^+$ & $Z_2^-$ & $T_1$ & $T_2$ & $J$ & $S$ \\
\hline
$a_{12}$ & & & & + & & & + & & & & & & + & & & & & & & & & & \\
\hline
$a_{13}$ & & + & & + & & + & & & & & & & & & + & & & & & & & & \\
\hline
$a_{14}$ & + & & & & & + & & & & & & & + & & & & & & & & & & \\
\hline
$a_{23}$ & & + & & & + & & & + & & & & & & & & & & & & & & & + \\
\hline
$a_{24}$ & + & & + & & + & & & & & & & & & + & & & & & & & & & \\
\hline
$a_{34}$ & & & + & & & & & & + & & & & + & & & & & & & & & & \\
\hline
$b_{12}$ & + & & & & & & + & & & + & & & & & & & & & & & & & \\
\hline
\end{tabular}\label{page4to3}
\vspace{2mm}

\noindent
{\bf Boundary operator $\partial_4: C_4 \to C_3$, left top}
}

{
\tiny

\noindent
\begin{tabular}{| l || c | c | c | c|c|c|c|c|c|c|c|c|c|c|c|c|c|c|c|c|c|c|}
\hline
& $\bar a_{12} \vphantom{\overline{(}}$ & $\bar a_{13}$ & $\bar a_{14}$ & $\bar a_{23}$ & $\bar a_{24}$ & $\bar a_{34}$ & $\bar b_{12}$ & $\bar b_{13}$ & $\bar b_{14}$ & $\bar b_{21}$ & $\bar b_{23}$ & $\bar b_{24}$ & $\bar b_{31}$ & $\bar b_{32}$ & $\bar b_{34}$ & $\bar b_{41}$ & $\bar b_{42}$ & $\bar b_{43}$ & $\bar c_1$ & $\bar c_2$ & $\bar c_3$ & $\bar c_4$ \\
\hline
$a_{12}$ & + & & + & & & & & & & & & & & & & & & & & & & \\
\hline
$a_{13}$ & & + & & & + & & & & & & & & & & & & & & & & & \\
\hline
$a_{14}$ & & & + & & & + & & & & & & & & & & & & & & & & \\
\hline
$a_{23}$ & + & & & + & & & & & & & & & & & & & & & & & & \\
\hline
$a_{24}$ & & + & & & + & & & & & & & & & & & & & & & & & \\
\hline
$a_{34}$ & & & & + & & + & & & & & & & & & & & & & & & & \\
\hline
$b_{12}$ & & & & & & & + & & & & & & & & & + & & & & & & \\
\hline
\end{tabular}
\vspace{2mm}

\noindent
{\bf Boundary operator $\partial_4: C_4 \to C_3$, center top}
}

{\tiny

\noindent
\begin{tabular}{| l || c | c | c | c|c|c|c|c|c|c|c|c|c|c|c|c|c|c|c|c|c|c|}
\hline
& $\bar a_{12} \vphantom{\overline{(}}$ & $\bar a_{13}$ & $\bar a_{14}$ & $\bar a_{23}$ & $\bar a_{24}$ & $\bar a_{34}$ & $\bar b_{12}$ & $\bar b_{13}$ & $\bar b_{14}$ & $\bar b_{21}$ & $\bar b_{23}$ & $\bar b_{24}$ & $\bar b_{31}$ & $\bar b_{32}$ & $\bar b_{34}$ & $\bar b_{41}$ & $\bar b_{42}$ & $\bar b_{43}$ & $\bar c_1$ & $\bar c_2$ & $\bar c_3$ & $\bar c_4$ \\
\hline
$\bar B_1^1 \vphantom{\overline{(}}$ & + & & + & & & & + & & & & & & & & & & & + & + & & & \\
\hline
$\bar B_1^2 \vphantom{\overline{(}}$ & & & & & & & & + & & & & & & & & & + & & + & & & \\
\hline
$\bar B_1^3 \vphantom{\overline{(}}$ & & + & & & & & & & + & & & & & & & + & & & & & & \\
\hline
$\bar B_2^1 \vphantom{\overline{(}}$ & & & & & + & & + & & & + & & & & & & & & & + & + & & \\
\hline
$\bar B_2^2 \vphantom{\overline{(}}$ & & & & & & & & + & & & & + & & & & & & & + & & & \\
\hline
$\bar B_2^3 \vphantom{\overline{(}}$ & & & + & + & & & & & + & & + & & & & & & & & & + & & \\
\hline
$\bar B_3^1 \vphantom{\overline{(}}$ & + & & & & & + & & & & + & & & & & + & & & & & + & & \\
\hline
\end{tabular}
\vspace{2mm}

\noindent
{\bf Boundary operator $\partial_4: C_4 \to C_3$, center bottom}
}

{
\tiny
\noindent
\begin{tabular}{| l || c | c | c | c|c|c|c|c|c|c|c|c|c|c|c|c|c|c|c|c|c|c|}
\hline
& $\bar d \vphantom{\overline{(}}$ & $\bar e_1^+$ & $\bar e_2^+$ & $\bar e_3^+$ & $\bar e_1^-$ & $\bar e_2^-$ & $\bar e_3^-$ & $\bar k_1^+$ & $\bar k_2^+$ & $\bar k_3^+$ & $\bar k_1^-$ & $\bar k_2^-$ & $\bar k_3^-$ & $\bar g_1$ & $\bar g_2$ & $\bar g_3$ & $\bar h^+$ & $\bar h^-$ & $\bar i_1^+$ & $\bar i_1^-$ & $\bar i_2^+$ & $\bar i_2^-$ \\
\hline
$a_{12}$ & & & & & & & & & & & & & & & & & & & & & & \\
\hline
$a_{13}$ & & & & & & & & & & & & & & & & & & & & & & \\
\hline
$a_{14}$ & & & & & & & & & & & & & & & & & & & & & & \\
\hline
$a_{23}$ & & & & & & & & & & & & & & & & & & & & & & \\
\hline
$a_{24}$ & & & & & & & & & & & & & & & & & & & & & & \\
\hline
$a_{34}$ & & & & & & & & & & & & & & & & & & & & & & \\
\hline
$b_{12}$ & & & & & & & & & & & & & & & & & & & & & & \\
\hline
\end{tabular}
\vspace{2mm}

\noindent
{\bf Boundary operator $\partial_4: C_4 \to C_3$, right top}
}

{\tiny
\noindent
\begin{tabular}{| l || c | c | c | c|c|c|c|c|c|c|c|c|c|c|c|c|c|c|c|c|c|c|}
\hline
& $\bar d \vphantom{\overline{(}}$ & $\bar e_1^+$ & $\bar e_2^+$ & $\bar e_3^+$ & $\bar e_1^-$ & $\bar e_2^-$ & $\bar e_3^-$ & $\bar k_1^+$ & $\bar k_2^+$ & $\bar k_3^+$ & $\bar k_1^-$ & $\bar k_2^-$ & $\bar k_3^-$ & $\bar g_1$ & $\bar g_2$ & $\bar g_3$ & $\bar h^+$ & $\bar h^-$ & $\bar i_1^+$ & $\bar i_1^-$ & $\bar i_2^+$ & $\bar i_2^-$ \\
\hline
$\bar B_1^1 \vphantom{\overline{(}}$ & & & & & & & & & & & & & & & & & & & & & + & \\
\hline
$\bar B_1^2 \vphantom{\overline{(}}$ & & + & & & & & & & & & & & & & & & & & & & & + \\
\hline
$\bar B_1^3 \vphantom{\overline{(}}$ & & & & & + & & & & & & & & & & & + & & & & & & \\
\hline
$\bar B_2^1 \vphantom{\overline{(}}$ & & & & & + & & & & & & & & & + & & & & & & & & \\
\hline
$\bar B_2^2 \vphantom{\overline{(}}$ & & + & & & & & & & & & & & & & & & & & & + & & \\
\hline
$\bar B_2^3 \vphantom{\overline{(}}$ & & & & & & & & & & & & & & & & & & & + & & & \\
\hline
$\bar B_3^1 \vphantom{\overline{(}}$ & & & & & & & & & & & & & & & & & & & + & & & \\
\hline
\end{tabular}
\label{page4to3e}
\vspace{2mm}

\noindent
{\bf Boundary map $\partial_4 : C_4 \to C_3$, right bottom}
}

{
\tiny
\label{page5to4} 

\noindent
\begin{tabular}{| l || c | c | c | c|c|c|c|c|c|c|c|c|c|c|c|c|c|c|c|c|c|}
\hline
& $a_{12}$ & $a_{13}$ & $a_{14}$ & $a_{23}$ & $a_{24}$ & $a_{34}$ & $b_{12}$ & $b_{13}$ & $b_{14}$ & $b_{21}$ & $b_{23}$ & $b_{24}$ & $b_{31}$ & $b_{32}$ & $b_{34}$ & $b_{41}$ & $b_{42}$ & $b_{43}$ & $c_1$ & $c_2$ & $c_3$ \\
\hline
$B_1^1$ & + & & + & & & & + & & & & & & & & & & & & + & & \\
\hline
$B_1^2$ & & & & & & & & + & & & & & & & & & & & + & & \\
\hline
$B_1^3$ & & + & & & & & & & + & & & & & & & & & & & & \\
\hline
$B_2^1$ & & & & & + & & + & & & + & & & & & & & & & & + & \\
\hline
$B_2^2$ & & & & & & & & + & & & & + & & & & & & & & & \\
\hline
$B_2^3$ & & & & + & & & & & + & & + & & & & & & & & & + & \\
\hline
$B_3^1$ & + & & & & & + & & & & + & & & & & + & & & & & & \\
\end{tabular}
\vspace{2mm}
{\bf Boundary operator $C_5 \to C_4$, first page}


\noindent
\begin{tabular}{| l || c | c | c|c|c|c|c|c|c|c|c|c|c|c|c|c|c|c|c|c|c|c|c|}
\hline
 & $c_{4}$ & $d$ & $e_1^+$ & $e_1^-$ & $e_{2}^+$ & $e_2^-$ & $e_3^+$ & $e_3^-$ & $k_1^+$ & $k_1^-$ & $k_2^+$ & $k_2^-$ & $k_3^+$ & $k_3^-$ & $g_1$ & $g_2$ & $g_3$ & $h^+$ & $h^-$ & $i_1^+$ & $i_1^-$ & $i_2^+$ & $i_2^-$ \\
\hline
$B_1^1$ & & & & & & & & & & & & & & & & & & & & & & & \\
\hline
$B_1^2$ & & & + & & & & & & & & & & & & & & & & & & & &\\
\hline
$B_1^3$ & & & & + & & & & & & & & & & & & & & & & & & & \\
\hline
$B_2^1$ & & & & & & & & & & & & & & & + & & & & & & & & \\
\hline
$B_2^2$ & & & & & & & & & & & & & & & & & & & & & + & & \\
\hline
$B_2^3$ & & & & & & & & & & & & & & & & & & & & + & & & \\
\hline
$B_3^1$ & & & & & & & & & & & & & & & & & & & & & & & \\
\hline
\end{tabular}
\vspace{2mm}
{\bf Boundary operator $C_5 \to C_4$, second page}

\noindent
\begin{tabular}{| l || c | c | c | c|c|c|c|c|c|c|c|c|c|c|c|c|c|c|c|c|c|c|}
\hline
$\vphantom{\overline{(}} $ & $\bar B_1^1$ & $\bar B_1^2$ & $\bar B_1^3$ & $\bar B_2^1$ & $\bar B_2^2$ & $\bar B_2^3$ & $\bar B_3^1$ & $\bar B_3^2$ & $\bar B_3^3$ & $\bar B_4^1$ & $\bar B_4^2$ & $\bar B_4^3$ & $\bar B_5^1$ & $\bar B_5^2$ & $\bar B_5^3$ & $\bar E_+$ & $\bar E_1$ & $\bar E_2$ & $ \bar E_3$ & $\bar C_{12}$ & $\bar C_{13}$ & $\bar C_{14}$ \\
\hline
$B_1^1$ & + & & & & & & & & & & & & + & & & & & & & & &\\
\hline
$B_1^2$ & & + & & & & & & & & & & & & + & & & & & & & &\\
\hline
$B_1^3$ & & & + & & & & & & & & & & & & + & & & & & & &\\
\hline
$B_2^1$ & & & + & + & & & & & & & & & & & & & & & & & &\\
\hline
$B_2^2$ & & + & & & + & & & & & & & & & & & & & & & & &\\
\hline
$B_2^3$ & + & & & & & + & & & & & & & & & & & & & & & &\\
\hline
$B_3^1$ & & & & & & + & + & & & & & & & & & & & & & & &\\
\hline
\end{tabular}
\vspace{2mm}
{\bf Boundary operator $C_5 \to C_4$, third page}


\noindent
\begin{tabular}{| l || c | c | c | c|c|c|c|c|c|c|c|c|c|c|c|c|c|c|c|}
\hline
$\vphantom{\overline{(}} $ & $\bar C_{15}$ & $\bar C_{23}$ & $\bar C_{24}$ & $\bar C_{25}$ & $\bar C_{34}$ & $\bar C_{35}$ & $\bar C_{45}$ & $\bar D_{12}^+$ & $\bar D_{12}^-$ & $\bar D_{13}^+$ & $\bar D_{13}^-$ & $\bar D_{14}^+$ & $\bar D_{14}^-$ & $\bar D_{23}^+$ & $\bar D_{23}^-$ & $\bar D_{24}^+$ & $\bar D_{24}^-$ & $\bar D_{34}^+$ & $\bar D_{34}^-$ \\
\hline
$B_1^1$ & & & & & & & & & & & & & & & & & & & \\
\hline
$B_1^2$ & & & & & & & & & & & & & & & & & & & \\
\hline
$B_1^3$ & & & & & & & & & & & & & & & & & & & \\
\hline
$B_2^1$ & & & & & & & & & & & & & & & & & & & \\
\hline
$B_2^2$ & & & & & & & & & & & & & & & & & & & \\
\hline
$B_2^3$ & & & & & & & & & & & & & & & & & & & \\
\hline
$B_3^1$ & & & & & & & & & & & & & & & & & & & \\
\hline
\end{tabular}
}
\label{page5to4e}
\vspace{2mm}
{\bf Boundary operator $C_5 \to C_4$, fourth and last page}

{
\tiny

\noindent
\begin{tabular}{| l || c | c | c | c|c|c|c|c|c|c|c|c|c|c|c|}
\hline
& $A_{11}$ & $A_{12}$ & $A_{13}$ & $A_{21}$ & $A_{22}$ & $A_{23}$ & $A_{31}$ & $A_{32}$ & $A_{33}$ & $A_{41}$ & $A_{42}$ & $A_{43}$ & $A_{51}$ & $A_{52}$ & $A_{53}$ \\
\hline
$B_1^1$ & + & & & & & & & & & & & & & &\\
\hline
$B_1^2$ & & + & & & & & & & & & & & & & \\
\hline
$B_1^3$ & & & + & & & & & & & & & & & &\\
\hline
$B_2^1$ & & & & & & & + & & & & & & & &\\
\hline
$B_2^2$ & & & & & & & & & & + & & & & &\\
\hline
$B_2^3$ & & & & & & & & & & & & & + & &\\
\hline
$B_3^1$ & + & & & & & & & & & & & & & &\\
\hline
\end{tabular}
\label{page6to5}
} \\
{\bf Boundary operator $C_6 \to C_5$}


\begin{thebibliography}{99}

\bibitem{A} V.I.~Arnold, {\it Plane curves, their invariants, perestroikas and classifications.} In: Singularities
and Bifurcations. Providence, RI: AMS, 1994. 33-91. (Adv. in Sov. Math., 21)


\bibitem{borel} A.~Borel, {\it La cohomologie mod 2 de certains espaces homogenes}, Comm. Math. Helvetici {\bf 27} (1953), 165--197.

\bibitem{bl} J.S.~Birman and X.-S.~Lin, {\it Knot polynomials and Vassiliev's invariants}, Invent. Math. 111 (1993) 225-270.

\bibitem{CDM} S.~Chmutov, S.~Duzhin, and J.~Mostovoy, {\it Introduction to Vassiliev knot invariants}, Cambridge University Press, Cambridge, 2012, xvi+504 pp.


\bibitem{FT} R.~Fenn and P.~Taylor, {\it Introducing doodles}. In: Topology of low-dimensional manifolds, R. Fenn (ed.), Lect. Notes Math., v.722, 1977, p. 37-43.

\bibitem{fuks} D.B.~Fuchs, {\it Cohomologies of the braid group mod 2,} Funct. Anal. Appl., 4:2 (1970), 143—151.

\bibitem{GG} M.~Golubitsky, V.~Guillemin; {\it Stable mappings and their
singularities.} Springer, Berlin a.o., 1973.

\bibitem{hir} H.~Hironaka, {\it Resolution of singularities of an algebraic variety over a field of characteristic zero}, Ann. Math. 79 109–326 (1964).


\bibitem{Merx} A.B.~Merkov, {\it Vassiliev Invariants of Doodles, Ornaments, Etc.}
Publications (N.S.) de l'Institut Mathématique, Beograd 66 (80), 1999, 101--126.

\bibitem{M} J.W.~Milnor, {\it Singular points of complex hypersurfaces.} Princeton University Press, 1968.

\bibitem{MS} J.W.~Milnor, J.D.~Stasheff; {\it Characteristic classes.}
Princeton Univ. Press and Univ. of Tokyo Press. Princeton,
New Jersey, 1974.

\bibitem{spi} M.~Spivakovsky, {\it Resolution of Singularities: An Introduction.} In: Cisneros-Molina, J.L., Lê, D.T., Seade, J. (eds) Handbook of Geometry and Topology of Singularities I. Springer, Cham 2020, 183--242.

\bibitem{ks} V.A.~Vassiliev, {\it Cohomology of knot spaces}, in Theory of Singularities and its Applications (Providence) (V. I. Arnold, ed.), Amer. Math. Soc., Providence, 1990, 23--69.



\bibitem{MD} V.A. Vassiliev, {\it On finite order invariants of triple point free plane curves}, Differential topology, infinite-dimensional Lie algebras, and applications, Amer. Math. Soc. Transl. Ser. 2, 194, Amer. Math. Soc., Providence, RI, 1999, 275–300, arxiv: 1407.7227 

\bibitem{pacific} V.A.~Vassiliev, {\it Varieties of chord diagrams, braid group cohomology and degeneration of equality conditions}, Pacific Journal of Mathematics, 326:1 (2023), 135–160, arXiv: 2108.00463

 

\end{thebibliography}
\end{document}